\documentclass[leqno,12pt]{article}
\usepackage{amsmath,amssymb,amscd,amsthm,latexsym,amsfonts, epsfig, graphicx}

\begin{document}

\numberwithin{equation}{section} 

\newtheorem{theorem}{Theorem}[section]
\newtheorem{question}{Question}
\newtheorem{conjecture}[theorem]{Conjecture} 
\newtheorem{lemma}[theorem]{Lemma}
\newtheorem*{claim}{Claim}
\newtheorem{corollary}[theorem]{Corollary}
\newtheorem{proposition}[theorem]{Proposition}

\theoremstyle{remark}
\newtheorem*{remark}{Remark}
\newtheorem*{remarks}{Remarks}
\newtheorem*{merci}{Acknowledgements}
\newtheorem*{defi}{Definition}

\newcommand{\dv}{\operatorname{div}}
\newcommand{\R}{\operatorname{Re}}
\newcommand{\supp}{\operatorname{supp}}
\newcommand{\dist}{\operatorname{dist}}
\newcommand{\Lip}{\operatorname{Lip}}
\newcommand{\diam}{\operatorname{diam}}
\newcommand{\epi}{\operatorname{Epi}}

\newcommand{\NN}{\mathbb{N}}
\newcommand{\RR}{\mathbb{R}}
\newcommand{\ZZ}{\mathbb{Z}}
\newcommand{\QQ}{\mathbb{Q}}
\newcommand{\CC}{\mathbb{C}}
\newcommand{\KK}{\mathbb{K}}
\newcommand{\rn}{\RR^n}
\newcommand{\srn}{{\scriptscriptstyle \RR}^n}
\newcommand{\rplus}{\RR_+}
\newcommand{\rplusbar}{\ol{\RR_+}}

\newcommand{\ep}{\varepsilon}
\newcommand{\si}{\sigma}
\newcommand{\dxdtovert}{\frac{dxdt}{t}}
\newcommand{\dtovert}{\frac{dt}{t}}
\newcommand{\comp}{{}^{\textnormal{c}}}
\newcommand{\lims}{\ol\lim}
\newcommand{\limi}{\ul\lim}
\newcommand{\adelta}{\calA(L)}
\newcommand{\om}{\Omega}
\newcommand{\pom}{\partial\om}
\newcommand{\dive}{\mathrm{div}\,}
\newcommand{\divagrad}{\,\dive\, (A \nabla)}
\newcommand{\kloc}{(K$_{\mathrm{loc}}$)}
\newcommand{\gloc}{(G$_{\mathrm{loc}}$)}
\newcommand{\dloc}{(D$_{\mathrm{loc}}$)}
\newcommand{\CZ}{Calder\'on-Zygmund }

\newcommand{\ul}{\underline}
\newcommand{\ol}{\overline}
\newcommand{\exist}{\;\exists\,}
\newcommand{\foral}{\;\forall\,}

\newcommand{\calA}{\mathcal{A}}
\newcommand{\calB}{\mathcal{B}}
\newcommand{\calC}{\mathcal{C}}
\newcommand{\calD}{\mathcal{D}}
\newcommand{\calE}{\mathcal{E}}
\newcommand{\calF}{\mathcal{F}}
\newcommand{\calG}{\mathcal{G}} 
\newcommand{\calH}{\mathcal{H}} 
\newcommand{\calI}{\mathcal{I}} 
\newcommand{\calJ}{\mathcal{J}} 
\newcommand{\calK}{\mathcal{K}} 
\newcommand{\calL}{\mathcal{L}}
\newcommand{\calM}{\mathcal{M}}
\newcommand{\calN}{\mathcal{N}}
\newcommand{\calO}{\mathcal{O}}
\newcommand{\calP}{\mathcal{P}}
\newcommand{\calQ}{\mathcal{Q}}
\newcommand{\calR}{\mathcal{R}}
\newcommand{\calS}{\mathcal{S}}
\newcommand{\calT}{\mathcal{T}}
\newcommand{\calU}{\mathcal{U}}
\newcommand{\calV}{\mathcal{V}}
\newcommand{\calW}{\mathcal{W}}
\newcommand{\calX}{\mathcal{X}}
\newcommand{\calY}{\mathcal{Y}}
\newcommand{\calZ}{\mathcal{Z}}

\title{
On necessary and sufficient conditions for
  $L^p$-estimates of Riesz transforms
associated to  
elliptic operators on
$\RR^n$ and related estimates}

\author{Pascal Auscher~\footnote{Universit\'e  Paris-Sud and CNRS UMR 
8628, Department de Math\'ematiques, 91405 Orsay Cedex (France) Email: 
\texttt{pascal.auscher@math.u-psud.fr}}}  

\date{January 6, 2004, Revised May 30, 2005}
   
\maketitle

\begin{abstract} This article focuses on   $L^p$ estimates for objects associated to elliptic operators in
divergence form:  its semigroup, the gradient of the semigroup, functional calculus, square functions and 
Riesz transforms.  We introduce four critical numbers associated to the semigroup and its gradient that
completely rule  the ranges of exponents for the $L^p$ estimates. It appears that the case $p<2$
already treated earlier is radically  different from  the case $p>2$ which is new.  We thus recover in a
unified and coherent way many $L^p$ estimates and give further applications. The key tools from harmonic
analysis are two criteria for $L^p$ boundedness, one for $p<2$ and the other for $p>2$ but in ranges
different from the usual intervals $(1,2)$ and $(2,\infty)$.
\end{abstract}

\noindent{\bf Key words:}  elliptic operators, divergence form, semigroup, $L^p$
estimates, Calder\'on-Zygmund theory, good lambda inequalities,  hypercontractivity, Riesz
transforms, holomorphic functional calculus, Littlewood-Paley-Stein estimates.
\bigskip 

\noindent{\bf MSC 2000 numbers:}   42B20, 42B25,  47F05, 47B44, 35J15, 35J30, 35J45.

\vfill\break

\tableofcontents

\vfill\break    

\section*{Introduction} 
 
Although the matter of this text   applies \textit{in extenso} to elliptic operators or systems
in divergence form to any order in $\RR^n$, we focus on second order operators in $\RR^n$.  A section 
will be devoted to these more general classes.   
 
Let $A=A(x)$ be an $n\times n$ matrix of complex, $L^\infty$
coefficients,  defined on $\mathbb{R}^n$, and satisfying the
ellipticity (or ``accretivity'') condition
\begin{equation*}
\lambda |\xi |^2\leq \R  A
\xi \cdot
\overline\xi 
\ \textrm{and}\ | A \xi \cdot \overline \zeta 
| \leq \Lambda |\xi ||\zeta|,\end{equation*} for $\xi,\zeta \in
\mathbb{C}^n$ and for some $\lambda ,\Lambda$ such that $0<\lambda \leq \Lambda
<\infty$. 
 We define a second order divergence form operator
\begin{equation*}
Lf\equiv -\dv (A\nabla f),\end{equation*} which we interpret in the sense of maximal
accretive operators via a sesquilinear form.  Here, $\nabla$ denotes the array of
first order partial derivatives. 

The maximal accretivity condition implies the existence of an analytic contraction semigroup on $L^2(\RR^n)$ generated by $-L$.
It also implies the existence of a holomorphic functional calculus that has the expected stability under commutation and convergence,
allowing for example to define
 fractional powers. This in turn yields the possibility of  defining various objects of interest both from functional and harmonic analysis points of view.
Let us mention Littlewood-Paley-Stein type  functionals  such as 
$$g_{L}(f)(x)= \left( \int_0^\infty |(L^{1/2}e^{-tL}f)(x)|^2 \, {dt}\right)^{1/2}
$$
and 
$$
G_{L}(f)(x)= \left( \int_0^\infty |(\nabla e^{-tL}f)(x)|^2 \, {dt}\right)^{1/2}.
$$
The  ``singular integral'' pending to the latter square function is
 the so-called Riesz transform associated to $L$ given for example by 
$$
\nabla L^{-1/2} = \frac 1 {\sqrt \pi} \int_0^\infty \sqrt t\, \nabla e^{-tL} \,
\frac {dt} t.
$$ Other objects of interest are 1)
the operator of maximal regularity for the parabolic equation associated to $L$, 2)
Riesz means and $L^p$-multipliers when $L$ is, in addition, self-adjoint \dots They
can be treated by the methods presented here but we have chosen not to go into such
developments.  

When the coefficients are constant, \textit{e.g.}  the Laplacian, one finds
classical objects   in harmonic analysis: multipliers, the
Littlewood-Paley-Stein functionals and the original Riesz transforms.  They
belong to the well-understood class of Calder\'on-Zygmund operators. If the
coefficients of $L$ still have some smoothness, then  the  tools of
pseudo-differential calculus  or of Calder\'on-Zygmund theory can still be used.
In absence of regularity of the coefficients, these operators fall beyond the
above classes and this participate to Calder\'on's program~\footnote{\,  See
\cite{C}.} of defining  algebras of differential operators with minimal
smoothness.

The first step of that study is the action on $L^2(\RR^n)$. First, there is a bounded holomorphic functional calculus on $L^2$
basically as a consequence maximal accretivity and Von Neumann's inequality. One has
$$
\|\varphi(L)f\|_2 \le \|\varphi\|_\infty \|f\|_2
$$ for $\varphi$   bounded holomorphic in the open right half plane.   
Next,  $g_L$ and $G_L$ are $L^2$ bounded  (see Section \ref{sec:RTsf}) and
$$
\|g_L(f)\|_2 \sim \|f\|_2 \sim \|G_L(f)\|_2. \ 
\footnote{\, Here $\sim$ is the
equivalence in the sense of norms, with implicit constants 
$C$ depending only on
$n$, $\lambda$ and $\Lambda$.}
$$ The $L^2$ boundedness of the 
Riesz transform has been proved recently and in fact, one has in all dimensions
\begin{equation*}
\| L^{1/2}f\|_{2}\sim
\| \nabla f\|_{2}.
\end{equation*}
This implies that the domain of $L^{1/2}$ is  the Sobolev 
space $W^{1,2}$, which was known as Kato's conjecture.~\footnote{\, The one
dimensional  Kato's conjecture (that is the case $n=1$)) is first proved
by  Coifman, M$^{\rm c}$Intosh \& Meyer
\cite[Th\'eor\`eme X]{CMcM} the two dimensional case by Hofmann \&
M$^{\rm c}$Intosh \cite[Theorem 1.4]{HM} and  the general case in any dimension  by
Hofmann, Lacey, M$^{\rm c}$Intosh \& Tchamitchian along with the author
\cite[Theorem 1.4]{AHLMcT}. We
 refer to the latter reference for historical remarks and connections with other
problems.}

The second step is the action on $L^p(\RR^n)$ with $1<p<\infty $ but $p\ne 2$. 
The bounded holomorphic functional calculus on $L^p$ consists in proving $L^p$ boundedness of $\varphi(L)$ for an appropriate class of  bounded  holomorphic
functions $\varphi$. It is completely understood~\footnote{\, This is
essentially due to Blunck \& Kunstmann \cite[Proposition 2.3]{BK3} although
the authors did not introduce the numbers $p_\pm(L)$.} and one has
$$
\|\varphi(L)(f)\|_p \lesssim \|\varphi\|_\infty \|f\|_p \quad 
\textrm{whenever $p_-(L)<p<p_+(L)$}, ~\footnote{\, Here $\lesssim$ is
the comparison in the sense of norms, with implicit constant 
$C$ that may depend  on $L$ through ellipticity, dimension, its type and  $p_\pm(L)$.}
$$
where $p_-(L)$ and $p_+(L)$ are the two critical exponents for the $L^p$ uniform boundedness of the semigroup 
$(e^{-tL})_{t>0}$. It is clear that this  interval is the largest  open range
of exponents for which such an inequality holds as $\varphi$ can be an exponential.
The
$L^p$ theory for square functions consists in comparing  the $L^p$ norms of
$g_L(f)$, $f$ and  
$G_L(f)$. For $g_L$, what  happens is completely understood in terms of functional
calculus:~
\footnote{\, This follows by combining 
works of Blunck \& Kunstmann \cite[Proposition 2.3]{BK3} and Le Merdy
\cite[Theorem 3]{LeM}.} one has  
$$
\|g_L(f)\|_p \sim \|f\|_p \quad \textrm{whenever $p_-(L)<p<p_+(L)$}.
$$  It turns out that this  interval  is the largest  open
range of exponents for which this equivalence holds. 
 The comparison between
the
$L^p$ norms 
$f$ and  
$G_L(f)$ has not been done in general so far~\footnote{\, After this paper was submitted, the author learned of a work by Yan \cite{Yan} where the inequality $\lesssim$ is obtained for $2n/(n+2) <p\le 2$.} and we shall see that
$$
\|G_L(f)\|_p \sim \|f\|_p \quad \textrm{whenever $q_-(L)<p<q_+(L)$},
$$  where $q_-(L)$ and $q_+(L)$ are the two critical exponents for the $L^p$ uniform 
boundedness  of the gradient of the semigroup 
$(\sqrt t\, \nabla e^{-tL})_{t>0}$, and that this open range is optimal. We also study the corresponding
non-tangential Littlewood-Paley-Stein functionals  (See Section 7).

The $L^p$ theory for   square roots  consists in comparing
$ L^{1/2}f$ and
$ 
\nabla f$ in
$L^p$ norms.~\footnote{\, This program was initialised by the author and P.
Tchamitchian in
\cite{AT}  for this class of complex operators. It arose from a different
perspective towards applications to boundary value problems in the works of
Dalbergh, Jerison, Kenig and their collaborators (see \cite[problem
3.3.16]{Ke}).} There are two issues here, namely the Riesz transform
$L^p$ boundedness, that is an inequality $\| \nabla f\|_{p}\lesssim
\| L^{1/2}f\|_{p}$, and its reverse $\| L^{1/2}f\|_{p}\lesssim
\| \nabla f\|_{p}$. It turns out  the ranges of $p$ may be different.
The state of the art for this class of operators $L$ is as follows.~\footnote{\,
Some of the results were obtained  prior to the Kato conjecture by making the
$L^2$ result an assumption.} One has the Riesz transforms
estimates~\footnote{\,   For
$n=1$, this is due to the author and Tchamitchian
\cite[Th\'eor\`eme A]{AT1} for $n=2$  to the author and Tchamitchian  
\cite[Chapter IV, Theorem 1]{AT} combining the Gaussian estimates of the
author, M$^{\rm c}$Intosh and Tchamitchian \cite[Theorem 3.5]{AMcT2} and the
$L^2$ result 
\cite[Theorem 1.4]{AHLMcT} and for  $n \ge 3$ and
$ 
p_n=\frac{2n}{n+2} < p< 2$,  independently to  Blunck \& Kunstmann  \cite[Theorem 1.2]{BK3},
  and 
to Hofmann \& Martell \cite[Theorem 1.2]{HM}. The enlargement of the range below
$p_n$ is due to the author \cite[Proposition 1]{A} and above 2  is a
consequence of the method  of  the author and Tchamitchian once the reverse
inequality is established (see
\cite[Chapter IV, Proposition 20]{AT} and \cite[Corollary 4]{A}).}
\begin{equation*}
 \| \nabla f\|_{p}\lesssim
\| L^{1/2}f\|_{p} \quad \textrm{if}\quad 
\begin{cases} n=1 &\textrm{and  }\quad 1<p<\infty\\
n=2 &\textrm{and  }\quad 1<p<2+\ep\\
n\ge 3 &\textrm{and  }\quad \frac{2n}{n+2} -\ep < p< 2+\ep
\end{cases},
\end{equation*}
and the reverse inequalities~\footnote{\, They
  are due, for
$n=1$, to the author and Tchamitchian 
(\cite{AT1}, Th\'eor\`eme A), for
$n=2$,  to the author and Tchamitchian
(\cite{AT}, Chapter IV, Theorem 1), again combining  \cite{AMcT2}, Theorem 3.5
with \cite{AHLMcT}, Theorem 1.4,  for
$n=3,4$ to  the author, Hofmann, Lacey, M$^{\rm c}$Intosh \& 
Tchamitchian
(\cite{AHLMcT}, Proposition 6.2) and, for $n\ge 5$ to the author (\cite{A},
Theorem 2).}
\begin{equation*}
 \| L^{1/2}f\|_{p}\lesssim
\| \nabla f\|_{p} \quad \textrm{if}\quad 
\begin{cases} n=1, 2 &\textrm{and  }\quad 1<p<\infty\\
n= 3,4 &\textrm{and  }\quad 1 < p< \frac{2n}{n-2}+\ep\\
n\ge 5 &\textrm{and  }\quad \frac{2n}{n+4} -\ep < p< \frac{2n}{n-2}+\ep
\end{cases}.
\end{equation*}
 Of course, if $L$  possesses more properties then the ranges of exponents  $p$ improve.  For example for
constant coefficients operators these inequalities holds when $1<p< \infty$, and for real operators in dimensions $n\ge 3$, the Riesz transform $L^p$ boundedness is
valid for 
$1<p<2+\ep$ and the reverse inequality for $1<p<\infty$.~\footnote{\, For
constant coefficients, this goes back to Calder\'on-Zygmund original work
\cite{CZ}  and for real operators, this is due to the author and Tchamitchian
\cite[Chapter IV]{AT}.} Hence, it is interesting to have a theory that works for
any single operator.  In fact, the conclusion of the story for the Riesz transform
is
\footnote{\, This problematic  of finding the ``smallest'' exponent $p$ is
implicit in
\cite{BK3}  and we present the counterpart for the ``largest'', which turns out to
require different arguments. After  this paper was submitted, Shen \cite{Sh} informed me of his independent and simultaneous work on the same problem for $p>2$ when $L$ is real. He obtains a different characterization of $q_{+}(L)$ in terms of reverse H\"older estimates for weak solutions of $Lu=0$.  We convince ourselves by e-mail discussions that this approach can be adapted to  complex $L$. }  
$$
\| \nabla f\|_{p}\lesssim
\| L^{1/2}f\|_{p}\quad \textrm{if and only if $ q_-(L)<p<q _+(L)$}
$$
and we also show that
$$
\| L^{1/2}f\|_{p}\lesssim
\| \nabla f\|_{p} \quad \textrm{whenever $\tilde p_-(L)=\sup (1, (p_-(L))_*)<p<p_+(L)$}.
$$
This  encapsulates all the
above mentioned estimates (See Section 1 for the notation $p_*$). Concerning the latter range, we show
that $p_+(L)$ is best possible in some sense, while we only know a lower bound on $\tilde
p_-(L)$.~\footnote{\, As said above, all applies to higher-order operators. For that extended class, this
lower bound is optimal due to some existing appropriate counter-examples in the theory. Similar
counter-examples for the second order case are not known.}  Staring at the formula given above  for
computing the Riesz transform this result seems to say that the integral yields a bounded operator on
$L^p$ if and only if  the integrands are uniformly bounded operators on $L^p$. Said
like this, the sufficiency looks   astonishingly simple. But this is not quite the
truth as there is a play on exponents in the proof. Note also the range of exponents for $L^p$-boundedness
of the Riesz transform $R(L)=\nabla L^{-1/2}$ is characterized. In particular, it is an open set. We also
show that 
$2<p<q_+(L)$ if and only if the Hodge projector   $\nabla L^{-1} \dv$ (or alternately, the second order
Riesz transform $R(L)R(L^*)^*$) is bounded on
$L^p$. For $p<2$ the $L^p$ boundedness of the Hodge projector is stronger than that of $R(L)$. 
   
\bigskip

The objective of this paper is to present a complete, 
coherent and unified theory for all these topics. We present works of others and
also original contributions. In particular, we have tried to be self-contained.  Our main observation
is  the following:
{\bf 
 four critical numbers~\footnote{\, In fact they reduce to two: $p_-(L)$ and $q_+(L)$.} 
rule the
$L^p$ behavior.}  These are $p_\pm(L)$, the limits of the
interval of exponents $p
\in [1,\infty]$ for which  the semigroup  $( e^{-tL})_{t>0}$ is $L^p$ bounded, and  
 $q_\pm(L)$, the limits of the interval of exponents $p \in
[1,\infty]$ for which  the family   $(\sqrt t \, \nabla  e^{-tL})_{t>0}$ is $L^p$
bounded. We make a thorough study of these numbers, their inner relationships and 
their values in terms of dimension for the whole class of such $L$ (Section 4). The
key stone of this theory, which makes use of the divergence structure of our
operators,  is that the $L^p$ boundedness of their semigroups  (resp. of the gradient of their semigroups)
is {\bf equivalent} to some  {\bf off-diagonal estimates} and also to some  {\bf
hypercontractivity}.~\footnote{\, This equivalence for values of $p$ different from
1, although not explicitely stated in the literature to our knowledge, is in the air
of a number of works dealing with semigroups of elliptic operators. It appears first
in Davies's work
\cite{Da}. See also \cite{LSV} and the references therein.} To be precise, we
sometimes have  to loosen  the exponent $p$, but this is harmless for the kind of
results we are after. Let us mention here that this equivalence is not powerful
enough for treating 
$L^\infty$ decay of the semigroup kernels whenever there is $L^1$ or $L^\infty$
boundedness of the semigroup. But, again, this is enough for our needs.    
 
Next, we turn to studying the harmonic analysis objects (Sections 5, 6 \& 7).  
On the one hand, finding necessary conditions on $p$ for which one has $L^p$ bounds
for the functional calculus, the square functions, the Riesz transforms is 
intuitively easy and the critical numbers appear then. On the other hand, it is not
clear at all why these conditions  {\bf alone} suffice. For this, appropriate 
criteria for $L^p$  boundedness with minimal hypotheses are needed.  

The $L^p$ estimates obtained in Sections 5, 6 \& 7 depend on the critical numbers of a given operator.
Thus,  they are individual results with sharp ranges of $p$'s, the operator norms depending on dimension,
the ellipticity constants and some of the  critical numbers. But when the critical numbers can be
estimated   for operators in some class, they become $L^p$ estimates for the whole class. In this case,
the optimality of the range of $p$'s  is function of the sharpness of the bounds on the critical numbers.
This is discussed in Section
\ref{sec:sharpness} for second order operators and in  Section   
\ref{sec:higherorder} for higher order.

As the reader may guess, the various
critical numbers  have no reason to be 1 or $\infty$. Hence, we have a class of
operators that lies {\bf beyond the class of Calder\'on-Zygmund operators}. We wish
here to present an appropriate machinery to obtain $L^p$ boundedness without  caring
about kernels of operators and for ranges of $p$ different from the usual intervals
$(1,2)$ or $(2,\infty)$ (See Section
\ref{sec:CZ} for more).

\section*{Acknowledgements}  It is a pleasure to thank here many colleagues with whom
I discussed these topics over the years and from whom I learnt a lot. Let me mention S. Blunck, T.
Coulhon, G. David, X.T. Duong, S. Hofmann, A. McIntosh,  E. Ouhabaz, P. Tchamitchian. Special thanks go to
E. Ouhabaz who helped me with some historical comments and to S. Blunck for providing me with recent
unpublished work. I also indebted to J.-M.~Martell for a careful reading that helped correcting
some misprints in the submitted version.
\section{Notation}\label{sec:notation}

We constantly work on $\RR^n$, $n\ge 1$, equipped with its usual Lebesgue measure. If $E$ is a measurable set in $\RR^n$, we write
$$
\|f\|_{L^p(E)}  = \left(\int_E |f|^p\right)^{1/p}
$$for the norm in the Lebesgue space $L^p(E)$, $1\le p \le \infty$, with the usual modification if $p=\infty$.
We do not indicate the integration variable and
the measure unless this is necessary for comprehension. Also we drop $E$ in the lower limit if $E=\RR^n$ and set $L^p$ and  $\|f\|_p$ for
$L^p(\RR^n)$ and $\|f\|_{L^p(\RR^n)}$ unless the value of $n$ is of matter.  For Hilbert-valued functions, $|f|$ is replaced by the norm in
the Hilbert space,
$|f|_H$, but we do not introduce a specific notation for $L^p$ as the situation will
make it clear. 

We use the notation $p'$ for the dual exponent to $p$: $p'=\frac p {p-1}$. 

The Sobolev space $W^{m,p}(\RR^n )$, $m\in \NN$, $1\le p\le \infty$, is the space of those  $L^p$ functions $f$ for which all derivatives
up to and including order $m$ are in $L^p$. The norm of $f$ is the sum of the $L^p$ norms of $f$ and all its derivatives. 

The homogeneous Sobolev space  $\dot W^{m,p}(\RR^n )$, $m\in \NN$, $1\le p< \infty$, is the closure of $C_0^\infty(\RR^n)$ for the seminorm being the sum of the
$L^p$ norms of all  derivatives of order $m$. 

We are essentially interested in Sobolev spaces of order $m=1$. The well-known Sobolev inequalities say that 
$$
\dot W^{1,p}(\RR^n) \subset L^q(\RR^n)
$$
whenever, $n>1$,
$$
1\le p<n, \qquad \frac{n}{n-1} \le q <\infty \qquad  \frac n q = \frac n p -1
$$
and 
$$
\|f \|_q \le C(n,p,q)\|\nabla f \|_p. 
$$
We  use the notation $p^*$ ($p$ upper star)  for the Sobolev exponent of $p$, that is 
$$
p^*=\frac {np}{n-p}
$$
with the convention that  $p^*=\infty$ if $p\ge n$ and  
$q_*$ ($q$ lower star) for the reverse Sobolev exponent of $q$, that is $$q_*= \frac
{nq}{n+q}$$

Note that $(p_*)'=(p')^*$ whenever $p_*\ge 1$ and $(p^*)'=(p')_*$ whenever $p^*<\infty$.  

For $1\le p,q\le
\infty$ we set 
$$
\gamma_{pq}=\left | \frac n q - \frac n p \right|
$$
and in the special case where $q=2$
$$
\gamma_p = \left | \frac n 2 - \frac n p \right|.
$$

As usual, we use positive  constants which do not depend on the parameters at stake and whose value  
change at each occurence. Often, we do not mention about such constants as their meaning is self-explanatory.

\section{Beyond Calder\'on-Zygmund operators}\label{sec:CZ}

 By definition,~\footnote{\, We take this   terminology from Meyer \cite[Chapter VII]{Me}.} a Calder\'on-Zygmund operator on $\RR^n$~\footnote{\, or more
generally on a space of homogeneous type
(see \cite{CW}); we shall not be concerned with the development on non homogeneous
spaces described in the excellent review by Verdera \cite{V} although 
 extensions to this setting of the results presented here would be interesting.}
is a bounded operator on
$L^2$ which is associated with a kernel possessing some size and regularity properties, 
the latter being called H\"ormander's condition. 
The fundamental result~\footnote{\, This is due to Calder\'on \& Zygmund 
\cite{CZ}
and H\"ormander \cite[Theorem 2.1]{Ho} 
in the convolution case and is  extended to non
convolution operators in 
\cite{CW}.} is that such operators are 
  weak type $(1,1)$, hence strong type $(p,p)$ when $1<p<2$ by the interpolation theorem of 
Marcinkiewicz and eventually strong type $(p,p)$ for $p>2$ by a duality argument. 
Another route is to begin with their $L^\infty-BMO$
boundedness,~\footnote{\, This is attributed to    Peetre-Spanne-Stein (see
\cite{St2}, p. 191)}  interpolation between $L^2$ and $BMO$
\footnote{\, This result is due to Stampacchia \cite{Sta}.} for $2<p<\infty$, and
duality for $1<p<2$. However, one should not forget that this interpolation is more
involved than the Calder\'on-Zygmund d\'ecomposition used for weak type $(1,1)$.

In applications, this is enough for numerous operators going from convolution operators such as the Hilbert transform and  the classical Riesz transforms (the
prototypes of Calder\'on-Zygmund operators)  to the Cauchy integral on a Lipschitz curve and the double layer operator on a Lipschitz
domain.

However, recently some interesting operators were found to be out of this class. That is they are strong type $(2,2)$ but the other 
properties fail. Some reasons  are 
\begin{enumerate}
\item their  kernel does not possess regularity properties such as the H\"ormander condition
\item they do not possess a kernel in any
reasonable sense (but the distribution sense)
\item they are found to be strong type $(p,p)$  for a range of values of $p$ different from $1<p<2$ or $2<p<\infty$ or their unions
\item  duality does not apply
\end{enumerate}

It is natural to ask the following question: {\bf is there a general machinery to
handle the $L^p$ theory of  such operators? }

The answer is yes. It turns out that the cases $p<2$ and $p>2$ are treated by different methods, which is useful when duality
is not available.

Let us come now to statements. For simplicity, we work in the framework of $\RR^n$ equipped with the Lebesgue measure, although  the original 
results are presented in spaces of homogeneous types.

We  denote as above by  $B(x,r)$ the open (Euclidean) ball of radius $r>0$ and center
$x\in \RR^n$,  and set $|E|$ the measure of a set $E$.

Here is some further notation used throughout the paper. For a ball $B$, we let $\lambda B$ be the ball with same center and radius $\lambda$ times that of $B$. 
We  set 
$$C_1(B)=4B  \quad \mathrm{and}\quad C_j(B)=2^{j+1}B \setminus 2^{j}B, \ \mathrm{ if}\ j\ge 2.$$  
We adopt the similar  notation, $\lambda Q$ and $C_j(Q)$, for  any cube $Q$
with sides parallel to the axes.

Denote by  ${ M}$  the 
Hardy-Littlewood maximal operator
$${ M}f(x)=\sup_{B\ni x}\frac{1}{|B|}\int_B|f|,$$
where $B$ ranges over all open balls (or cubes)  containing $x$.

\begin{theorem}\label{lpp<2}~\footnote{\, This is due to Blunck \& Kunstmann
\cite[Theorem 1.1]{BK2} generalizing earlier work of  Duong \& M$^{\rm c}$Intosh
\cite[Theorem 1]{DMc} who obtained weak type $(1,1)$ under a weakened
H\"ormander condition, still assuming  reasonable  pointwise estimates for kernels
but no regularity in the classical sense. The statement and  proof here
simplify that of 
\cite{BK2}. Note that  when $p_0=1$ the assumptions are slightly different
 than those in \cite{DMc}. One can find in Fefferman \cite{Fef} weak type results for values of $p$ not
equal to one but no general statement is made.}  
  Let $p_0\in [1,2)$.
Suppose
that $T$ is sublinear operator  of strong type $(2,2)$, and let  $A_r$, $r>0$,
be  a family of linear operators acting on $L^2$.
     Assume for $j\ge 2$   
\begin{equation} \label{pointwisep<2}
\left(\frac{1}{|2^{j+1}B|}\int_{C_j(B)}|T(I-A_{r(B)})f|^2\right)^{1/2}  \le  g(j) \left(\frac{1}{|B|}\int_B |f|^{p_0}\right)^{1/p_0}
\end{equation}
and for $j\ge 1$
\begin{equation} \label{domap<2}
\left(\frac{1}{|2^{j+1}B|}\int_{C_j(B)} |A_{r(B)}f|^{2}\right)^{1/2}  \le  g(j) \left(\frac{1}{|B|}\int_B |f|^{p_0}\right)^{1/p_0}
\end{equation}
for   all ball $B$ with $r(B)$ the radius of $B$ and all $f$ supported in $B$.
If $\Sigma= \sum g(j) 2^{nj} <\infty$, then  $T$ is of weak type  $(p_0,p_0)$, with a bound depending only on the strong type $(2,2)$ bound of $T$,   $p_0$ and
$\Sigma$, hence bounded on
$L^p$ for $p_0<p<2$.
\end{theorem}

\begin{theorem}\label{lpp>2}\footnote{\, This is due to  the
author, Coulhon, Duong \& Hofmann
\cite[Theorem 2.1]{ACDH} using ideas of
  Martell in \cite{Ma} who developed a variant of the sharp function theory
of Fefferman-Stein \cite{FS} in the spirit of \cite{DMc}, again assuming there
are reasonable pointwise  estimates for kernels but no regularity in the classical
sense. The proof here simplifies the exposition in \cite{ACDH}.  Shen  independently proved a  boundedness result similar in spirit by analogous methods \cite[Theorem 3.1]{Sh} which he attributes to Caffarelli and Peral \cite{CP}. In fact, it is easy to recover Shen's theorem as a consequence of this one.}  Let 
$p_0\in (2,\infty]$. Suppose
that $T$ is sublinear operator acting on $L^2$, and let  $A_r$, $r>0$,
be  a family of linear operators acting on $L^2$.
     Assume   
\begin{equation}
\left(\frac{1}{|B|}\int_B|T(I-A_{r(B)})f|^2 \right)^{1/2}  \le  C\big({M}
(|f|^2)\big)^{1/2}(y) ,\label{pointwise}
\end{equation}
and 
\begin{equation}
\left(\frac{1}{|B|}\int_B |TA_{r(B)}f|^{p_0}\right)^{1/p_0}  \le  C\big({M}
(|Tf|^2)\big)^{1/2} (y),\label{doma}
\end{equation}
for all $f\in L^2$, all ball $B$ and all $y \in B$ where $r(B)$ is the radius of $B$.
If  $2<p<p_0$ and  $Tf \in L^p$ when $f\in L^p$  then $T$ is strong type $(p,p)$.
More precisely,  for all $f\in L^p\cap L^2$,
$$
\|Tf\|_p \le c \|f\|_p
$$
where $c$ depends only on $n$, $p$ and $p_0$
and $C$.
\end{theorem}

\begin{remarks}
\begin{enumerate}

\item The operators $A_r$ play the role of approximate identities (as $r\to 0$) eventhough $A_r(1)=1$ is not assumed.
The boundedness of $A_r$ on $L^2$ is a consequence of linearity but the $L^2$ bounds are not explicitely needed.
In applications, the $L^2$ bounds are uniform in $r$ and used to check the
hypotheses. The improvement in the exponents from $p_0$ to $2$ in \eqref{domap<2}
and  from 2 to $p_0$ in \eqref{doma}  expresses the regularizing effect of $A_r$.
When $p_0=\infty$, the left hand side of \eqref{doma} is understood as the 
essential supremum on $B$.

\item Possible weakening of Theorem \ref{lpp<2} is as follows: the exponent 2 in
\eqref{pointwisep<2} can be changed to $1$ and the exponent $2$ in \eqref{domap<2}
can be changed to $p>p_0$.

\item As we shall see, Theorem \ref{lpp>2} has little to do with operators but
rather with decomposition of functions in the spirit of Fefferman-Stein's argument 
for the sharp function and, in fact, it is an extension of it. This is why 
the
regularised version $TA_r$ of
$T$ is controlled by the maximal function of $|Tf|^2$,
 which may be surprising at first
sight.

\item Define, for $f\in
L^2$, 
$${\cal M}^{\#}_{T,A}f(x)=\sup_{B\ni
x}\left(\frac{1}{|B|}\int_B|T(I-A_{r(B)})f|^2\right)^{1/2},$$ where the supremum is taken
over all balls $B$ in $\RR^n$ containing $x$,  and $r(B)$ is  the radius of $B$.  The assumption is that ${\cal M}^{\#}_{T,A}f$ is
controlled pointwise by $({ M}
(|f|^2))^{1/2}$. In fact, rather than the exact form of the control what matters is that ${\cal M}^{\#}_{T,A}$ is
strong type
$(p,p)$ for the desired values of $p$.

\item The family $(A_r)$ indexed by positive $r$ could be replaced by a family $(A_B)$ indexed by
balls. Then $A_{r(B)}$ is replaced by $A_B$ in the statements. An example of such a family  is given by
mean values operators 
$A_Bf=\frac{1}{|B|} \int_B f$.
 
\item Note that in Theorem \ref{lpp>2}, $T$  acts on $L^2$ but its 
boundedness is not needed in the proof. However, it is used in applications to check 
\eqref{doma} and \eqref{pointwise}. Note also that $T$  already acts on
$L^p$ and the purpose of the statement is to bound its norm. 
In concrete situations, this theorem
 is applied to suitable approximations of
$T$, the uniformity of the bounds allowing a 
limiting argument to deduce $L^p$ boundedness. 
So an argument to conclude for generic $L^p$ functions is not needed here.

\item Both theorems are valid in the vector-valued context, that is when  $f$ is valued in a Banach space $B_1$ and $Tf$ is valued in 
another Banach space $B_2$. We leave to the reader the care of checking details. We apply this for square function estimates in Section \ref{sec:RTsf}.  

\end{enumerate}

\end{remarks}

\paragraph{Proof of Theorem \ref{lpp<2}:}
 It begins with the classical Calder\'on-Zygmund decomposition which we
recall.~\footnote{\, Many good references for this tool. One is \cite{St1}.}

\begin{lemma}\label{lemmaCZDnormal} Let $n\ge 1$,  $1\le p\le \infty$ and  
$\| f\|_p <\infty$. Let $\alpha>0$. Then, one can find a collection of cubes $(Q_i)$, functions $g$ and $b_i$  such that 
\begin{equation}\label{eqczdn1}
f= g+\sum_i b_i  \end{equation}
and the following properties hold:
\begin{equation}\label{eqczdn2}
\| g\|_\infty \le C\alpha, \end{equation}
\begin{equation}\label{eqczdn3}
\supp b_i \subset Q_i \ \text{and} \ \int_{Q_i} | b_i|^p \le C\alpha^p |Q_i|, \end{equation}
 \begin{equation}\label{eqczdn4}
\sum_i |Q_i| \le C\alpha^{-p} \int_{\RR^n} | f|^p , \end{equation}
\begin{equation}\label{eqczdn5}
\sum_i {\bf 1}_{Q_i} \le N, \end{equation}
where $C$ and
 $N$ depends only on dimension and $p$.
\end{lemma}

Let $f \in L^{p_0} \cap L^2$. We have  to  prove that  
 for any $\alpha>0$, 
$$
\bigg|\bigg\{x\in \RR^n; |Tf(x)| >  \alpha \bigg\}\bigg| \le  \frac C{\alpha^{p_0}}\int_{\RR^n}
|f|^{p_0}.
$$
Write $f=g+\sum_i b_i$ by the Calder\'on-Zygmund decomposition at height $\alpha>0$. The construction of this decomposition implies that 
$g \in L^2$ with $\int |g|^2 \le C \alpha^{2-p_0} \int |f|^{p_0}$. Then $Tg \in L^2$ with 
$\|Tg\|_2 \le C\|g\|_2$ by the assumption in Theorem \ref{lpp<2}. This with \eqref{eqczdn2} yield
$$
\bigg|\bigg\{x\in \RR^n; |Tg(x)| > \frac \alpha 3\bigg\}\bigg| \le \frac C{\alpha^2} \int_{\RR^n} |g|^2 \le \frac C{\alpha^{p_0}}\int_{\RR^n}
|f|^{p_0}.
$$
To handle the remaining term,  introduce for $r\ge 0$  the operator
$B_r=I-A_r$  and let $r_i$ be the radius of  $Q_i$. 
Since $$
\bigg|T\bigg(\sum_i b_i\bigg)\bigg| \le \sum_i |TB_{ri}b_i| + \bigg|T\bigg(\sum_i
A_{r_i}b_i\bigg)\bigg|
$$
 it is enough to 
estimate $A=|\{x \in \RR^n; \sum |TB_{r_i}b_i(x)| >\alpha/3\}|$ and 
$B=|\{x \in \RR^n; |T(\sum_i A_{r_i}b_i)(x)| >\alpha/3\}|$. 
Let us  bound the first term.  
First,  
$$A \le |\cup_i 4Q_i| + \bigg|\bigg\{x \in \RR^n \setminus \cup_i  4Q_i ; 
\sum_i |TB_{r_i}b_i(x)| >\frac \alpha 3\bigg\}\bigg|,$$
and by \eqref{eqczdn4}, $|\cup_i 4Q_i| \le \frac{C}{\alpha^{p_0}}\int | f|^{p_0}$.
To handle the other term, we observe that by Tchebytchev inequality, 
$$
\bigg|\bigg\{x \in \RR^n \setminus \cup_i 4Q_i ; \sum_i
|TB_{r_i}b_i(x)| >\frac \alpha 3\bigg\}\bigg|
\le
\frac{C}{\alpha^2}\int \bigg|\sum_i h_i\bigg|^2$$ with $h_i = {\bf
1}_{(4Q_i)^c}|TB_{r_i}b_i|$. To estimate the $L^2$ norm, we dualize against
$u\in L^{2}$ with $\|u\|_{2}=1$. Write 
$$
\int  |u| \sum_i h_i = \sum_i\sum_{j=2}^\infty A_{ij}
$$  
where 
$$
A_{ij}= \int_{C_j(Q_i)} |TB_{r_i}b_i||u|.
$$  
By  \eqref{pointwisep<2} and \eqref{eqczdn3},
\begin{align*}\| TB_{r_i}b_i \|_{L^2(C_j(Q_i))}  
&\le  |2^{j+1}Q_i|^{1/2} g(j) \left(\frac{1}{|Q_i|}\int_{Q_i} |b_i|^{p_0}\right)^{1/p_0}
\\
& \le C |2^{j+1}Q_i|^{1/2} g(j) \alpha
\end{align*}
for some appropriate constant $C$.
Now remark that for any $y \in Q_i$ and any $j\ge 2$, 
$$
\left( \int_{C_j(Q_i)} |u|^{2}\right)^{1/2} \le \left( \int_{2^{j+1}Q_i} |u|^{2}\right)^{1/2} \le  |2^{j+1}Q_i|^{1/2}
\big(M(|u|^{2})(y)\big)^{1/2}.
$$
Applying H\"older inequality,     one obtains
$$
A_{ij} \le  C\alpha 2^{nj} g(j) |Q_i| \big(M(|u|^{2})(y)\big)^{1/2}.
$$
Averaging over $Q_i$ yields
$$
A_{ij} \le C\alpha 2^{nj} g(j) \int_{Q_i} \big(M(|u|^{2})(y)\big)^{1/2}\, dy.$$
Summing over $j\ge 2$ and $i$, we have
\begin{align*}
\int| u| \sum_i h_i &\le C \alpha \int \sum_i {\bf 1}_{Q_i}(y) \big(M(|u|^{2})(y)\big)^{1/2}\, dy
\\
& \le C N\alpha \int_{ \cup_i Q_i} \big(M(|u|^{2})(y)\big)^{1/2}\, dy
\\
& \le C'N\alpha \big| \cup_i Q_i \big|^{1/2} \||u|^2\|_1^{1/2}.
\end{align*}
In the next to last inequality, we used \eqref{eqczdn5}, and in the  last 
inequality, we used Kolmogorov's lemma and the weak type $(1,1)$ of the maximal
function.~\footnote{\, This idea is borrowed from Hofmann \& Martell 
\cite{HM}.} Hence from
\eqref{eqczdn4} 
$$
A \le  \frac {C'' N }{ \alpha^{p_0}} \int_{\RR^n} | f|^{p_0}.
$$

It remains to estimate the term $B$. To this end, we use that $T$ is bounded on $L^2$ to obtain
$$
B
\le
\frac{C}{\alpha^2}\int \bigg|T\bigg(\sum_i  A_{r_i}b_i\bigg)\bigg|^2
\le \frac{C}{\alpha^2}\int \bigg|\sum_i  A_{r_i}b_i \bigg|^2. $$
To estimate the $L^2$ norm, we dualize against
$u\in L^{2}$ with $\|u\|_{2}=1$ and write
$$ 
\int| u| \sum_i |A_{r_i}b_i| = \sum_i\sum_{j=1}^\infty B_{ij}
$$  
where 
$$
B_{ij}= \int_{C_j(Q_i)} |A_{r_i}b_i||u|.
$$ 
Using \eqref{domap<2} 	and  \eqref{eqczdn3},  
\begin{align*}\| A_{r_i}b_i \|_{L^2(F_{ij})}  
&\le  |2^{j+1}Q_i|^{1/2} g(j) \left(\frac{1}{|Q_i|}\int_{Q_i} |b_i|^{p_0}\right)^{1/p_0}
\\
& \le C |2^{j+1}Q_i|^{1/2} g(j) \alpha
\end{align*}
 for $j \ge 1$. 
From here, we may argue as before and conclude that $B$ is bounded by 
$\frac {C }{ \alpha^{p_0}} \int| f|^{p_0}$ as desired.  
\qed

\begin{remark} As the reader can see, there is absolutely no use of mean value
properties of the $b_i$'s. 
\end{remark}

\paragraph{Proof of Theorem \ref{lpp>2}:}

 We begin with a useful localisation lemma.

\begin{lemma} There is $K_0$ depending only on dimension such that the
following holds. For every  
$f\in L^1_{loc}$  and  every cube 
$Q$ and every $\lambda>0$ for which there exists $\bar x \in 4Q$ for which
$Mf(\bar x)
\le
\lambda$, then for every $K\ge K_0$,
$$
\{\chi_QMf >K\lambda\} \subset  \{M(f\chi_{8Q}) >\frac K {K_0}\lambda\}
$$ 
\end{lemma}

\paragraph{Proof:}  We use that $M$ is comparable to the
centered maximal function
$M_c$: there is $K_0$ depending only on the doubling constant such that
$M \le K_0 M_c$. 

Let $x \in Q$ with $Mf(x)>K\lambda$. Then  $M_cf(x) > \frac K {K_0}
\lambda$. Hence, there is a cube  centered at $x$ with radius $r$ such that
$$
\frac {1}{|Q(x,r)|}\int_{Q(x,r)} |f| > \frac K {K_0}
\lambda.
$$
If $\frac K {K_0}
>1$,  $\bar x \notin Q(x,r)$ since $Mf(\bar x) \le \lambda$. 
The conditions $x\in Q$, $\bar x \in 4Q$ and $\bar x \notin Q(x,r)$ imply
$Q(x,r) \subset 8Q$. Hence,
$$
 \frac K {K_0}
\lambda < \frac {1}{|Q(x,r)|}\int_{Q(x,r)} |f\chi_{8Q}| \le
M(f\chi_{8Q})(x).
$$  
\qed

We continue with a two parameters family of good lambda inequalities.

\begin{proposition}\label{goodlambda} Fix $1<q\le \infty$ and $a>1$. Then, there exist $C=C(q,
n, a)$ and $K'_0=K_0'(n,a)$ with the following property: 
If $F,G$ are  nonnegative measurable functions such that for every
cube  $Q$ there exist non negative functions $G_Q,H_Q $ with 
$$
F\le G_Q+H_Q \quad a.e. \ on \ Q,
$$
$$
 \left( \frac {1}{|Q|}\int_{Q}   {H_Q}^{q}\,
\right)^{1/q} \le a MF(x) \quad {\mathrm for\ all}\ x \in Q
$$
$$
 \frac {1}{|Q|}\int_{Q} G_Q\,
 \le G(x)\quad {\mathrm for\ all}\ x \in Q.
$$
Then for all $\lambda>0$, for all $K>K_0'$ and $\gamma<1$, 
$$
|\{MF>K\lambda, G\le \gamma\lambda\}| \le C\left(\frac{1}{K^q} + \frac
\gamma K\right) |\{MF>\lambda\}|.
$$
If $q=\infty$, we understand the average in $L^q$ as an essential supremum. In this
case, $\frac{1}{K^q} =0$.

\end{proposition}

\paragraph{Proof:}   Let
$E_\lambda=\{MF>\lambda\}$. We assume this is a proper  subset in $\RR^n$ otherwise there is nothing to
prove. Since $E_\lambda$ is open, the Whitney  decomposition  
yields a family of non overlapping dyadic cubes $Q_i$ 
such that $E_\lambda=\cup_iQ_i$ and $4Q_i$ contains at least one point
$\overline{x_i}$ outside
$E_\lambda$, that is
$${ M}F(\overline{x_i})\le \lambda.$$

Let 
$B_\lambda=\{MF>K\lambda, G\le \gamma\lambda\}$. If $K\ge 1$ then 
$B_\lambda \subset E_\lambda$, hence
$$
|B_\lambda| \le \sum_i |B_\lambda \cap Q_i| .
$$
Fix $i$. If $B_\lambda
\cap Q_i=\emptyset$, we have nothing to do. If not, there is a point
$\overline{y_i} \in Q_i$ such that 
$$G(\overline{y_i}) \le
\gamma\lambda.$$
By the localisation lemma applied to $F$ on $Q_i$, if $K\ge K_0$, then
$$
|B_\lambda
\cap Q_i|\le |\{MF>K\lambda\} \cap Q_i| \le
|\{M(F\chi_{8Q_i})>
\frac K {K_0} \lambda\}|.
$$
Now use $F\le G_i+H_i$ with $G_i=G_{8Q_i}$ and $H_i=H_{8Q_i}$ to deduce
$$|\{M(F\chi_{8Q_i})>
\frac K {K_0} \lambda\}| \le 
|\{M(G_i\chi_{8Q_i})>
\frac K {2K_0} \lambda\}| + |\{M(H_i\chi_{8Q_i})>
\frac K {2K_0} \lambda\}|.
$$
Now by using the weak type $(1,1)$ and $(q,q)$ of the maximal operator
with respective constant $c_1$ and $c_q$, we have
\begin{align*}
|\{M(G_i\chi_{8Q_i})>
\frac K {2K_0} \lambda\}| &\le \frac {2K_0c_1}{K\lambda} \int_{8Q_i} G_i
\\
& \le \frac
{2K_0c_1}{K\lambda} |8Q_i| G(\overline{y_i})
\\
& \le 
\frac {2K_0c_1}{K} |8Q_i| \gamma,
\end{align*}
and, if $q<\infty$
\begin{align*}
|\{M(H_i\chi_{8Q_i})>
\frac K {2K_0} \lambda\}| &\le \left(\frac {2K_0c_q}{K\lambda}\right)^q
\int_{8Q_i} H_i^q
\\
&
\le \left(\frac {2K_0c_q}{K\lambda}\right)^q
|8Q_i| (a MF(\overline{x_i}))^q \\
&
\le 
\left(\frac {2K_0c_q}{K}\right)^q |8Q_i| a^q.
\end{align*}
Hence, summing over $i$ yields
$$|B_\lambda|
\le 
\sum_i C\left(\frac{a^q}{K^q} + \frac
\gamma K\right) |8Q_i| \le C8^n\left(\frac{a^q}{K^q} + \frac
\gamma K\right) |E_\lambda|. 
$$
If $q=\infty$, then 
$$
\|M(H_i\chi_{8Q_i})\|_\infty \le  \|H_i\chi_{8Q_i}\|_\infty  \le a
MF(\overline{x_i}) \le a \lambda,
$$
so that, choosing $K\ge 2K_0a$ leads us to $\{M(H_i\chi_{8Q_i})>
\frac K {2K_0} \lambda\}=\emptyset$.

\qed

Let us now assume that $F,G$ are so that the conclusion of the proposition 
holds. Let $0<p<q$. Then, we have
$$
\|MF\|_p^p \le CK^p\left(\frac{a^q}{K^q} + \frac
\gamma K\right) \|MF\|_p^p + \frac {K^p} {\gamma^p} \|G\|_p^p.
$$
Hence, if furthermore $\|MF\|_p<\infty$, since $p<q$, one can choose $K$
large enough and $\gamma$ small enough so that 
$$
CK^p\left(\frac{a^q}{K^q} + \frac
\gamma K\right) \le 1 - 2^{-p}.
$$
This choice depends therefore on $p$, $q$, $n$, $a$ and we have
$$
\|MF\|_p^p \le \frac {(2K)^p} {\gamma^p} \|G\|_p^p.
$$

Let us now prove Theorem \ref{lpp>2}.
Let $f\in
L^2\cap L^p$. We let $q=\frac{p_0} 2$ and set $F=|Tf|^2\in L^{p/2}$. 
By sublinearity of $T$, we have for any cube  $Q$ that  
$F\le G_Q + H_Q$ with $G_Q=2 |T(I-A_{r(Q)})f|^2$ and
 $H_Q= 2|TA_{r(Q)}f|^2
$. Hence the hypotheses of the proposition apply with $a=2C^2$ and 
$G=2C^2{ M}(f^2)$. Let $2<p<p_0$. Since we know that $MF \in
L^{p/2}$ from the hypothesis, we obtain
$$
\|Tf\|_p^2 \le \|MF\|_{p/2} \le C \|G\|_{p/2}  \le C'
\|f\|_p^2.
$$
In the last inequality, we have used  the assumption on $T(I-A_r)$.
\qed

\begin{remark} 
With $q=\infty$  and given $F$, we let $G_Q=|F-m_QF|$ and $H_Q=|m_QF|$,
$a=1$ and $G=M^\#F$, the sharp function of  Fefferman-Stein. We obtain if $0<p<\infty$
and
$\|MF\|_p <\infty$ that 
$$
\|MF\|_p \le C_p \|M^\#F\|_p,
$$
hence this argument contains in particular Fefferman-Stein's. To recover  the sharp
function introduced by Martell, we take $G_Q=|F-A_{r(Q)}F|$ and $H_Q=|A_{r(Q)}F|$
and, as proved Martell,   $\sup_{Q\ni x} H_Q \le a MF(x)$ under kernel upper bounds
on $A_r$.

\end{remark}

\section{Basic $L^2$ theory for elliptic operators}\label{sec:basicL2}

\subsection{Definition}\label{sec:def}

Let $A=A(x)$ be an $n\times n$ matrix of complex, $L^\infty$
coefficients,  defined on $\mathbb{R}^n$, and satisfying the
ellipticity (or ``accretivity'') condition
\begin{equation}\label{eq1} \lambda |\xi |^2\leq \R  A
\xi \cdot
\overline\xi 
\ \textrm{and}\ | A \xi \cdot \overline \zeta 
| \leq \Lambda |\xi ||\zeta|,\end{equation} for $\xi,\zeta \in
\mathbb{C}^n$ and for some $\lambda ,\Lambda$ such that $0<\lambda \leq \Lambda
<\infty$. 
 We define a second order divergence form operator
\begin{equation}\label{eq2}
Lf\equiv -\dv (A\nabla f),\end{equation} which we first interpret in the sense of maximal accretive operators via a 
sesquilinear form.
 That is,  $\calD(L)$ is the largest subspace  contained in $W^{1,2}$ for which
$$
\left|\int_{\RR^n} A \nabla f \cdot \nabla g\right| \le C \|g\|_2
$$
for all $g \in W^{1,2}$ and we set  $Lf$ by 
$$
\langle Lf, g\rangle = \int_{\RR^n} A \nabla f \cdot \overline{\nabla g}
$$
for $f \in \calD(L)$ and $g\in W^{1,2}$. 
Thus defined, $L$ is maximal-accretive operator on $L^2$ 
and $\calD(L)$ is dense in $W^{1,2}$.~\footnote{\, For precise
definitions, details and proofs, see Kato's book \cite[Chapter VI]{K2}.} 

Since $W^{1,2}$ is
dense in its homogeneous version $\dot W^{1,2}$ (for the semi-norm $\|\nabla
f\|_2$),  
 $L$ extends to a bounded operator invertible operator from $\dot W^{1,2}$ into its dual space
$\dot W^{-1,2}$, which justifies the divergence notation in \eqref{eq2}. In particular, one has
\begin{equation}\label{hodgep=2}
\|\nabla	L^{-1}\dv f\|_2 \le c\|f\|_2.
\end{equation}

\subsection{Holomorphic functional calculus on $L^2$}\label{sec:fcL2}

Let $L$ be as above. There exists $\omega\in [0,{\frac \pi 2})$ 
depending only on the ellipticity constants  such that
for all $f \in {\cal D}(L)$,
\begin{equation}\label{type}
| \arg \langle Lf, f \rangle| \le \omega.
\end{equation}
We fix the smallest such $\omega$ 
and the following $L^2$ estimate is easily proved: for all $\mu \in (\omega, \pi)$
and all  complex numbers $\lambda \in \Sigma_{\pi-\mu}$,
$$
\|(L+\lambda)^{-1}f\|_2 \le \frac{c\|f\|_2}{|\lambda|},
$$
where we have set $\Sigma_\mu=\{z\in \CC^*; |\arg z\, | <\mu\}$.~\footnote{\,
See, for example, \cite[Preliminary Chapter]{AT} for a proof.} 
Hence,  $L$ is  of type $\omega$ on
$L^2$. Several consequences follow.


In particular,  $-L$ generates a semigroup $(e^{-tL})_{t>0}$
which has an  analytic extension to a complex semigroup $(e^{-zL})_{z \in \Sigma_{{\frac \pi 2} -\omega}}$ of contractions on
$L^2$.

Since $L$ is also maximal-accretive operator,
it has a bounded holomorphic functional calculus on $L^2$. 
In particular, 
for any $\mu\in (\omega,\pi)$ and any $\varphi$ holomorphic and bounded in $\Sigma_\mu$, the operator $\varphi(L)$ is
bounded on
$L^2$ with the estimate
$$
\|\varphi(L)f\|_2 \le  c \|\varphi\|_\infty \|f\|_2,
$$
the constant $c$ depending only on $\omega $ and $\mu$. 
If $\varphi$ satisfies the technical condition 
\begin{equation}\label{eqtechnical}
|\varphi(\zeta)| \le c |\zeta|^s(1+|\zeta|)^{-2s}
\end{equation}
for all $\zeta \in \Sigma_\mu$ for some positive constants $c$, $s$, then $\varphi(L)$ can be computed 
 using the semigroup. Let
$\omega<\theta<\nu<\mu<\frac \pi 2$. One has
\begin{equation}\label{eqrepresentation}
\varphi(L)= \int_{\Gamma_+} e^{-zL} \eta_+(z)\, dz + \int_{\Gamma_-} e^{-zL} \eta_-(z)\, dz 
\end{equation}
where $\Gamma_\pm$ is the half-ray $\RR^+ e^{\pm i (\frac \pi 2-\theta)}$, 
\begin{equation}\label{eqeta}
\eta_\pm(z)= 
\frac 1 {2\pi i} \int_{\gamma_\pm} e^{\zeta z} \varphi(\zeta)\, d\zeta, \quad z
\in \Gamma_\pm,
\end{equation}
with $\gamma_\pm$  being the half-ray $ \RR^+ e^{\pm i \nu}$ (the orientation on paths is irrelevant in the arguments where this 
representation is used so we do not
insist on that).
For general bounded holomorphic functions, $\varphi(L)$ is defined by a limiting procedure which we do not 
need in this work.

Finally, one can define  unbounded operators $\varphi(L)$ for $\varphi$ holomorphic in $\Sigma_\mu$, $\omega<\mu<\pi$, satisfying
\begin{equation}\label{eqtechnicals}
|\varphi(\zeta)| \le c \sup( |\zeta|^s,|\zeta|^{-s'})
\end{equation}
for some $c,s,s'\ge 0$. This includes fractional powers of $L$.  
We call $\calF(\Sigma_\mu)$ the class of such holomorphic functions.~\footnote{\,
For definitions, we refer to \cite{Mc} and \cite{CDMcY}.}

\subsection{$L^2$ off-diagonal estimates}\label{sec:offdiagL2}

A very important ingredient in this paper is the off-diagonal estimates of Gaffney type. They are crucial to 
our analysis because when dealing with complex operators, we do not have at our disposal contractivity of
the semigroup on  all $L^p$ spaces, and in fact this is false in general.

\begin{defi} Let $\calT=(T_t)_{t>0}$ be a family of operators. We say that $\calT$ satisfies
$L^2$ off-diagonal estimates if for some constants $C\ge 0$ and $\alpha>0$ 
for all closed sets
$E$ and
$F$, all
$h
\in L^2$ with support in 
$E$ and all $t>0$ we have
\begin{equation}\label{eq:offL2}
\|T_t  h\|_{L^2(F)} \le C e^{-\frac{c d(E,F)^2}{t}} \|h\|_2.
\end{equation}
Here, and subsequently, $d(E,F)$ denotes the semi-distance induced on sets by the
Euclidean distance.  
\end{defi}

In case $\calT=(T_t)_{t>0}$ is replaced by a family $\calT=(T_z)_{z\in \Sigma_\mu}$  defined on a complex sector  $\Sigma_\mu$ with $0\le \mu<\frac \pi 2$,
then  we adopt the same definition and replacing $t$ by $|z|$ in the right hand side of the inequality. In this case, the constants
$C$ and $\alpha$ may depend on the angle $\mu$.

\begin{proposition} For all $\mu \in (0, \frac \pi 2-\omega)$, the families  $(e^{-zL})_{z\in \Sigma_\mu}$,  $(zLe^{-zL})_{z\in \Sigma_\mu}$ and $(\sqrt z\, \nabla
e^{-zL})_{z\in \Sigma_\mu}$ satisfies
$L^2$ off-diagonal estimates.  
\end{proposition}  

\begin{proof} We begin with the case of real times $t>0$. Let
$\varphi$ be a Lipschitz function on $\RR^n$ with Lipschitz norm 1 and $\rho>0$. 
We may define $L_\rho=e^{\rho\varphi} L e^{-\rho\varphi}$ by the form method.  
This operator is of second order type with same principal term as
$L$ and lower order terms with bounded coefficients. More precisely, let
$Q_\rho$ be the associated form. Then it is bounded on $W^{1,2}$ and  one can
find $c$ depending only on dimension and the ellipticity constants of $L$ (not on
$\varphi$)
 such
that
$$
\Re Q_\rho(f) \ge \frac \lambda 2 \|\nabla f\|_2^2 -c\rho^2 \|f\|_2^2, \ f \in
W^{1,2}.
$$
The construction guarantees that $L_\rho + c\rho^2$ is maximal-accretive on
$L^2$.  Hence, the semigroup
$(e^{-tL_{\rho}})_{t>0}$  exists on $L^2$ and its analyticity gives us: 
$$
\| e^{-tL_{\rho}}f\|_2 + \|t\frac d {dt} e^{-tL_{\rho}}f\|_2 + \|\sqrt t \,
\nabla e^{-tL_{\rho}}f\|_2 \le C e^{c\rho^2 t} \|f\|_2
$$
for all $t>0$ where $C$ depends  on the ellipticity constants  of $L$ and
dimension only. Let $E$ and $F$ be two closed sets and $f \in L^2$, with
compact support contained in
$E$. Choose
$\varphi(x)=d(x,E)$. Then, 
 \begin{equation}\label{Lrho}
e^{-tL}f  = e^{-\rho \varphi} e^{-tL_{\rho}}f .
\end{equation}
Hence, for all $t>0$, all $\rho>0$
$$\|e^{-tL}f\|_{L^2(F)} \le C e^{-\rho d(E,F)}e^{c\rho^2 t} \|f\|_2.
$$
Optimizing with respect to $\rho >0$ yields 
$$
\|e^{-tL}f\|_{L^2(F)} \le C e^{-\frac{d(E,F)^2} {4ct}} \|f\|_2.
$$
Next, differentiating \eqref{Lrho}, one has
$$ 
\frac d {dt} e^{-tL}f = e^{-\rho \varphi} \frac d {dt}e^{-tL_{\rho}}f
$$
and the same argument applies. Eventually, applying the gradient operator to 
\eqref{Lrho} yields
$$
\nabla e^{-tL}f = -\rho   (e^{-\rho \varphi} e^{-tL_{\rho}}f)(\nabla \varphi ) +
e^{-\rho \varphi}
\nabla 
e^{-tL_{\rho}}f.
$$
Hence,
$$
\| \nabla 
e^{-tL}f\|_{L^2(F)} \le   C\rho  e^{-\frac{d(E,F)^2} {4ct}} \|f\|_2 +   C
t^{-1/2}     e^{-\rho d(E,F)}e^{c\rho^2 t}
\|f\|_2
$$
and choosing $\rho =\frac{ d(E,F)}{2ct}$ yields
$$
\|\nabla e^{-tL}f\|_{L^2(F)} \le \frac{C}{\sqrt t} \left( 1 + \frac {d(E,F)}{\sqrt
t}\right) e^{-\frac{d(E,F)^2} {4ct}}
\|f\|_2.
$$
A density argument (since $f$ was supposed with compact support in $E$)
concludes the proof for real times. To go to complex times, we notice that this
applies to
$e^{i\alpha}L$, which is an operator in the same class as $L$ (with coefficients
$e^{i\alpha} A(x)$) as long as
$|\alpha| <
\frac \pi 2 -\omega$. Hence, the above estimates apply and the constants remain
uniform when
$\alpha$ is contained in a compact subset of $(-\frac \pi 2 +\omega, \frac\pi 2
-\omega)$. Finally, the desired estimates follow easily from the observation
that $e^{-zL}=e^{-t(e^{i\alpha} L)}$ when 
$z=te^{i\alpha}$. 
 \end{proof}

\subsection{Square root}\label{sec:squareroot}

As $L$ is a maximal-accretive operator, it has a square root, which we denote by $L^{1/2}$, defined as the unique maximal-accretive operator 
such that 
$$
L^{1/2}L^{1/2}=L
$$
as unbounded operators.~\footnote{\, For an explicit construction, we refer the
reader to Kato's book \cite[p. 281]{K2}   or to Meyer and
Coifman's book \cite[Chapter XIV]{MeCo}.}  Many formulas can be used to compute
$L^{1/2}$. The one we are going to use is
$$
L^{1/2}f= \pi^{-1/2}\int_0^\infty e^{-t L} Lf \, \frac{dt}{\sqrt t}.
$$
 This equality is valid as a Bochner integral when $f \in \calD(L)$ and 
as the limit in $L^2$ of the truncated Bochner integrals $\int_{\ep}^{R} \ldots$ as $\ep \downarrow 0$ and $R\uparrow \infty$
when $f \in \calD(L^{1/2})$. Also, this construction implies
that ${\cal D}(L)$ is dense in  ${\cal D}(L^{1/2})$.

The determination of the domain  of the square root of $L$ has become known as the
Kato square root problem and it is now a theorem in all dimensions, as recalled in the Introduction, that
${\cal D}(L^{1/2})=W^{1,2}$ with 
\begin{equation}\label{eq5}
\| L^{1/2} f \|_2 \sim \| \nabla f\|_2
\end{equation}
for all $f \in {\cal D}(L)$, hence by density in $\dot W^{1,2}$. In particular,
$L^{1/2}$ extends  to an isomorphism from $\dot W^{1,2}$ to $L^2$ and the formula 
$$
g =L^{1/2}L^{-1/2}g
$$ 
extends to all $L^2$ functions $g$.  

 This implies the following representation formula

\begin{lemma}\label{lem:rep} If $f, h \in \dot W^{1,2}$ then 
$$ \langle {(L^*)}^{1/2}f, L^{1/2}h \rangle  = \int_{\RR^n} \nabla f \cdot \ol{ A\nabla h} .
$$
\end{lemma}

\begin{proof} Since $f,h \in \dot W^{1,2}$, $(L^*)^{1/2}f, L^{1/2}h  \in L^2$.  Hence both sides of the equality are well-defined. It suffices to obtain  the
equality  if, in addition,
$h
\in {\cal D}(L)$, as ${\cal D}(L)$ is dense  in $\dot W^{1,2}$. 

In this  case, $L^{1/2}h $ belongs to the domain of $L^{1/2}$, so 
$ \langle {(L^*)}^{1/2}f, L^{1/2}h \rangle=  \langle f, Lh \rangle $ and the latter is equal to
$\int_{\RR^n} \nabla f \cdot \ol{ A\nabla h}$ by construction of $L$.  
\end{proof}

\subsection{The conservation property}\label{sec:conservation}

For real operators, the semigroup is contracting on $L^\infty$ and 
the conservation property $e^{-tL}1=1$ is a classical consequence of the
probabilistic interpretation of the semigroup or of the maximum principle for
parabolic equations.   But for complex operators,  the semigroup may not act
from 
$L^\infty$ into
$L^\infty$ (see the section on $L^p$ theory for the semigroup). Yet, thanks to
$L^2$ off-diagonal estimates, the action of the semigroup on $L^\infty$ can be
defined in the $L^2_{loc}$ sense and the conservation property still holds in
this sense:~\footnote{\, This is proved under $L^1-L^\infty$ off diagonal estimates
of the semigroup in \cite{AT} (Chapter I, Proposition 25) and is a consequence of
Corollary 4.6 in
\cite{ABBO} under the weaker
$L^\infty$-boundedness. It is mentioned in all generality but without proof in
\cite{AHLMcT}, p. 638.} 
\begin{equation}\label{conservation}
e^{-tL}1=1, \quad t>0.
\end{equation}

The first step is to show that $e^{-tL}$ maps $L^2$ functions with compact supports
in $L^1$. Fix $t>0$ and $\phi \in L^2$ supported in a cube $Q$. Cover
$\RR^n$ with a family of nonoverlapping cubes $(Q_k)$ of constant size with $Q_0 \supset 2Q$. 
Using the $L^2$ off-diagonal estimates, 
\begin{equation}\label{qk}
\int_{Q_k} |e^{-tL}\phi| \le |Q_k|^{1/2} e^{-\frac{c d(Q_k,Q)^2}{t}} \|\phi\|_2
\end{equation}
so that  summing in $k$ gives us the result.

Applying the first step to $L^*$ means that one can define
$e^{-tL}1$ in $L^2_{loc}$ by
$$
\int_{\RR^n} e^{-tL}1\,\, \ol\phi  = \int_{\RR^n}\, \ol{e^{-tL^*}\phi} 
$$
for all $\phi$ in $L^2$ with compact support. 

Next, let $\calX$ be a smooth function with $\calX(x)=1$ if $|x| \le 1$
and 
$\calX(x)=0$ if $|x| \ge 2$. Let $\calX_R(x)= \calX(x/R)$ for $R>0$. If $\phi$
 is an $L^2$ compactly supported function, then for $R>0$ and $t>0$
 \begin{equation}\label{eqetL1}
\int_{\RR^n} e^{-tL}1\,\, \ol\phi = \int_{\RR^n} \calX_R\, \ol{e^{-tL^*}\phi} + \int_{\RR^n} (1-\calX_R)\, \ol{e^{-tL^*}\phi}. 
 \end{equation}
 We use this representation twice, first
to show that the left hand side does not depend on $t>0$ and, second, to find
$\int_{\RR^n} \ol{\phi}$ as its value. This, indeed, shows that $e^{-tL} 1=
1$ in the sense of $L^2_{loc}$.

Let us begin with differentiating (\ref{eqetL1}) with respect to $t$. Indeed, the
first step applies also to $\frac{d}{dt}  e^{-tL^*}$ by the $L^2$ off-diagonal
estimates and this allows us to use the Lebesgue differentiation theorem to see
that 
$$\frac{d}{dt} \int_{\RR^n} e^{-tL}1\,\, \ol\phi  = 
\int_{\RR^n} \calX_R\,  \ol{\frac{d}{dt} e^{-tL^*}\phi} + \int_{\RR^n} (1-\calX_R)\, \ol{\frac{d}{dt}e^{-tL^*}\phi}.
$$
Fix $t>0$ and let $R\to \infty$. By Lebesgue dominated convergence, the latter integral  tends to 0.
Now,  since $\nabla\calX_R \in L^2$ and $ \frac{d}{dt} e^{-tL^*}\phi \in  {\cal D}(L^*)$,  we have that
$$\int_{\RR^n} \calX_R\, \ol{\frac{d}{dt} e^{-tL^*}\phi} =
\int_{\RR^n} A\nabla\calX_R\cdot \ol{\nabla e^{-tL^*}\phi}.
 $$
Again, using the $L^2$ off-diagonal estimates for $\nabla e^{-tL^*}$ and arguing as for \eqref{qk},
this integral is bounded by $CR^{n/2-1} e^{-cR^2/t}\|\phi\|_2$ for $R$ large enough so that it tends to 0.
This shows that the left hand side of (\ref{eqetL1}) is independent of $t>0$. In 
the right hand side,  choose and fix $R$ large enough so that the supports of
$\phi$ and
$(1-\calX_R)$ are far apart. It follows from \eqref{qk} that 
$\int_{\RR^n} (1-\calX_R)\, \ol{e^{-tL^*}\phi}$ tends to $0$ with $t$  by
dominated convergence.  Eventually, since $e^{-tL^*}$ is a strongly continuous  in
$L^2$ at $t=0$, we obtain that $\int_{\RR^n} \calX_R\, \ol{e^{-tL^*}\phi}$ tends to
$\int_{\RR^n} \calX_R\, \ol{\phi}= \int_{\RR^n}  \ol{\phi}$ as $t$ tends to 0. This
proves (\ref{conservation}).

\section{$L^p$ theory for the semigroup}\label{sec:Lpsemigroup}

This section is devoted to establishing the basic properties concerning uniform
boundedness of the semigroup $(e^{-tL})_{t>0}$ and of the family $(\sqrt t\, \nabla e^{-tL})_{t>0}$
on $L^p$ spaces.

\subsection{Hypercontractivity and uniform boundedness}\label{sec:hyper}

The point of this section is to present a general statement allowing to pass from 
hypercontractivity properties for the semigroup to uniform boundedness properties.
The bridge between both are  off-diagonal estimates.~\footnote{\, This tool was
introduced for that purpose by Davies \cite{Da}. See \cite{LSV} for further
historical comments. This notion turns out to be equivalent to some $L^p-L^q$
boundedness for perturbed semigroups by exponential weights (See \cite{LSV}). We
do not need to go into this here.}

We introduce a few definitions.  

\begin{defi}\label{LpLqbdd} Let $\calT=(T_t)_{t>0}$ be a family of uniformly bounded operators on $L^2$. 
 We say that $\calT$ is $L^p-L^q$ bounded for some $p,q \in [1, \infty]$ with $p\le q$ if for some constant $C$, for all $t>0$
and all $h\in L^p \cap L^2$
\begin{equation}\label{eq:Spq}
\|T_t  h\|_{q} \le C t^{-\gamma_{pq}/2}  \|h\|_p.
\end{equation} We say that $\calT$ satisfies $L^p-L^q$ off-diagonal estimates for some $p,q \in [1, \infty]$ with $p\le q$ if  
 if for some constants $C,c>0$,
for all closed sets
$E$ and
$F$, all
$h
\in L^p \cap L^2$ with support in
$E$ and all $t>0$ we have
\begin{equation}\label{eq:offLpLq}
\|T_t  h\|_{L^q(F)} \le C t^{-\gamma_{pq}/2} e^{-\frac{c d(E,F)^2}{t}} \|h\|_p.
\end{equation} 
\end{defi}

Recall that the numbers $\gamma_{pq}$ are defined in the notation section. Such estimates depend on dimension and on the ``parabolic'' character of the family
through this number.  If
$p=q$ we speak of
$L^p$ boundedness and
$L^p$ off-diagonal estimates.  Remark that these notions dualize. A family $(T_t)_{t>0}$ is $L^p-L^q$ bounded (resp. satisfies 
$L^p-L^q$ off-diagonal estimates if and only if the dual family $((T_t)^*)_{t>0}$ is $L^{q'}-L^{p'}$ bounded
(resp. satisfies $L^{q'}-L^{p'}$ off-diagonal estimates).

Let us also state a useful result whose easy proof is skipped.

\begin{proposition}\label{propcomposition} If $(T_t)_{t>0}$ satisfies $L^p-L^q$ boundedness (resp. off-diagonal estimates) and $(S_t)_{t>0}$ satisfies  $L^q-L^r$
 boundedness (resp. off-diagonal estimates) then $(S_tT_t)_{t>0}$ satisfies $L^p-L^r$ boundedness (resp. off-diagonal estimates).
\end{proposition}

For a semigroup, the terminology \textbf{hypercontractivity} is often used for  $L^p-L^q$ boundedness for some $p<q$.
The relation between hypercontractivity, boundedness and off-diagonal estimates is the following result.

\begin{proposition} \label{propequiv} Let $p\in [1,2)$ and $n\ge 1$. Let
$\calS=(e^{-tL})_{t>0}$. 
\begin{enumerate}
\item If $\calS$ is $L^p$ bounded    then it is 
$L^p-L^2$ bounded.  
\item If $\calS$ is $L^p-L^2$ bounded,  then for all $q\in (p,2)$ 
it satisfies $L^q-L^2$ off-diagonal
estimates. 
\item If $\calS$ satisfies $L^p-L^2$ off-diagonal estimates then it is $L^p$ bounded. 
\end{enumerate}
\end{proposition}

\begin{remark} The result applies when $2<p\le \infty$  by duality:
replace $L^p-L^2$ by $L^2-L^{p}$ everywhere. We have privileged the central role of $L^2$ for reasons of simplicity and usefulness.
$L^2$ could be replaced by $L^q$ for $q$ larger than $2$ if necessary. It occurs in Section \ref{sec:hodge}.
\end{remark}

\begin{proof} The proof of item 1 is obtained from Nash type
inequalities.~\footnote{\, The proof follows that of \cite[Theorem II.3.2]{VSC} 
done for $p=1$.} We start from the Gagliardo-Nirenberg inequality 
$$
\|f\|_2^2 \le C\|\nabla f\|_2^{2\alpha} \|f\|_p^{2\beta}
$$ with $$\alpha+\beta=1 \quad \mathrm{and}\quad (1+\gamma_p) 
\alpha =\gamma_p.$$ This yields the Nash  inequality
$$
\|e^{-tL}f\|_2^2 \le C\|\nabla e^{-tL}f\|_2^{2\alpha} \|e^{-tL}f\|_p^{2\beta}
$$
for all $t>0$ and $f \in L^2\cap L^p$.
By ellipticity, one has
$$
\|\nabla e^{-tL}f\|_2^{2} \le \lambda^{-1} \Re \int_{\RR^n} A\nabla e^{-tL}f \cdot \ol{ \nabla e^{-tL}f}  = -  (2\lambda)^{-1} \frac d {dt} \|
e^{-tL}f\|^2_2.
$$
Assume $f \in L^2\cap L^p$ with $\|f\|_p=1$.
Using $L^p$ boundedness of the semigroup in the Nash inequality,   one obtains the differential inequality
 $$\varphi(t) ^{1/ \alpha} \le  - C\varphi' (t),$$
where $\varphi(t)= \|
e^{-tL}f\|^2_2$. 
Integrating between $t$ and $2t$ and using that $\varphi$ is nonincreasing, one finds easily that
$$\varphi(t) \le C t^{-\frac \alpha{\alpha -1}}=Ct^{-\gamma_{p}},$$ 
which is the desired estimate. 

The proof of item 2 consists in interpolating by the Riesz-Thorin theorem    the
$L^p-L^2$ boundedness assumption with the $L^2$ off-diagonal estimates, once we fix the sets $E$ and $F$ in the definition
of the off-diagonal estimates.

The proof of item 3  can be seen by invoking the following simple lemma which has nothing to do with semigroups.

\begin{lemma}\label{lemmaDa}\footnote{\, This result is implicit in the 
proof of \cite[Theorem 25]{Da}.} If $1\le p\le q \le \infty$ and 
 $T$ is a linear operator which satisfies $L^p-L^q$ off-diagonal estimates in the form
$\|Tf\|_{L^q(F)} \le g(d(E,F)) \|f\|_{L^p(E)}$ whenever $E,F$ are closed cubes and $f$ is supported in 
$E$ and $g$ is some function.  Then   $T$ is bounded on $L^p$ with norm bounded by $s^{\gamma_{pq}} \sum_{k\in \ZZ^n} g(\sup(|k|-1,0)s)$
for any $s>0$ provided this sum is finite.
\end{lemma}

\begin{proof} We may assume $p\le q<\infty$.    Let $(Q_k)_{k \in \ZZ^n}$ be a partition of $\RR^n$ by  cubes
having sidelength
$s$ and
${\cal X}_\ell$ be the indicator function of $Q_\ell$. Then
\begin{align*}
\|Tf\|_q^q& = \sum_k\left\|\sum_\ell T(f{\cal X}_\ell)\right\|_{L^q(Q_k)}^q 
\\
&\le  \sum_k \left(\sum_\ell g(\sup(|k-\ell|-1,0)s) \|f\|_{L^p(Q_\ell)}\right)^q
\\
 &\le  \sum_k \left(\sum_\ell g(\sup(|k-\ell|-1,0)s) s^{\gamma_{pq}} \|f\|_{L^q(Q_\ell)}\right)^q
\\
& \le  \left(\sum_{k} g(\sup(|k-\ell|-1,0)s) s^{\gamma_{pq}} \right)^q\sum_\ell \|f\|_{L^q(Q_\ell)}^q
\\
& =  \left(\sum_{k} g(\sup(|k|-1,0)s) s^{\gamma_{pq}}\right)^q\|f\|_q^q
\end{align*}
where we  used that the 
discrete convolution with an $\ell^1$ sequence  is bounded on $\ell^q$.
\end{proof}

We come back to the proof of the proposition. It suffices to apply this lemma to $T=e^{-tL}$ from $L^p$ to $L^2$ with $g(u) = c  t^{-\gamma_{p}/
2} e^{-cu^2/t}$ and choose $s= t^{1/2}$. This yields $L^p$ boundedness of $\calS$.
\end{proof}

\begin{remark} For the \textbf{same} $p$,  $L^p-L^2$ off-diagonal estimates
$\Longrightarrow$
$L^p$ boundedness $\Longrightarrow$ $L^p-L^2$ boundedness for the semigroup 
${\cal S}$. We do not know the status of the converses  for this class of semigroups. It would be of great interest when  $p=1$.
\end{remark}

\subsection{$W^{1,p}$ elliptic estimates and hypercontractivity}\label{sec:w1p}

In this section, we show how to obtain hypercontractivity properties 
 using $W^{1,p}$ elliptic estimates. We proceed independently of dimension although
specific arguments for dimensions 1 and 2 yield much better results.

\begin{lemma}\label{lemmaw1p}\footnote{\, In this generality, this is  \cite[Proposition 3.1]{AMcT2}. See also \cite{BLP} and \cite{IS}. 
Note that when $n=1$, invertibility holds for $1\le p \le
\infty$ \cite[Theorem 2.2]{AMcT2}. If $n\ge 2$, $r$ may
be arbitrary close to
$2$.  
} There is
an
$r\in [1,2)$ depending on dimension and the ellipticity constants only, such that
 $I+L$ extends to a bounded and invertible operator   from 
$W^{1,p}$ onto $ W^{-1,p}$  for $|\frac  1 2
-\frac 1 p| <  |\frac  1 2 -\frac 1 r|$.~\footnote{\, There is also the
corresponding homogeneous statement for $L$  from
$\dot W^{1,p}$ onto $\dot W^{-1,p}$. See Section \ref{sec:hodge}.}   
\end{lemma}

\begin{proof}  That $I+L$ is bounded from 
$W^{1,p}$ into $ W^{-1,p}$ for all $p$ with $1\le p\le \infty$ is obvious.

Let $A$ denote the  operator of pointwise multiplication with $A(x)$ and
$\|A\|$ its norm acting on $L^p$ spaces of 
$\CC^n$-valued functions (endowing matrices with the norm inherited from the
Hermitian structure on $\CC^n$). This is the same number for all $1\le p\le \infty$.
By ellipticity, there exists a large constant
$c$ such that
$\|A-c\,I\| <c$. Thus, $A=c(I+M)$ where
$\|M\|<1$. This means that 
$$L=(I-c\Delta - c\dv M \nabla )=J (I- J^{-1}\dv M \nabla J^{-1})
J,$$ where
$\Delta$ denotes the ordinary Laplacian and $J=c^{-1/2}(I-c\Delta)^{1/2}$. 

Remark that  by standard multiplier theorems (or kernel estimates and
Calder\'on-Zygmund theory) the array of  operators
$\nabla J^{-1}$  is bounded from  $L^p$  (of $\CC$-valued functions) to $L^p$
(of $\CC^n$-valued functions) for $1<p<\infty$ ($p=1$ and $p=\infty$ included if
$n=1$). Moreover, a Fourier transform argument shows that $c_2$=1 if $c_p$ denotes
the norm on $L^p$.

Thus $R=J^{-1}\dv M \nabla J^{-1}$  is bounded on $L^p$ for  $1 < p < \infty$  with
norm bounded above by $\|M\| c_p^2$.  Since 
$c_p$ is controlled by a convex function, the  operator norm  of
$R$ on
$L^p$ remains  less than 1 provided $p$ is close to 2. Therefore,
one can invert $I -R$ by a converging Neumann series
in the space bounded operators on $L^p$ for $p$ close to 2.  Since $J$ is bounded
and invertible from $W^{s,p}$ onto $W^{s-1,p}$  for $s=0,1$ and $1<p<\infty$, this
proves the invertibility of $I+L$ from 
$W^{1,p}$ onto $ W^{-1,p}$ for $p$ close to 2.
\end{proof}

\begin{corollary}\label{corohyper} Let $r_*=\frac{nr}{n+r}$ with $r$ as in Lemma \ref{lemmaw1p}. Let $p\in [1,2]$ be such that   
$$
\begin{cases} p>r_*,& \mathrm{if}\
r_*\ge 1 \\ 
p=1& \mathrm{if}\ r_*<1.
\end{cases}
$$ 
The semigroup $(e^{-tL})_{t>0}$ and its dual 
$(e^{-tL^*})_{t>0}$ are  $L^p-L^2$ bounded and the best  constant $C$ in 
 \begin{equation}\label{eq24}
\|e^{-tL}f\|_{2} + \|e^{-tL^*}f\|_{2}\le C{t}^{-\gamma_p/2}\|f\|_{p}, \quad f\in L^2\cap L^p
\end{equation}
depends only on dimension, ellipticity and  $p$. 
\end{corollary} 

\begin{proof} Assume first that $t=1$.
By Lemma \ref{lemmaw1p} and the Sobolev embedding theorem, in a finite number of steps $(1+L)^{-k}$ extends to a bounded map from
$L^p$ into 
$L^2$. Note that $k$  depends only on $r$, hence ellipticity and
dimension.  Let $f\in L^2\cap L^p$. Since $f$ is in
$L^2$, the equality 
$$e^{-L}f =e^{-L}(I+L)^{k}(I+L)^{-k}f $$
is justified. As
 $e^{-L}(I+L)^{k}$ extends to bounded operator on $L^2$ by analyticity of the semigroup on $L^2$,
we have obtained that 
$\|e^{-L}f\|_{2} \le C\| f\|_{p}$, with a constant $C$ that depends only on ellipticity, dimension and $p$. 

If $t\ne 1$, then the   affine change of variable in $\RR^n$ defined by $g(x)= f(t^{1/2}x)$ gives us 
$e^{-tL}f(x)=(e^{- L_{t} }g)(t^{-1/2}x)$ 
with $L_{t}$ the second order operator associated with the matrix of coefficients $A(t^{1/2}x)$. Since $L_{t}$ has same ellipticity constants 
as $L$, the previous bound applies and yields the desired estimate in function of $t$ by change of variables. 

The same argument applies to $L^*$.
\end{proof}  

\begin{corollary}\label{corolpbdd}
\footnote{\, For $n=1$, this is in \cite{AMcT2}, Theorem 2.21, for $n=2$ in
\cite{AMcT2}, Theorem 3.5. There, 
 $L^1$ and $L^\infty$ boundedness and even $L^1-L^\infty$
off-diagonal estimates are proved, but further arguments such as the so-called
Davies' trick are needed. This method, sufficient for our needs here, is not
powerful enough. For 
$n\ge 3$, this follows from combining the works of Davies \cite{Da} for $p_n\le p
\le (p_n)'$ and the perturbation method in 
\cite{AT}, Chapter~I.} If
$n=1$ or
$2$, then  the semigroup
$(e^{-tL})_{t>0}$ is $L^p-L^2$ bounded for $1\le p<2$ and $L^2-L^p$ bounded for $2< p\le \infty$. Hence, it is $L^p$-bounded for all $p \in (1, \infty)$.  
 If $n\ge 3$, there exists $\ep>0$ depending only on dimension and the ellipticity
constants   such that the semigroup $(e^{-tL})_{t>0}$
is $L^p$-bounded for all $p$ contained in the interval
$(p_n-\ep, (p_n-\ep)')$ for $p_n=\frac {2n}{n+2}$.
\end{corollary}

\begin{proof} In dimensions $n=1,2$, we have $r_*<1$. In dimension $n\ge 3$, the value of $r_*$ may or may not
be less than 1, but it is less than $2_*=\frac {2n}{n+2}$.      It suffices to combine Corollary \ref{corohyper} and Proposition 
\ref{propequiv} to finish the proof.
\end{proof}

We introduce here two critical exponents for $L$. Let $\calJ(L)$ denote the maximal 
interval of exponents $p$ in $[1,\infty]$ for which 
the semigroup $(e^{-tL})_{t>0}$ is $L^p$ bounded.~\footnote{\, In \cite{LSV}, a
systematic study of this interval is made for a class of real elliptic
operators. The methods heavily rely on real functions.}  We write
$\mathrm{int}\calJ(L)= (p_-(L), p_+(L))$. Note that 
$(p_+(L))'= p_-(L^*)$ and vice versa. One has shown $p_+(L)=\infty$ and $p_-(L)=1$ if $n=1,2$, and $p_+(L)> \frac {2n}{n-2}$ and $p_-(L)<\frac {2n}{n+2}$ if $n\ge
3$. In specific situations, much more can be said. For example, we have from the well-known formula for the
heat kernel
$p_-(-\Delta)=1$,
$p_+(-\Delta)=\infty$. From the maximum principle for real parabolic equations, one 
has $p_-(L)=1$, $p_+(L)=\infty$ if $L$ has real coefficients. 

\subsection{Gradient estimates}

Let us consider the possible estimates for $\nabla e^{-tL}$. Let $\calN(L)$ denote
the maximal  interval (if nonempty) of exponents $p$ in $[1,\infty]$ for which 
the family $(\sqrt t\, \nabla e^{-tL})_{t>0}$ is $L^p$ bounded. We write
$\mathrm{int}\calN(L)= (q_-(L), q_+(L))$.  A dichotomy  between the cases
$p>2$ and $p<2$ appears immediatly.

\begin{proposition}\label{propselfimprovp<2}Let $1\le p<2$. 
If $(\sqrt t\, \nabla e^{-tL})_{t>0}$ is
$L^p$-bounded, then $(e^{-tL})_{t>0}$ is $L^q$-bounded for
$p<q<2$. Conversely, if $(e^{-tL})_{t>0}$ is $L^p$-bounded then $(\sqrt t\, \nabla e^{-tL})_{t>0}$ is
$L^q$-bounded for
$p<q<2$. Hence $q_-(L)=p_-(L)$.
\end{proposition}

\begin{proof} Let us see the direct part first.  We have nothing to prove if $n\le 2$ by Corollary \ref{corolpbdd} as the conclusion holds true. If $n\ge 3$, by
interpolation  and  Sobolev embeddings we have that  $(e^{-tL})_{t>0}$ is $L^q-L^{q^*}$ bounded for $p\le q \le 2$. Let $p_k$ be defined by $p_0=p$, 
$p_{k+1}={(p_k)}^*$ and stop whenever
$2_*\le p_k<2$. By composition and the semigroup property, we have that $(e^{-tL})_{t>0}$ is $L^p-L^{p_k}$ bounded. By Corollary \ref{corohyper}, since $p_k\ge
2_*$,
$(e^{-tL})_{t>0}$ is $L^{p_k}-L^{2}$ bounded. By composition again, $(e^{-tL})_{t>0}$ is $L^p-L^2$ bounded. The conclusion follows from
Proposition \ref{propequiv}.  

For the converse, by Proposition \ref{propequiv}, $(e^{-tL})_{t>0}$ is $L^p-L^2$ bounded.
Since $(\sqrt t\, \nabla e^{-tL})_{t>0}$ is $L^2$ bounded,  this self-improves by composition to 
$L^p-L^2$ boundedness.  Lemma \ref{lemmaDa} applied to $T=\nabla e^{-tL}$  yields the conclusion.

The relation $q_-(L)=p_-(L)$ follows immediately.
\end{proof}

\begin{remark} \begin{enumerate}
\item We see that the semigroup acting on $L^p$ self-improves into $W^{1,p}$ if $p<2$ up to allowing the $p$'s to vary.  This is  false when
$p>2$. 
\item Although we do not need
such a refinement here, it would be interesting to know   whether the conclusions of $L^p$
boundedness hold at  the endpoint $p$.  The arguments show, nevertheless, that,  for the \textbf{same} $p$, 
$L^p-L^2$ boundedness for
$(e^{-tL})_{t>0}$  is equivalent to 
$L^p-L^2$ boundedness for $(\sqrt t\, \nabla e^{-tL})_{t>0}$,  and  $L^p-L^2$ off-diagonal estimates  for  $(e^{-tL})_{t>0}$  implies  $L^p-L^2$ off-diagonal
estimates for $(\sqrt t\, \nabla e^{-tL})_{t>0}$. The converse of the latter is not clear.~\footnote{\, In recent work with J.-M. Martell, we established the converse \cite{AM2}.}
\end{enumerate}
\end{remark}

Let us record here the following consequences of the above argument for later use.

\begin{corollary}\label{propequivgradp<2} Let $p\in [1,2)$ and $n\ge 1$. Let
$\calN=(\sqrt t\, \nabla e^{-tL})_{t>0}$. 
\begin{enumerate}
\item If $\calN$ is $L^p$ bounded    then it is 
$L^p-L^2$ bounded.  
\item If $\calN$ is $L^p-L^2$ bounded,  then 
it satisfies $L^q-L^2$ off-diagonal
estimates for $p<q<2$ . 
\item If $\calN$ satisfies $L^p-L^2$ off-diagonal estimates then it is $L^p$ bounded. 
\end{enumerate}

\end{corollary}

Next, let us consider the case $p>2$. We have the same statement as the corollary but with a different argument.

\begin{proposition} \label{propequivgradp>2} Let $p\in (2,\infty]$ and $n\ge 1$. Let
$\calN=(\sqrt t\, \nabla e^{-tL})_{t>0}$. 
\begin{enumerate}
\item If $\calN$ is $L^p$ bounded    then it is 
$L^2-L^p$ bounded.  
\item If $\calN$ is $L^2-L^p$ bounded,  then  
it satisfies $L^2-L^q$ off-diagonal
estimates for $2<q<p$. 
\item If $\calN$ satisfies $L^2-L^p$ off-diagonal estimates then it is $L^p$ bounded. 
\end{enumerate}
\end{proposition}

\begin{proof} Let us prove the first item. Assume that $n\ge 3$ and also that $p<\infty$
If $\nabla e^{-tL}$ is bounded on $L^p$ with bound $Ct^{-1/2}$, then the same is true for all $q\in [2,p]$. By Sobolev embeddings,
for all  $q\in [2,p]$ with $q<n$,
$$
e^{-tL}\colon L^q \to L^{q^*}$$
with bound $Ct^{-1/2}$.  Since we know that $
e^{-tL}\colon L^2 \to L^{q}$ for any $q \in \RR$ with $2\le q \le 2^*$ from Corollary \ref{corohyper},  it follows as in the proof
of Proposition \ref{propselfimprovp<2}
 that $e^{-tL}$ maps $L^2$ to $L^q$ for all $q \in \RR$ with $2\le q \le p^*$. 
Writing $\nabla e^{-tL}$ as $\nabla e^{-(t/2)L} e^{-(t/2)L}$ shows that this operator is bounded from 
$L^2$ into $L^p$ with the appropriate norm growth. 

Assume next $n\ge 3$ and $p=\infty$. We just need to prove that the semigroup is $L^2-L^\infty$ bounded and the rest of the argument applies.   
By Morrey's embedding, if $n<q<\infty$ and $\alpha=1 -\frac n q$, $t>0$ 
and $x,y \in \RR^n$, 
$$
|e^{-tL}f(x) -e^{-tL}f(y) | \lesssim \|\nabla e^{-tL}f\|_q |x-y|^\alpha \lesssim t^{-\frac 1 2 - \frac {\gamma_q } 2} \|f\|_2 |x-y|^\alpha.
$$
In the last inequality, we used that $\calN$ is $L^2-L^q$ bounded as we have just proved it. 
Fix $x$ and average the square of this inequality on a ball $B$ with center $x$ and
radius $r=\sqrt t$ to find
$$
|e^{tL}f(x)| \lesssim  |B|^{-1/2} \|e^{-tL}f\|_{L^2(B)} +  r^{\frac \alpha2}t^{-\frac 1 2 - \frac {\gamma_q } 2}\|f\|_2 \lesssim t^{-n/4} \|f\|_2
$$
by using the $L^2$ boundedness of $e^{-tL}$. This prove the $L^2-L^\infty$ boundedness of the semigroup. 

If $n\le 2$, we already know that $
e^{-tL}\colon L^2 \to L^{q}$ for any $q \in [2,\infty]$  from Corollary \ref{corohyper} and the rest of the above argument applies.

The second item is a consequence  of interpolation between the hypothesis and the $L^2$ off-diagonal estimates for $(\sqrt t\, \nabla e^{-tL})_{t>0}$.

  For the third item, it is enough  to apply Lemma \ref{lemmaDa} to $T= e^{-tL^*}\dv$ with $s=t^{1/2}$ and duality.
\end{proof}

\begin{corollary}\label{coroqplus>2} $q_+(L)>2$.
\end{corollary}

\begin{proof}
Let $p>2$ such that $|\frac  1 2
-\frac 1 p| <  |\frac  1 2 -\frac 1 r|$ where $r$ is given in Lemma \ref{lemmaw1p}. Let $f \in L^2$. After a finite 
number of steps (depending only on $r$) $(I+L)^{-k}$ maps $L^2$ into $W^{1,p}$. In
particular, since $(I+L)^k e^{-L}f \in L^2$ by analyticity,
$$
\nabla e^{-L}f=\nabla (I+L)^{-k} (I+L)^k e^{-L}f \in L^p.
$$
This proves the boundedness of $\nabla e^{-L}$ from $L^2$ into $L^p$.
The $L^2-L^p$ boundedness of the family $(\sqrt t\, \nabla e^{-tL})_{t>0}$ follows by the rescaling argument
of the proof of Corollary \ref{corohyper}. Hence, $q_+(L)\ge p$.
\end{proof}

\begin{corollary}\label{coroglopen} $\calN(L)$ is not empty and contains a
neighborhood of 2.
\end{corollary}

\begin{proof} We have $q_-(L)=p_-(L)< 2$ and $q_+(L)>2$.
\end{proof}

More is true in dimension 1.

\begin{proposition}\label{propqplusn1} If $n=1$, we have $q_+(L)=\infty$.
\end{proposition}

\begin{proof}  In dimension 1, the operator $L$ takes the form 
$-  \frac{d}{dx}\big( a \frac{d}{dx}\big)$.  By Proposition
\ref{propequivgradp>2}, it suffices to establish  that
$(\sqrt t\,
\frac{d}{dx} e^{-tL})_{t>0}$ is
$L^2-L^\infty$ bounded.~\footnote{\, In fact, this family satisfies the
$L^1-L^\infty$ off-diagonal estimates, which implies all the other boundedness
properties. See \cite[p.~728]{AT1} where it is mentioned, and \cite[Theorem 2.21]{AMcT2}
for details.} 
Let
$f\in L^2$ and
$t>0$ and set
$g = \sqrt t\, a \frac{d}{dx} e^{-tL}f$. We know that $g \in L^2$ with norm $O(1)$, while $g' = -\sqrt t\, L e^{-tL}f \in L^\infty$ with
norm $O(t^{-1/2})$. The last fact can be seen from   writing $L e^{-tL}=e^{-(t/2)L} Le^{-(t/2)L}$ and using that
 $(tLe^{-tL})_{t>0} $ is $L^2$ bounded and $( e^{-tL})_{t>0}$ is $L^2-L^\infty$ bounded by Corollary  \ref{corolpbdd}. This implies that $g \in L^\infty$
with norm
$O(t^{-1/4})$. Indeed, let $I$ be an interval of size $t^{1/4}$ and $x,y \in I$. We have
$g(x)-g(y)=\int_y^xg'(s)\, ds$. Hence, $|g(x)| \le |g(y)| + Ct^{-1/4}$. Averaging squares over $I$ with respect to $y$ yields
$|g(x)|^2 \le  Ct^{-1/2} \|g\|_2^2 + Ct^{-1/2} = O(t^{-1/2})$. Since $I$ and $x$ are arbitrary, this gives us the desired $L^\infty$ bound
on $\RR$. 
\end{proof}

\begin{proposition}\label{propselfimprovp>2} Assume $2<p  \le \infty$ and $(\sqrt t\, \nabla e^{-tL})_{t>0}$ is $L^p$-bounded.
\begin{enumerate} 
\item   $(e^{-tL})_{t>0}$ is $L^q$-bounded for   $2 \le q \le p^*$ except when $p=n$ for which we conclude only for $2\le q<\infty$. In particular, one has
$p_+(L) \ge (q_+(L))^*$.
\item  If $p<\infty$, $( t\, \nabla e^{-tL}\dv)_{t>0}$ is $L^p$-bounded. 
\end{enumerate}
\end{proposition}

\begin{proof} The first part has been seen in the proof of Proposition \ref{propequivgradp>2}. 
For the second part,  the   $L^p$ boundedness of $(\sqrt t\, \nabla e^{-tL})_{t>0}$ implies in particular
that  $(e^{-tL})_{t>0}$ is $L^q$ bounded for $q$ in a neighborhood of $p$, thus, that  $(e^{-tL^*})_{t>0}$ is $L^q$ bounded for $q$ in a neighborhood of $p'$.
Applying Proposition \ref{propselfimprovp<2} shows, in particular, that $(\sqrt t\, \nabla e^{-tL^*})_{t>0}$ is $L^{p'}$-bounded, hence,
that $( \sqrt t\,  e^{-tL}\dv)_{t>0}$ is $L^p$-bounded by duality. We conclude by composition.
\end{proof}

\begin{remark} If $n=1$, combining Corollary \ref{corolpbdd}, Proposition
\ref{propqplusn1} and Proposition \ref{propselfimprovp>2}, one has that the
semigroup
is bounded from $W^{-1,p}$ into $W^{1,p}$ for $1<p<\infty$ and $t>0$.~\footnote{\,
More is true, in fact, this holds for the endpoints $p=1$ and $p=\infty$ and
one even has $L^1-L^\infty$ off-diagonal estimates for the family $(t\frac d{dx}
e^{-tL} \frac d {dx})_{t>0}$. See [AMcT].} 
\end{remark}

\subsection{Summary}

This section is nothing but a summary of results obtained so far and 
gathered here for the reader's convenience. 

Recall that $\calJ(L)$ denotes the maximal 
interval of exponents $p$ in $[1,\infty]$ for which 
the semigroup $(e^{-tL})_{t>0}$ is $L^p$ bounded. We write $\mathrm{int}\calJ(L)= (p_-(L), p_+(L))$. Note that 
$(p_+(L))'= p_-(L^*)$ and vice versa. One has $p_+(L)=\infty$ and $p_-(L)=1$ if $n=1,2$ and $p_+(L)> \frac {2n}{n-2}$ and $p_-(L)<\frac {2n}{n+2}$ if $n\ge 3$.

Next,  $\calN(L)$ denotes the maximal  interval  of exponents $p$ in $[1,\infty]$
for which  the family $(\sqrt t\, \nabla e^{-tL})_{t>0}$ is $L^p$ bounded. We
write $\mathrm{int}\calN(L)= (q_-(L), q_+(L))$.   We have $q_-(L)=1$ and
$q_+(L)=\infty$ if $n=1$, $q_-(L)=1$ if $n=2$, $q_-(L)<\frac {2n}{n+2}$ if $n\ge
3$ and $q_+(L)>2$ if $n\ge 2$. In particular, $\mathrm{int}\calN(L)$ is a
neighborhood of 2.

The relation between $p_\pm(L)$ and $q_\pm(L)$ are as follows:
\begin{align*}
q_-(L)&=p_-(L), 
\\
p_+(L) &\ge (q_+(L))^*.
\end{align*}

\subsection{Sharpness issues}\label{sec:sharpness}

Denote by $\calE$ the class of all complex elliptic operators under study.

It is known  that if $n\ge 5$, there exists  $L\in \calE$ for which $p_-(L)>1$.~\footnote{\, This is
proved in \cite{ACT}, based on counterexamples built in
\cite{MNP}. See \cite{AT} for a shorter argument due to Davies.}

It is not known if the inequality $p_-(L) < \frac{2n}{n+2}$ ($n\ge 3$), or equivalently by taking the
adjoint
$p_+(L) > \frac{2n}{n-2}$ is sharp: that is, if given $p<\frac{2n}{n+2}$, there exists $L\in \calE$ for
which $p_-(L)>p$.

In view of this a natural conjecture is

\begin{conjecture} If $n\ge 3$, the   inequality  $p_-(L)<\frac{2n}{n +2}$ is sharp for the class $\calE$.
\end{conjecture}

It is known that  the inequality $q_+(L)>2$ is sharp in dimensions $n\ge 2$: for any $\ep>0$, there exists 
$L\in \calE$ (in fact, $L$ is real and symmetric) such that $q_+(L)\le 2+\ep$.~\footnote{\, This follows
from Meyers' example. See \cite{AT}.} 
As a consequence, there is no  general upper bound of $p_+(L)$ in terms of $q_+(L)$ since $p_+(L)=\infty$
if $n=2$.

 It is not known whether the inequality $p_+(L) \ge q_+(L)^{*}$ is best possible ($n\ge 3$):
that is, if given
$\ep>0$, one can find $L\in \calE$ with $p_+(L)<q_+(L)^{*} +\ep$.

\subsection{Analytic extension}

For technical reasons, we often need to apply the above results to the analytic extension of the semigroup
associated to $L$. We come to this now. 

 The  definition \ref{LpLqbdd} of $L^p-L^q$ boundedness and $L^p-L^q$ off-diagonal estimates 
 applies to a  family  $\calT=(T_z)_{z\in \Sigma_\beta}$  defined on a
complex sector 
$\Sigma_\beta$ with
$0\le
\beta<{\frac \pi 2}$,  replacing $t$ by $|z|$ in the right hand side of the
inequalities \eqref{eq:Spq} and \eqref{eq:offLpLq}. 

 Proposition  \ref{propcomposition} extends right away to such complex families. The statement of Proposition \ref{propequiv}
is also true for the analytic extension of the semigroup. 
However, more is true. It suffices to 
know boundedness of the semigroup to obtain 
 all properties for its analytic extension with optimal angles of the
sectors, that is, the sectors are 
$p$-{\bf independent}.~\footnote{\, The history of $p$-independence of sector of
holomorphy for $1\le p<\infty$ begins with that of $p$-independence of spectra. 
 Hempel and Voigt \cite{HV} proved that the spectrum  is $p$-independent for a large
class of Schr\"odinger operators acting on $L^p(\RR^n)$. 
Arendt \cite{Are} extended this to elliptic operators on domains of $\RR^n$  under
the assumption that the semi-group of the given operator satisfies
$L^1-L^\infty$ off-diagonal estimates.  Ouhabaz realized in his PhD thesis
(published later in \cite{Ou}) that the same assumption  yields holomorphy of the
semi-group up to $L^1$ with optimal angles in the case of a self-adjoint operator
on a subset of $\RR^n$. Independently, Arendt \& ter Elst \cite[Theorem 5.4]{AE},  and
Hieber
\cite{Hi} removed the self-adjointness  assumption. Finally, Davies \cite{Da2}
extended Ouhabaz argument to the setting of doubling spaces (for non-negative
self-adjoint operators) to obtain in a simpler manner the result of $p$-independence
of  spectrum of Arendt.  } 

Recall that  $\calJ(L)$ (resp. $\calN(L)$) is the maximal  interval of exponents
$p$ in $[1,\infty]$ for which the semigroup $(e^{-tL})_{t>0}$ (resp. the family
$(\sqrt t\, \nabla e^{-tL})_{t>0}$)
is 
$L^p$ bounded.

\begin{proposition}\label{propequivz} Let $\omega$ 
be the type of $L$.  Then the semigroup $(e^{-tL})_{t>0}$ has an analytic extension to $\Sigma_{{\frac \pi 2}-\omega}$ on $L^p$ for $p\in \mathrm{int}\calJ(L)$.
Moreover, for 
$p\in \mathrm{int}\calJ(L)$ and all $\beta \in (0, {\frac \pi 2}-\omega)$, the families $(e^{-zL})_{z\in \Sigma_\beta}$ and $(zLe^{-zL})_{z\in \Sigma_\beta}$ are $L^p$
bounded, satisfy
$L^p-L^2$ (resp. $L^2-L^p$) off-diagonal estimates and are $L^p-L^2$ (resp. $L^2-L^p$) bounded if $p<2$ (resp. if $p>2$). 

\end{proposition}

\begin{proof} The analyticity  is a consequence of  Stein's complex interpolation theorem together with holomorphy of the semigroup on $L^2$. 
However, this does not yield the best angle for the sector of holomorphy. Here is a better argument. 

 For  $\alpha \in (-{\frac \pi 2}+\omega,{\frac \pi 2}-\omega)$,  set $L_\alpha=e^{i\alpha}L$. It is an operator in the same class as $L$, associated to the matrix
of coefficients
$e^{i\alpha}A(x)$. Hence, Proposition \ref{propequiv} applies to the semigroup associated to $L_\alpha$. Furthermore,
a careful check of the proof shows that the various constants are independent of $\alpha$ as long as $\alpha $ is restricted to 
a compact subset of $(-{\frac \pi 2}+\omega,{\frac \pi 2}-\omega)$. Let $\beta\in (0,{\frac \pi 2}-\omega)$. If  $z=te^{i\alpha} \in
\Sigma_{\beta}$, then
$e^{-zL}=e^{-tL_\alpha}$ and the reasonning above shows that the statement of Proposition \ref{propequiv} extends to the complex family
$(e^{-zL})_{z\in \Sigma_\beta}$. 
Hence, for  $p\in \mathrm{int}\calJ(L)$, it remains to  showing that this family is $L^p-L^2$ (resp. $L^2-L^p$) bounded    if $p<2$ (resp. if $p>2$). 

The case
$p>2$ can be handled by duality, we restrict attention to $p<2$.  Let $z \in \Sigma_\beta$. Choose $\beta<\beta'<{\frac \pi 2}-\omega$. Elementary geometry shows that one
can  decompose $z=\zeta+t$ with $\zeta \in \Sigma_{\beta'}$, $t>0$ and $|z|\sim |\zeta| \sim t$ where the implicit constants only depend on $\beta$ and $\beta'$.
 By assumption and Proposition \ref{propequiv}, $(e^{-tL})_{t>0}$ is $L^p-L^2$ bounded. 
Since $(e^{-\zeta L})_{\zeta \in \Sigma_{\beta'}}$ is $L^2$ bounded, writing 
$e^{-zL}= e^{-\zeta L}e^{-tL}$, this shows that $(e^{-zL})_{z\in \Sigma_\beta}$ is $L^p-L^2$ bounded. 
\end{proof}

\begin{proposition}\label{propequivgrafdz} Let $\omega$ be the type of $L$.  Then, for 
$p\in \mathrm{int}\calN(L)$ and all $\beta \in (0, {\frac \pi 2}-\omega)$, its analytic extension  $(\sqrt z\, \nabla e^{-zL})_{z\in \Sigma_\beta}$ is
$L^p$ bounded, satisfies
$L^p-L^2$ (resp. $L^2-L^p$) off-diagonal estimates and is $L^p-L^2$ (resp. $L^2-L^p$) bounded if $p<2$ (resp. if $p>2$). 

\end{proposition}

It suffices to adapt the  above proof. Further details are left to the
reader.

\section{$L^p$ theory for square roots}\label{sec:lptheorysquareroots}

Let $L$ be as the introduction. We study here the following sets:

\begin{enumerate}
\item the maximal interval of exponents $p$ in $(1,\infty)$ for which 
one has the $L^p$ boundedness of the Riesz transform $\nabla L^{-1/2}$, which
we call  $\calI(L)$. 

\item
 the maximal interval
 of exponents $p$ in $(1,\infty)$ for which 
one has the \textit{a priori} inequality  $\| L^{1/2}f\|_{p}\lesssim
\| \nabla f\|_{p}
$ for $f \in  C^\infty_0$, which we call $\calR(L)$.

\end{enumerate}

We  characterize the limits of $\calI(L)$ and obtain bounds on the limits of
$\calR(L)$.

\subsection{Riesz transforms on $L^p$}\label{sec:RT}
We prove here the following theorem. 
 
\begin{theorem} The interior of $\calI(L)$ equals $(p_-(L), q_+(L))$. 
\end{theorem}

Recall that $p_-(L)$ is the lower limit of both $\calJ(L)$ and $\calN(L)$, and
that $q_+(L)$ is the upper limit of $\calN(L)$. The cases $p<2$ and $p>2$ are
treated by different methods. 

\subsubsection{The case $p<2$}\label{sec:RTp<2}

  We introduce the following sets.
\begin{align*} 
{\cal I}_-(L)&= \calI(L) \cap (1,2)
\\
{\cal J}_-(L)&= \calJ(L) \cap [1,2)
\\
{\cal K}_-(L)&= \{1\le p< 2;  (e^{-tL})_{t>0}
\ \mathrm{ is }\  L^p-L^2 \ \mathrm{ bounded} \}
\\
{\cal M}_-(L)&= \{1\le p< 2; (e^{-tL})_{t>0} \ \mathrm { satisfies }\  L^p-L^2 \ \mathrm{ off-diagonal\ estimates } \}
\end{align*}

  It is quite clear from 
interpolation that these sets are intervals with 2 as upper limit, provided they are nonempty. 
It follows directly from Proposition \ref{propequiv} that ${\cal J}_-(L)$, ${\cal K}_-(L)$ and ${\cal M}_-(L)$ have the same interiors
(that is, the same lower limit).
Also, by Corollary \ref{corohyper} the interior of ${\cal J}_-(L)$ is not empty.

\begin{theorem}\label{thrieszp<2}\footnote{\,
 The fact that off-diagonal estimates implies boundedness of the Riesz transforms
is proved in \cite{BK3} and this idea is also the main tool in \cite{HM}. The
converse was not noticed in these works.} 
 All the above sets are intervals with common interiors.
\end{theorem}

Let us derive a corollary for results in the range $1<p<2$.

\begin{proposition}\label{propRTallp<2} $\nabla L^{-1/2}$ is bounded on $L^p$ for $1<p<2$ if, and only if, 
$(e^{-tL})_{t>0}$ is $L^p-L^2$ bounded for $1<p<2$. In particular,
 if $(e^{-tL})_{t>0}$ is $L^1-L^2$ bounded, then $\nabla L^{-1/2}$ is bounded on $L^p$ for $1<p<2$. 
\end{proposition}

\begin{proof} The equivalence is contained in the theorem above since $int{\calI}_-(L)=int{\calK}_-(L)$. 
Next, if
$(e^{-tL})_{t>0}$ is
$L^1-L^2$ bounded, then by interpolation it is also $L^p-L^2$ bounded for $1<p<2$
\end{proof}

\begin{remark} Usually, the $L^p$ boundedness of the Riesz transform for $1<p<2$ is obtained under the
stronger assumption that the semigroup kernel has good pointwise upper bounds.\footnote{\, In \cite[Chapter IV]{AT},
 this is obtained via an $\calH^1$ estimate under $L^1-L^\infty$ off diagonal estimates plus
H\"older regularity on the heat kernel. The regularity was removed in \cite{DMc1} to obtain weak type
(1,1). In
\cite{CD}, the use of weighted $L^2$ estimates similar to the $L^1-L^2$ boundedness was stressed in the
context of the Riesz transform for the Laplace-Beltrami on Riemannian manifolds.} 
\end{remark}

Let us prove Theorem \ref{thrieszp<2}. According to the above remark, the following  steps suffice to prove
this result:

\paragraph{Step 1. $p \in {\calI}_-(L)$ implies $p \in {\calK}_-(L)$.}

\paragraph{Step 2. $p_0 \in {\calM}_-(L)$ and $p_0<p<2$ imply $p \in {\calI}_-(L)$.}

\paragraph{Step 1. $p \in {\cal I}_-(L)$ implies $p \in {\cal K}_-(L)$.}

First the conclusion is true if $n\le 2$ from Corollary \ref{corolpbdd}. We assume, therefore, that $n\ge 3$. 

Let $p\in {\cal I}_-(L)$. If $p\ge 2_*=\frac{2n}{n+2}$, we already know that $(e^{-tL})_{t>0}$
 is $L^p-L^2$. Assume now that $p<2_*$. Remark that any $q\in [p,2]$ belongs to ${\calI}_-(L)$. By Sobolev embedding's, we  see
that 
$L^{-1/2}\colon L^{q} \to L^{q^*}.$
Let $p_0=p$ and $p_{k} =(p_{k-1})^*$ and stop when $k$ is the largest integer that satisfies
$p_k< 2$.  We have that  
$L^{-k/2}$ is bounded from $L^p$ into $L^{p_k}$. Write for $f \in L^2\cap L^p$,
$$
e^{-tL}f= (e^{-(t/2)L}L^{k/2}) e^{-(t/2)L} L^{-k/2}f,$$
the equality being justified by the fact that $f \in L^2$. We successively have 
that $g=L^{-k/2}f \in L^{p_k}$, then $h=e^{-(t/2)L} g\in L^2$ by Corollary \ref{corohyper} since $p_k\ge \frac{2n}{n+2}$,
and $e^{-(t/2)L}L^{k/2} h \in L^2$ by the bounded holomorphic calculus of $L$ on $L^2$.
Hence, $e^{-tL}$ is bounded from $L^p$ into $L^2$ by density and
we obtain the right bound  for its operator norm by keeping track of the bounds from each step.

\paragraph{Step 2. $p_0 \in {\calM}_-(L)$ and $p_0<p<2$ imply $p \in {\calI}_-(L)$.}
We   apply Theorem \ref{lpp<2} to the operator $T=\nabla L^{-1/2}$ to obtain weak type $(p_0,p_0)$.  We first introduce 
the operators $A_r=I- (I-e^{-r^2L})^m$ where $m$ is some integer to be specified later.

Observe that $A_r =\sum_{k=1}^m c_k e^{-kr^2L}$ for some numbers $c_k$. A direct consequence of the assumption
that $L$ satisfies $L^{p_0}-L^2$ off-diagonal estimates is that 
if $B$ is any ball, with radius $r$, and $f \in L^2\cap 
L^{p_0}$  with support in  $B$, and $j\ge 1$, then
$$
\left(\int_{C_j(B)}| A_r f|^2\right)^{1/2} \le C r^{-\gamma} e^{-\alpha 4^j}  \left(\int_{B}|  f|^{p_0}\right)^{1/p_0}
$$
where $\gamma=\gamma_{p_0}= |\frac n 2 - \frac n {p_0}|.$ Recall that $C_j(B)$ denotes the ring $2^{j+1}B\setminus 2^jB$ if $j\ge 2$
and $C_1(B)=4B$. Hence the assumption
\eqref{domap<2} in Theorem  \ref{lpp<2} holds with $g(j)=C(m)2^{-nj/2} e^{-\alpha 4^j}$ for any $m\ge 1$. It remains to check the assumption 
\eqref{pointwisep<2} and this is where we use the role of $m$. 

\begin{lemma}\label{lemmagradient}  Let $p_0\in {\calM}_-(L)$.
There exists $C\ge 0,$ such that  for all balls $B$ with radius $r>0$   and  $f \in L^2\cap 
L^{p_0}$  with support in  $B$, and $j\ge 2$
  \begin{equation}\label{eqgradient}
\|\nabla L^{-1/2}(I-e^{-r^2L})^mf\|_{L^2(C_j(B))} \le Cr^{-\gamma} 2^{-j(2m+ \gamma )} \| f\|_{L^{p_0}(B)},
\end{equation}
where $\gamma=\gamma_{p_0}= |\frac n 2 - \frac n {p_0}|.$
\end{lemma}  

\begin{proof}     By expanding $(I-e^{-r^2L})^m$ in  the representation of the square root,
we obtain
\begin{eqnarray*}\nabla L^{-1/2}
(I-e^{-r^2L})^mf&=&\pi^{-1/2}\int_0^\infty  \nabla e^{-tL} (I-e^{-r^2L})^m f\,\frac{dt}{\sqrt
t}\\ &=&\pi^{-1/2}\int_0^\infty  g_{r^2}(t) \nabla e^{-tL} f \, {dt}\\
\end{eqnarray*}
where using the usual notation for the binomial coefficient, 
$$g_s(t)=\sum_{k=0}^m \binom m k (-1)^k \frac{\chi_{\{t>ks\}}}
{\sqrt{t-ks}},$$ and $\chi$ is the indicator function of $(0,\infty)$. 

Observe that $(\sqrt t\,\nabla e^{-tL})_{t>0}$ satisfies $L^{p_0}-L^2$ off-diagonal estimates by applying the composition
to the families $(\sqrt t\,\nabla e^{-(t/2)L})_{t>0}$ and $(e^{-(t/2)L})_{t>0}$ which satisfy respectively
$L^2$ and $L^{p_0}-L^2$ off-diagonal estimates. By Minkow\-ski integral inequality,  we
have that 
$$
\| \nabla L^{-1/2}
(I-e^{-r^2L})^mf\|_{L^2(C_j(B))} \le C \int_0^\infty |g_{r^2}(t)| e^{-\frac
{\alpha 4^jr^2} t} t^{-\gamma/2}\, \frac {dt} {\sqrt t} \ \| f\|_{L^{p_0}(B)}.
$$ 
The latter integral can be estimated as follows.  
Elementary analysis yields the following estimates for $g_{r^2}$:
$$
|g_{r^2}(t)| \le \frac C{\sqrt {t-kr^2}} \quad \textrm{if} \quad {kr^2<t \le (k+1)r^2 \le (m+1)r^2}
$$
and 
$$
|g_{r^2}(t)| \le Cr^{2m}t^{-m-\frac{1}{2}} \quad \textrm{if} \quad{t> (m+1)r^2}.
$$
The latter estimate comes from the inequality 
$$\left|\sum_{k=0}^m \binom m k (-1)^k v(t-ks) \right| \le C s^m \sup_{u\ge t/(m+1)}| v^{(m)}(u)| 
$$
for $t >  (m+1)s\ge 0$ after expanding $v(t-ks)$ using Taylor's formula about $t$ and using the
classical  relations $\sum_{k=0}^m \binom m k (-1)^k k^\ell =0$ for $\ell \in
\NN$, $\ell < m$.   This readily yields the estimates 
$$
\int_0^\infty |g_{r^2}(t)| e^{-\frac
{\alpha 4^ir^2} t}t^{-\gamma/2}\,  \frac {dt} {\sqrt t} \le C r^{-\gamma}2^{-j(2m+\gamma)}.
$$

\subparagraph{Alternate proof of the lemma:} By   the representation of the square root,
one obtains
\begin{equation}\label{eqRTregularise}
\nabla L^{-1/2}
(I-e^{-r^2L})^mf=\pi^{-1/2}\int_0^\infty  \sqrt t\, \nabla e^{-tL} (I-e^{-r^2L})^m f\,\frac{dt}{
t}.
\end{equation}
The  function $\varphi(z)=e^{-tz}(1-e^{r^2z})^m$ is holomorphic and  satisfies the technical condition \eqref{eqtechnical}
in any sector $\Sigma_\mu$, $\mu<{\frac \pi 2}$. Hence, one may use the representation \eqref{eqrepresentation} to compute $\varphi(L)$. With the same choices
of the parameters as in Section \ref{sec:fcL2},  for some positive constant $c$, the functions $\eta_\pm$ given by \eqref{eqeta} satisfy
$$|\eta_\pm(z)| \le \int_{\gamma_\pm} e^{-c|\zeta|(|z|+t)} |1-e^{-r^2\zeta}|^m\, |d\zeta|,
$$
and one finds
\begin{equation}\label{eqetaRT}
|\eta_\pm(z)| \le \frac {C}{|z|+t} \inf \left(1 , \frac {r^{2m}}{(|z|+t)^m}\right), \quad z\in \Gamma_\pm.
\end{equation}
Observe that for any $0<\beta< \beta' <{\frac \pi 2}-\omega$,  $(\sqrt {z}\,\nabla e^{-zL})_{z \in \Sigma_{\beta}}$ satisfies $L^{p_0}-L^2$ off-diagonal estimates by
applying the composition lemma to the families $(\sqrt {z}\,\nabla e^{-zL})_{z \in \Sigma_{\beta'}}$ and $(e^{-tL})_{t>0}$ which satisfy respectively
$L^2$ and $L^{p_0}-L^2$ off-diagonal estimates. It suffices to decompose $ z \in \Sigma_{\beta}$ into  $z'+t$ with $z' \in \Sigma_{\beta'}$ and $t>0$
with $|z|\sim |z'|\sim t$. Using this in \eqref{eqrepresentation} and the estimate for $\eta_\pm$,  $\| \sqrt t\, \nabla e^{-tL}
(I-e^{-r^2L})^mf\|_{L^2(C_j(B))}$ is bounded by 
$$
  C \int_{\Gamma_+} e^{-\tfrac{c4^jr^2} {|z|}}\frac { t^{1/2}}{|z|^{1/2}} \frac 1 {|z|^{\gamma/2}} \frac {1}
{(|z|+t)} \frac {r^{2m}}{(|z|+t)^{m}}\, |dz| \ \| f\|_{L^{p_0}(B)}
$$
plus the similar term corresponding to integration  on $\Gamma_-$. Here $\gamma=\gamma_{p_0}=|\frac n 2 -\frac n {p_0}|$.  Using the inequality \eqref{eqintegrale} below for the integral, this
gives us the bound
$$
  \frac{C}{4^{jm}(2^jr)^{\gamma}} \inf \left( \left(\frac {t} {4^jr^2}\right)^{1/2}, 
 \left(\frac  {4^jr^2}{t}\right)^{m-1/2}\,\, \right)\ \| f\|_{L^{p_0}(B)}.
$$
Integrating with respect to  $t$ in \eqref{eqRTregularise}, one finds \eqref{eqgradient} by Minkowski inequality
 since $m\ge 1>1/2$. 
\end{proof}

Thus, \eqref{pointwisep<2} holds with $g(j)= 2^{-nj/2}2^{-j(2m+\gamma)}$. Hence, 
the summability condition $\sum_{j\ge 1} 2^{nj} g(j)<\infty$ is granted if
$2m+\gamma>n/2$. This finishes the proof modulo the proof of \eqref{eqintegrale} which we do in the following lemma.

\begin{lemma}
Let $\gamma \ge 0$, $\alpha\ge 0$, $m> 0$ be fixed parameters, and $c$ a positive constant.  For  some $C$ independent of  $j\in
\NN$,
$r>0$ and
$t>0$,  the integral
$$
I= \int_0^\infty  e^{-\tfrac{c4^jr^2} {s}} \frac 1 {s^{\gamma/2}} \frac {t^\alpha}
{(s+t)^{1+\alpha}} \frac {r^{2m}}{(s+t)^{m}}\, ds
$$
satisfies   the estimate
 \begin{equation}\label{eqintegrale}
I \le  \frac{C}{4^{jm}(2^jr)^\gamma} \inf \left( \left(\frac {t} {4^jr^2}\right)^{\alpha}, 
 \left(\frac  {4^jr^2}{t}\right)^{m}\,\, \right).
\end{equation}
\end{lemma} 

\begin{proof} Set $\beta= \frac  {4^jr^2}{t}$ for the argument. 
Split $I=I_1+I_2$ where $I_1=\int_t^\infty \ldots$ and $I_2=\int_0^t \ldots$.

If $t\le s$, write $\frac {t^\alpha}
{(s+t)^{1+\alpha}} \le \frac {t^\alpha}
{s^{1+\alpha}}$, $\frac {r^{2m}}{(s+t)^{m}} \le \frac {r^{2m}}{s^{m}}$ and change variable by setting
$\frac{4^jr^2} {s}=u$. Then,
$$
I_1\le  \int_0^{\frac{4^jr^2} {t}} e^{-cu} \left( \frac {u}{4^jr^2}\right)^{\gamma/2}\left( \frac
{tu}{4^jr^2}\right)^{\alpha} \frac {u^m}{4^{jm}}\, \frac {du} u =  \frac{\beta^{-\alpha}}{4^{jm}(2^jr)^\gamma } \int_0^\beta e^{-cu} u^{\gamma/2 + \alpha + m} \,  \frac {du} u.
$$
As $\gamma/2 + \alpha +m>0$, we obtain
$$
I_1 \le   \frac{C}{4^{jm}(2^jr)^\gamma}  \inf \big( \beta^{\gamma/2 + m}, \beta^{-\alpha}\big). 
$$

If $t\ge s$ and $\gamma >0$, write 
$\frac {t^\alpha}
{(s+t)^{1+\alpha}} \le \frac {1}
{s}$, $\frac {r^{2m}}{(s+t)^{m}} \le \frac {r^{2m}}{t^{m}}$ and 
change variable by setting
$\frac{4^jr^2} {s}=u$. Then,
$$
I_2\le  \int_{\frac{4^jr^2} {t}}^\infty e^{-cu} \left( \frac {u}{4^jr^2}\right)^{\gamma/2} \frac
{r^{2m}}{t^{m}}\, \frac {du} u = \frac{\beta^m}{4^{jm}(2^jr)^\gamma} \int_{\beta}^\infty e^{-cu} 
u^{\gamma/2}\,  \frac {du} u.
$$
Since $\gamma>0$, we obtain  
$$
I_2 \le  \frac{C}{4^{jm}(2^jr)^\gamma}  \inf\big( \beta^m, e^{-c\beta}\big)$$
where the value of $c$ has changed. The conclusion follows readily. 

To treat the case  $\gamma=0$  when $\beta \le 1$, then use instead 
$\frac {t^\alpha}
{(s+t)^{1+\alpha}} \le \frac {1}
{t}$ and that  $\int_\ep^\infty e^{-cu} \frac {du}{u^2}$ is asymptotic to $c\ep^{-1}$ as $\ep $ tends to
$0$. Further details are left to the reader.
\end{proof}

\subsubsection{The case $p>2$}\label{sec:RTp>2}

 Introduce the following  sets:
\begin{align*} 
{\cal I}_+(L)&= \calI(L) \cap (2,\infty)
\\
{\calN}_+(L)&= \calN(L) \cap (2,\infty]
\\
{\calK}_+(L)&= \{2< p\le\infty; (\sqrt t\, \nabla e^{-tL})_{t>0}\  \mathrm{ is }\ L^2-L^p \ \mathrm{ bounded } \}
\\
{\calM}_+(L)&= \{2< p \le\infty; (\sqrt t\, \nabla e^{-tL})_{t>0}\  \mathrm{ satisfies }\ L^2-L^p\ \mathrm{off-diagonal\ estimates } \}
\end{align*}
 
Observe that these sets, if nonempty,  are intervals with 2 as lower 
limit. We know that $\calN_+(L)$ is not empty and that 
it has the same interior as $\calK_+(L)$ and $\calM_+(L)$.

\begin{theorem}\label{thrieszp>2} The above sets are intervals with common interiors.
\end{theorem}

Let us  state a result for the range $2<p<\infty$.

\begin{proposition}\label{propRTallp>2} $\nabla L^{-1/2}$ is bounded on $L^p$ for $2<p<\infty$
if, and only if,   $(\sqrt t\,\nabla e^{-tL})_{t>0}$ is $L^2-L^p$ bounded for $2<p<\infty$.
 In particular, if $(\sqrt t\,\nabla e^{-tL})_{t>0}$ is $L^2-L^\infty$ bounded, 
then $\nabla L^{-1/2}$ is bounded on $L^p$ for $2<p<\infty$.
\end{proposition}

\begin{proof} The equivalence is contained in the theorem above. 
If $(\sqrt t\,\nabla e^{-tL})_{t>0}$ is $L^2-L^\infty$ bounded then by interpolation with the
$L^2$-boundedness, it is also $L^2-L^p$ bounded for $2<p<\infty$.   
\end{proof}

Let us prove Theorem
\ref{thrieszp>2}.
In view of the above remarks, it is enough to show   

\paragraph{Step 1. $p \in {\calI}_+(L)$ implies $p \in {\calK}_+(L)$.} 

\paragraph{Step 2. $p_0\in {\calM}_+(L)$ and $2<p<p_0$ imply $p \in {\calI}_+(L)$.}

\paragraph{Step 1.} \textbf{ $p \in {\calI}_+(L)$ implies $p \in {\calK}_+(L)$.}

Assume first that $n\ge 3$. By Sobolev embedding's, 
 if $q\in [2,p]$ with $q<n$ then
$$
L^{-1/2}\colon L^q \to L^{q^*}.$$
Following the method of Section \ref{sec:RTp<2} applied to the dual semigroup, we obtain
that if $q \in \RR$ with $2\le q \le p^*$, 
$$
e^{-tL}\colon L^2 \to L^q
$$
with bound $Ct^{-\gamma_q}/2$. Since,  $L^{1/2} e^{-(t/2)L}$ is bounded on $L^2$ with bound $Ct^{-1/2}$, we obtain 
$\nabla e^{-tL}=  \nabla L^{-1/2} e^{-(t/2) L} L^{1/2} e^{-(t/2)L}$ is bounded from $L^2$ to  $L^p$ with the desired bound on the operator norm. 

If $n\le2$, we already know that $$
e^{-tL}\colon L^2 \to L^q
$$ for all $q\in [2, \infty]$ (See Section \ref{sec:Lpsemigroup}). The argument is then similar.

\paragraph{Step 2. $p_0\in {\calM}_+(L)$ and $2<p<p_0$ imply $p \in {\calI}_+(L)$:}

We apply Theorem \ref{lpp>2} to $T=\nabla L^{-1/2}$.  We set again $A_r=I -(I-e^{-r^2L})^m$ for some integral number $m$ to be chosen. 
We have to check \eqref{pointwise} and \eqref{doma}.

We begin with  the following lemma.  

\begin{lemma}\label{lempointwise}  For every ball $B$ with radius $r>0$,  
\begin{equation}\label{pointwiseriesz}
\| \nabla L^{-1/2}
(I-e^{-r^2L})^mf\|_{L^2(B)} \le   |B|^{1/2}  \sum_{j\ge 1} g(j) \left( \frac 1{|2^{j+1}B|} \int_{2^{j+1}B} |f|^2 \right)^{1/2} 
\end{equation}
with $g(j) =C 2^{(n/2)j} 4^{-mj}$.
  \end{lemma}

Hence,  (\ref{pointwise}) follows provided  $m>n/4$, since $\left( \frac 1{|2^{j+1}B|} \int_{2^{j+1}B} |f|^2 \right)^{1/2}$
is controlled by $\big({M}(|f|^2)(y)\big)^{1/2}$ for any $y\in B$. 

\begin{proof}  Fix a ball $B$ and $r=r(B)$ its radius.   Decompose $f$ as $f_1+f_2+f_3+\ldots$
with 
$f_j=f\chi_{C_j}$ where $\chi_{C_j}$ is the indicator function of $C_j=C_j(B)$. By Minkowski inequality we have that 
$$
\| \nabla L^{-1/2} (I-e^{-r^2L})^mf\|_{L^2(B)} \le \sum_{j\ge 1} 
\| \nabla L^{-1/2}
(I-e^{-r^2L})^mf_j\|_{L^2(B)}.
$$
For $j=1$ we use $L^2$ boundedness of $\nabla L^{-1/2}
(I-e^{-r^2L})^m$: 
$$
\| \nabla L^{-1/2}
(I-e^{-r^2L})^mf_1\|_{L^2(B)} \lesssim \| f \|_{L^2(4B)} \le |4B|^{1/2}
 \left( \frac 1{|4B|} \int_{4B} |f|^2 \right)^{1/2}.
$$
\subparagraph{First proof in the case $j\ge 2$:} Expanding $(I-e^{-r^2L})^m$ in  the representation of the square root, we find
\begin{eqnarray*}\nabla L^{-1/2}
(I-e^{-r^2L})^mf_j&=&\pi^{-1/2}\int_0^\infty  \nabla e^{-tL} (I-e^{-r^2L})^mf_j\,\frac{dt}{\sqrt
t}\\ &=&\pi^{-1/2}\int_0^\infty  g_{r^2}(t) \nabla e^{-tL}f_j \, {dt}\\
\end{eqnarray*}
where, as in Section \ref{sec:RTp<2}, 
$$g_s(t)=\sum_{k=0}^m \binom m k (-1)^k \frac{\chi_{\{t>ks\}}}
{\sqrt{t-ks}}.$$  

By Minkowski integral inequality and the $L^2$ off-diagonal estimates for $(\sqrt t \nabla e^{-tL})_{t>0}$ using the support of $f_j$, we
have that 
$$
\| \nabla L^{-1/2}
(I-e^{-r^2L})^mf_j\|_{L^2(B)} \le C \int_0^\infty |g_{r^2}(t)| e^{-\frac
{\alpha 4^jr^2} t}  \frac {dt} {\sqrt t} \ \| f\|_{L^2(C_j)}.
$$ 
As in Section \ref{sec:RTp<2}, the latter integral is bounded above by $C4^{-jm}$ uniformly over $r>0$.
Next,
$$\| f\|_{L^2(C_j)}  \le |2^{j+1}B|^{1/2}
 \left( \frac 1{|2^{j+1}B|} \int_{2^{j+1}B} |f|^2 \right)^{1/2}
$$
and we obtain \eqref{pointwiseriesz}.

\subparagraph{Second proof in the case $j\ge 2$:} Write again \eqref{eqRTregularise}:
$$\nabla L^{-1/2}
(I-e^{-r^2L})^mf_j=\pi^{-1/2}\int_0^\infty  \sqrt t\, \nabla e^{-tL} (I-e^{-r^2L})^mf_j\,\frac{dt}{
t}.
$$ 
As in Section \ref{sec:RTp<2},   one may use the representation \eqref{eqrepresentation} with the  function $\varphi(z)=e^{-tz}(1-e^{r^2z})^m$. The functions 
$\eta_\pm$ in \eqref{eqeta} satisfy the estimates \eqref{eqetaRT}. Since 
for any $0<\beta< {\frac \pi 2}-\omega$,  $(\sqrt z\,\nabla e^{-zL})_{z \in \Sigma_{\beta}}$ satisfies $L^{2}$ off-diagonal estimates, using 
\eqref{eqrepresentation} and the estimate for $\eta_\pm$,  $\| \sqrt t\, \nabla e^{-tL} (I-e^{-r^2L})^mf_j\|_{L^2(B)}$ is bounded by 
$$
  C \int_{\Gamma_+} e^{-\tfrac{c4^jr^2} {|z|}}\frac { t^{1/2}}{|z|^{1/2}}  \frac {1}
{(|z|+t)} \frac {r^{2m}}{(|z|+t)^{m}}\, |dz| \ \| f\|_{L^{2}(C_j)}
$$
plus the similar term corresponding to integration  on $\Gamma_-$.   Using \eqref{eqintegrale}, this gives us the bound
$$
  \frac{C}{4^{jm}}  \inf \left( \left(\frac {t} {4^jr^2}\right)^{1/2}, 
 \left(\frac  {4^jr^2}{t}\right)^{m-1/2}\,\, \right)  \ \| f\|_{L^{2}(C_j)}.
$$
Integrating with respect to  $t$ in \eqref{eqRTregularise}, one finds \eqref{pointwiseriesz} by Minkowski inequality
 since $m\ge 1>1/2$. 
\end{proof}

We now show that (\ref{doma}) holds. By expanding $A_r=I-(I-e^{-r^2L})^m$ it suffices to show
\begin{equation}\label{maxp}
\left(\frac{1}{|B|}\int_B |\nabla e^{-kr^2L}f|^{p_0} \right)^{1/p_0}\le  \sum_{j\ge 1} g(j) \left( \frac 1{|2^{j+1}B|} \int_{2^{j+1}B} |\nabla f|^2 \right)^{1/2}
\end{equation}
for every ball $B$  with $r=r(B)$,  and 
$k=1,2, \ldots, m$ with $\sum g(j) <\infty$. Recall that $m$ is  chosen larger than $n/4$. This, applied to $f=L^{-1/2}g$ for appropriate $g$,
gives us \eqref{doma}.

Using the conservation property \eqref{conservation} of the semigroup, we have 
$$
\nabla e^{-kr^2L}f = \nabla e^{-kr^2L}(f- f_{4B})$$
where $f_{E}$ is the mean of $f$ over $E$.
Write $f- f_{4B}= f_1+f_2+f_3+\ldots$ where 
 $f_j= (f- f_{4B})\chi_{C_j(B)}$.
For $j=1$, we use the fact that $p_0 \in {\cal K}_+$, that is that
$(\sqrt t\, \nabla e^{-tL})_{t>0}$ satisfies $L^2-L^{p_0}$ off-diagonal estimates, and Poincar\'e inequalities to obtain
\begin{align*}
\left( \frac 1{|B|}\int_B   |\nabla e^{-kr^2L}f_1|^{p_0} \right)^{1/p_0} 
&\le C\left(\frac {1} {r^2|4B|} \int_{4B}
|f-f_{4B}|^2  \right)^{1/2}
\\
&
\le C\left(\frac {1} {|4B|} \int_{4B}
|\nabla f|^2  \right)^{1/2}
\end{align*}
For $j\ge 2$, we have similarly by the assumption on $p_0$ and the various support assumptions,
$$
\left( \frac 1{|B|}\int_B   |\nabla e^{-kr^2L}f_{j}|^{p_0} \right)^{1/p_0} 
\le \frac {Ce^{- \alpha 4^j}2^{jn/2}}{r} \left(\frac {1} {|2^{j+1}B|} \int_{C_j(B)}
|f_j|^2 \right)^{1/2}.
$$
But 
$$
\int_{C_j(B)}
|f_j|^2 \le \int_{2^{j+1}B}
|f - f_{4B}|^2,
$$
$$
|f - f_{4B}| \le |f-f_{2^{j+1}B}|
+\sum_{\ell=2}^{j}|f_{2^\ell B}-f_{2^{\ell+1}B}|
$$
and observe that by Poincar\'e inequality,
$$
|f_{2^\ell B}-f_{2^{\ell+1}B}|^2  \le \frac {2^n} {|2^{\ell+1}B|} \int_{2^{\ell+1}B}
|f - f_{2^{\ell+1}B}|^2 \le   \frac {C(2^{\ell}r)^2 } {|2^{\ell+1}B|} \int_{2^{\ell+1}B}
|\nabla f|^2.
$$
Hence, by Minkowski inequality, we easily obtain
$$
\left( \frac 1{|B|}\int_B   |\nabla e^{-kr^2L}f_j|^{p_0} \right)^{1/p_0} 
\le  {Ce^{- \alpha 4^j}} 2^{jn/2} \sum_{\ell=1}^j  {2^{\ell} } \left(\frac 1 {|2^{\ell+1}B|} \int_{2^{\ell+1}B}
|\nabla f|^2\right)^{1/2}
$$
and summing over $j\ge 2$  gives us \eqref{maxp}.

This  concludes the proof of Theorem \ref{thrieszp>2}.

\subsection{Reverse inequalities} 

 In this section, we study the reverse inequality to the Riesz transform $L^p$ boundedness.
 Recall  the maximal interval,
$\calR(L)$, of exponents $p$ in $(1,\infty)$ for which 
one has the \textit{a priori} inequality  $\| L^{1/2}f\|_{p}\lesssim
\| \nabla f\|_{p}
$ for $f \in  C^\infty_0$
and write $\mathrm{int}{\calR(L)}= (r_-(L),r_+(L))$. We know so far that $\mathrm{int}\calI(L)=(p_-(L),q_+(L))$ and that
$p_+(L) \ge (q_+(L))^*$.
We show the following bounds. 

\begin{theorem}  \label{riesztoreverse} We have
\begin{align*}
r_-(L) &\le \sup (1, (p_-(L))_*),\\
r_+(L) &\ge p_+(L). 
\end{align*}
Hence,  ${\calR}(L)$ contains a neighborhood of the closure (in $\RR$) of ${\calI}(L)$.
\end{theorem}

We discuss the optimality of these bounds in Section \ref{sec:more}.

This result is a consequence of  Lemma \ref{lem:duality}  and
Proposition \ref{riesztoreversep<2} below.
We begin with a duality principle which applies for all $p$ in $(1,\infty)$ but
which  gives us the bound in Theorem \ref{riesztoreverse} only
for $r_+(L)$.

\begin{lemma}\label{lem:duality}  If $1<p<\infty$ and $\nabla L^{-1/2}$ is bounded on $L^p$ then   $\| {(L^*)}^{1/2}f\|_{p'}\lesssim
\| \nabla f\|_{p'}$ holds for $f \in C_0^\infty $. Hence, $r_+(L)\ge p_+(L)$.
\end{lemma}

\begin{proof} 
Let $f \in C_0^\infty $. Then $ {(L^*)}^{1/2}f$ is defined and belongs to $L^2$. We estimate $ \| {(L^*)}^{1/2}f\|_{p'}$ by testing
against
$g
\in L^2
\cap L^{p}$.  Since
$g
\in L^2$, we have
$g=L^{1/2}L^{-1/2}g$  and $h=  L^{-1/2}g \in \dot W^{1,2}$.
 Hence, by Lemma \ref{lem:rep},
$$\langle {(L^*)}^{1/2}f, g \rangle = \int_{\RR^n} \nabla f \cdot \ol{ A\nabla h }.
$$
Since $\|\nabla h\|_p \lesssim \|g\|_p$, it follows that $\| {(L^*)}^{1/2}f\|_{p'} \lesssim \| \nabla f\|_{p'}$
as desired.

Let $2<p<p_+(L)$. By duality,  $p_-(L^*)<p'<2$, hence, $\nabla (L^*)^{-1/2}$ is
bounded on $L^{p'}$ by Theorem \ref{thrieszp<2}  and $p \in \calR(L)$. Thus,
$r_+(L)>p$ and the conclusion follows.
\end{proof}

The above lemma does not give an interesting information for $r_-(L)$.
The next proposition yields a much better bound. Define $\tilde p_-(L)=\sup (1, (p_-(L))_*)$.

\begin{proposition}\label{riesztoreversep<2} Let $1<p<2$.  If $\nabla L^{-1/2}$ is bounded on $L^p$ and $p\ge \frac{n}{n-1}$  then  $\| L^{1/2}f\|_{q}\lesssim
\| \nabla f\|_{q}$ for $p_*<q<2$.
 If $\nabla L^{-1/2}$ is bounded on $L^p$ and $p< \frac{n}{n-1}$  then $\| L^{1/2}f\|_{q}\lesssim
\| \nabla f\|_{q}$ for $1<q<2$. Hence, $r_-(L) \le \tilde p_-(L)$.
\end{proposition}

The inequality on $r_-(L)$ follows from the definition of $r_-(L)$ and the identification of $p_-(L)$ as the lower limit of $\calI(L)$.

Let us remark that for $n\ge 2$ and $p<2$, 
there exists an $L$ such that $p \in {\calI}_-(L)$ but $p' \notin
{\calI}_+(L^*)$.~\footnote{\, This is due to Kenig. See
\cite[Chapter IV, Theorem 7]{AT}.} Hence, no duality argument can help us here. 
In fact, the argument will rely on a Calder\'on-Zygmund decomposition for Sobolev
functions which is proved in the Appendix.

\begin{lemma}\label{lemmaCZD} Let $n\ge 1$,  $1\le p\le \infty$ and $f\in \calD'(\RR^n)$ be such that 
$\|\nabla f\|_p <\infty$. Let $\alpha>0$. Then, one can find a collection of cubes $(Q_i)$, functions $g$ and $b_i$  such that 
\begin{equation}\label{eqcsds1}
f= g+\sum_i b_i  \end{equation}
and the following properties hold:
\begin{equation}\label{eqcsds2}
\|\nabla g\|_\infty \le C\alpha, \end{equation}
\begin{equation}\label{eqcsds3}
b_i \in W_0^{1,p}(Q_i)\ \text{and} \ \int_{Q_i} |\nabla b_i|^p \le C\alpha^p |Q_i|, \end{equation}
 \begin{equation}\label{eqcsds4}
\sum_i |Q_i| \le C\alpha^{-p} \int_{\RR^n} |\nabla f|^p , \end{equation}
\begin{equation}\label{eqcsds5}
\sum_i {\bf 1}_{Q_i} \le N, \end{equation}
where $C$ and
 $N$ depend only on dimension and $p$.
\end{lemma}

 The space $W^{1,p}_0(\Omega)$ denotes the closure of $C_0^\infty(\Omega)$ in $W^{1,p}(\Omega)$. The point is
in the fact that the functions
$b_i$ are supported in cubes as the original  Cal\-de\-r\'on-Zygmund decomposition applied to $\nabla f$ would not give this.

\paragraph{Proof of Proposition \ref{riesztoreversep<2}:} By Theorem \ref{thrieszp<2}
 we may transform the hypothesis on the Riesz transform, \textit{i.e.} 
$p \in {\calI}_-(L)$, into an hypothesis on the semigroup. 
Proposition \ref{riesztoreversep<2}  is 
therefore a consequence of  the next result combined with Marcinkiewicz
interpolation. 

\begin{lemma}\label{thsrho}  Let  $\rho \in {\calM}_-(L)$. Then we have
\begin{align}\label{eq20ter}
\| L^{1/2}f\|_{p,\infty}  &\lesssim \| \nabla f\|_{p},  \quad \text{ if } 1 \le \rho_*<p<2, 
\\
\| L^{1/2}f\|_{1,\infty}  &\lesssim \| \nabla f\|_{1},  \quad \text{ if } \rho_*<1.
 \end{align}
\end{lemma}

\begin{proof} Recall that $\rho_*=\frac {n\rho}{n+\rho}$. Of course, it is enough to pick $\rho$ and $p$ as small as possible.  If $n=1,2$, then $\rho$ may be
chosen with $\rho_*<1$ (as a consequence of Corollary \ref{corolpbdd}) and we set $p=1$ in the proof.  If $n\ge 3$, then one can always assume that $\rho<p_n=
\frac {2n}{n+2}$ and we pick $p$ so that 
$\rho_*< p<p_n$.
Let $f
\in C_0^\infty$. We have to establish the  weak type  estimate
\begin{equation}\label{eq21}
|\{x \in \RR^n; |L^{1/2}f(x)| >\alpha\}| \le 
\frac{C}{\alpha^p}\int |\nabla f|^p,
\end{equation}
for all $\alpha>0$.  We use the following resolution of $L^{1/2}$:
$$
L^{1/2}f= c\int_0^\infty e^{-t^2 L} Lf \, dt
$$
where $c=2\pi^{-1/2}$ is forgotten from now on. It suffices to obtain the 
result for the truncated integrals $\int_\ep^R\ldots$ with bounds independent of $\ep,R$, 
and then to let
$\ep\downarrow 0$ and $R\uparrow \infty$. 
For the truncated integrals, all the calculations are justified. We ignore
this step and argue directly on $L^{1/2}$. 
Apply the Cal\-de\-r\'on-Zygmund decomposition of Lemma \ref{lemmaCZD} 
to $f$ at height $\alpha^p$ and write $f=g+\sum_i b_i$.
 By construction, $\|\nabla g\|_p \le c\|\nabla f\|_p$. Interpolating with  
 \eqref{eqcsds2} yields $\int |\nabla g|^2 \le c\alpha^{2-p}
 \int |\nabla f|^p$.
Hence
 \begin{equation*}\label{eq22}
\bigg|\bigg\{x \in \RR^n; |L^{1/2}g(x)| >\frac \alpha 3 \bigg\}\bigg| \le 
\frac{C}{\alpha^2}\int
|L^{1/2} g|^2 \le \frac{C}{\alpha^2}\int |\nabla g|^2 \le 
\frac{C}{\alpha^p}\int |\nabla f|^p
\end{equation*}
where we used the $L^2$-estimate \eqref{eq5} for square roots. 
To compute $L^{1/2}b_i$, let $r_i=2^k$ if $2^k \le \ell_i= \ell(Q_i) < 2^{k+1}$ and
set 
$T_i= \int_0^{r_i} e^{-t^2 L} L \, dt$ and
 $U_i= \int_{r_i}^\infty e^{-t^2 L} L \, dt$. It is enough to 
estimate $A=|\{x \in \RR^n; |\sum_i T_ib_i(x)| >\alpha/3\}|$ and 
$B=|\{x \in \RR^n; |\sum_i U_ib_i(x)| >\alpha/3\}|$. 
Let us  bound the first term.  

First,  
$$A \le |\cup_i 4Q_i| + \bigg|\bigg\{x \in \RR^n \setminus \cup_i  4Q_i ; 
\bigg|\sum_i T_ib_i(x)\bigg| >\frac \alpha 3\bigg\}\bigg|,$$
and by \eqref{eqcsds4}, $|\cup_i 4Q_i| \le \frac{C}{\alpha^p}\int |\nabla f|^p$.

For the other term, we have
$$
\bigg|\bigg\{x \in \RR^n \setminus \cup_i 4Q_i ; \bigg|\sum_i
T_ib_i(x)\bigg| >\frac \alpha 3\bigg\}\bigg|
\le
\frac{C}{\alpha^2}\int \bigg|\sum_i h_i\bigg|^2$$ with $h_i = {\bf
1}_{(4Q_i)^c}|T_ib_i|$. To estimate the $L^2$ norm, we dualize against
$u\in L^{2}$ with $\|u\|_{2}=1$:
$$
\int| u| \sum_i h_i = \sum_i\sum_{j=2}^\infty A_{ij}
$$  
where 
$$
A_{ij}= \int_{C_j(Q_i)} |T_ib_i||u|.
$$   
Let $q=2$ if $n\le 2$ and $q= p^*= \frac{np}{n-p}$ if $n\ge 3$. Observe that 
 $\rho< q\le  2$ and that $p\le q\le p^*$.

Since $\rho<q\le 2$, the family 
 $(tLe^{-tL})_{t>0}$ satisfies $L^q-L^2$ off-diagonal estimates on combining Proposition \ref{propequiv} and  Proposition \ref{propequivz}.
 Hence, using also $r_i \sim \ell_i$,
$$
\|e^{-t^2 L} L b_i \|_{L^2(C_j(Q_i))} \le \frac{C}{t^{\gamma+2}} e^{-\frac{c 4^jr_i^2}{t^2}}\ \|b_i\|_q
$$
where $\gamma=\gamma_q=|\frac n 2 - \frac n q|$.  By Poincar\'e-Sobolev inequality (since $p\le q \le p^*$) and \eqref{eqcsds3}, 
$$\|b_i\|_q \le c \ell_i^{1-(\frac n p -\frac n q)}
\|\nabla b_i
\|_p\le c
\alpha
\ell_i^{1+
\frac n q}.$$
 Hence,
by Minkowski integral inequality, for some appropriate  positive constants $C,c$,
\begin{align*}\| T_ib_i \|_{L^2(C_j(Q_i))} &\le   \int_0^{r_i} \|e^{-t^2 L} L b_i \|_{L^2(C_j(Q_i))} \, dt \\
&\le C\alpha  e^{-c4^j} \ell_i^{\frac n  2},
\end{align*}

Now remark that for any $y \in Q_i$ and any $j\ge 2$, 
$$
\left( \int_{C_j(Q_i)} |u|^{2}\right)^{1/2} \le \left( \int_{2^{j+1}Q_i} |u|^{2}\right)^{1/2} \le  (2^{n(j+1)}|Q_i|)^{1/2}
\big(M(|u|^{2})(y)\big)^{1/2}.
$$
Applying H\"older inequality,     one obtains
$$
A_{ij} \le  C\alpha 2^{nj/2} e^{-c4^j} \ell_i^n \big(M(|u|^{2})(y)\big)^{1/2}.
$$
Averaging over $Q_i$ yields
$$
A_{ij} \le C\alpha 2^{nj/2} e^{-c4^j} \int_{Q_i} \big(M(|u|^{2})(y)\big)^{1/2}\, dy.$$
Summing over $j\ge 2$ and $i$, we have
$$
\int| u| \sum_i h_i \le C \alpha \int \sum_i {\bf 1}_{Q_i}(y) \big(M(|u|^{2})(y)\big)^{1/2}\, dy.$$
As in the proof of Theorem \ref{lpp<2} using finite overlap \eqref{eqcsds5} of the cubes $Q_i$ and Kolmogorov's inequality, one obtains
$$
\int| u| \sum_i h_i  \le C'N\alpha \big| \cup_i Q_i \big|^{1/2} \||u|^2\|_1^{1/2}.
$$
Hence 
$$
\bigg|\bigg\{x \in \RR^n \setminus \cup_i 4Q_i ; \bigg|\sum_i T_ib_i(x)\bigg| >\frac \alpha 3\bigg\}\bigg| \le C \big| \cup_i Q_i \big| \le 
\frac{C}{\alpha^p}\int |\nabla f|^p$$
by \eqref{eqcsds5} and \eqref{eqcsds4}.

It remains to handling  the term $B$. Using functional calculus for $L$ one can
compute 
$U_i$ as $r_i^{-1}\psi(r_i^2L)$ with $\psi$ the holomorphic function on the sector $|\arg z \,| < {\frac \pi 2}$ given
by 
\begin{equation}\label{eqpsi}
\psi(z)= \int_1^\infty e^{-t^2z} z\, dt.
\end{equation}
It is easy to show that $|\psi(z)| \le C|z|^{1/2} e^{-c|z|^2}$, uniformly on subsectors $|\arg z\,  | \le \mu < {\frac \pi 2}$.

We invoke the following lemma proved by duality in Section
\ref{sec:RTsf}.~\footnote{\, It can also be proved directly using the
vector-valued extension of Theorem \ref{lpp<2}.}

\begin{lemma}\label{lemmafc} If $\rho \in {\calM}_-$  then for $\rho<q\le 2$ 
\begin{equation}\label{eq23}
\left\| \sum_{k\in \ZZ} \psi(4^kL) \beta_k \right\|_q \lesssim \left\|\left(\sum_{k\in \ZZ}  |\beta_k|^2\right)^{1/2} \right\|_q,
\end{equation}
whenever the right hand side is finite.
\end{lemma}

To apply this lemma, observe that the definitions of $r_i$ and $U_i$ yield 
$$
\sum_i U_ib_i = \sum_{k\in \ZZ} \psi(4^kL) \beta_k
$$
with 
$$\beta_k = \sum_{i, r_i=2^k} \frac{b_i}{r_i}.
$$
Using the bounded overlap property \eqref{eqcsds5}, one has that 
$$
\left\|\left(\sum_{k\in \ZZ}  |\beta_k|^2\right)^{1/2} \right\|_q^q \le C\int \sum_{i}  \frac{|b_i|^q}{r_i^q} .
$$
By \eqref{lemmasumbi}, and $p\le q \le p^*$, together with $\ell_i\sim r_i$, 
$$
\int \sum_{i}  \frac{|b_i|^q}{r_i^q}  \le  C \alpha^q  \sum_{i}|Q_i| .$$
Hence, by  \eqref{eqcsds4}
$$
\bigg|\bigg\{x \in \RR^n; \bigg|\sum_i U_ib_i(x)\bigg| 
>\frac \alpha 3\bigg\}\bigg| \le  C   \sum_{i}|Q_i| \le \frac{C}{\alpha^p}\int
|\nabla f|^p.$$
\end{proof}

\subsection{Invertibility}\label{sec:invertibility}

We finish the study of square roots with their invertibility properties on $L^p$ spaces. This is summed up in the following theorem.
Recall that ${\cal I}(L)$ is the maximal interval  of exponents  $p\in (1,\infty)$ for which $\nabla L^{-1/2}$ is bounded on $L^p$.

\begin{theorem}\label{thinvertibility} Let $1<p<\infty$. Then $p\in {\cal I}(L)$ if
and only if
$L^{1/2}$, \textit{ a priori} defined from $C_0^\infty$ into $L^2$, extends to an isomorphism from $ \dot W^{1,p}$ onto $L^p$ with 
$
 \| L^{1/2}f\|_{p}\sim
\| \nabla f\|_{p}
$.
Furthermore, ${\calI}(L)$ is an open interval.
\end{theorem}

\begin{proof} For $p=2$, the equivalence is a consequence of 
the solution of the Kato's conjecture \eqref{eq5} as recalled in Section
\ref{sec:squareroot}. 
 
For other values of $p$, one implication is
straightforward. Conversely, we know from Theorem \ref{riesztoreverse} that ${\calR}(L)$, the maximal interval  of exponents $p\in (1,\infty)$ for which one has the
\textit{ a priori} inequality 
$ \| L^{1/2}f\|_{p}\lesssim
\| \nabla f\|_{p},
$ 
contains a neighborhood of the closure (in $\RR$) of ${\cal I}(L)$.

In particular, if $\nabla L^{-1/2}$ is bounded on $L^p$ then for the same $p$, $L^{1/2}$ can be extended boundedly from $\dot W^{1,p}$ into $L^p$. The
isomorphism property is  now  easy. Indeed, from $
 \| L^{1/2}f\|_{p}\sim
\| \nabla f\|_{p}
$, we deduce that this extension is one-one and has closed range in $L^p$.  It remains to establish the density of this range.  If
$g 
\in L^p
\cap L^2$, we have
$g=L^{1/2}L^{-1/2}g$ and by
$\|\nabla L^{-1/2}g\|_p\le c\|g\|_p<\infty$, we conclude that $g$ is in
$L^{1/2}(\dot W^{1,p})$. Thus $L^{1/2}(\dot W^{1,p})$ contains a dense subspace of $L^p$. 

The openness of ${\cal I}(L)$ is  a consequence of  the following result.

\begin{lemma}\footnote{\, This is due to I. Sneiberg (\cite{Sn}).}
\label{lemmasneiberg} 
Let $X^s,Y^s$,  be  two
scales of complex interpolation Banach spaces, $s$ describing an open interval $I$. If $T\colon
X^s \to Y^s$ is bounded  for each $s \in I$, then the set
of $s\in I$ for which $T$ is an isomorphism from $X^s$ onto $Y^s$ is open.

\end{lemma}

Apply this result (changing $s$ to $p$) with $I=\mathrm{int}{\cal R}(L)$, $X^p=\dot W^{1,p}$, $Y^p= L^p$
and $T=L^{1/2}$: we just proved that the set of $p \in \mathrm{int}{\cal R}(L)$ for which $L^{1/2}$ is an isomorphism
from $\dot W^{1,p}$ onto $L^p$ is ${\cal I}(L)$. Thus, ${\cal I}(L)$ is open.
\end{proof}

\subsection{Applications}\label{sec:applications}

We have developed  the necessary theory to reobtain all the results mentioned in the Introduction. 
Let us start with results holding in all generality.

\begin{proposition}\footnote{\, This is first proved in \cite{AT1}, Theorem A.
More is proved there: in that $a\frac{d}{dx} L^{-1/2}$ is a Calder\'on-Zygmund
operator and the boundedness properties at $p=1$ and $p=\infty$ are
studied.}\label{propRTn=1} If
$n=1$, we have  $
 \| L^{1/2}f\|_{p}\sim
\| \frac{d}{dx} f\|_{p}
$ and 
 $L^{1/2}$ extends to an isomorphism from $ \dot W^{1,p}$ onto $L^p$ when
$1<p<\infty$.
\end{proposition}

\begin{proof}  We know from Corollary \ref{corolpbdd}   that $p_-(L)=1$ and  from Proposition \ref{propqplusn1}  that $q_+(L)=\infty$. 
 Hence ${\cal I}(L)=(1,\infty)$ and the conclusion follows from Theorem \ref{thinvertibility}. 
\end{proof}

\begin{proposition}\footnote{\, This follows from  \cite{AT},
Chapter IV, Theorem 1, using \cite{AMcT2}, Theorem 3.5 and
\cite{AHLMcT}, Theorem 1.4.}\label{propRTn=2} If
$n=2$, we have  $\| \nabla L^{-1/2} f\|_{p}
\lesssim
 \| f\|_{p}
$ for $1<p<2+\ep$ and $
 \| L^{1/2}f\|_{p}\lesssim
\| \nabla f\|_{p}
$  for $1<p<\infty$. Furthermore,  
 $L^{1/2}$ extends to an isomorphism from $ \dot W^{1,p}$ onto $L^p$ when
$1<p<2+\ep$.
\end{proposition}

\begin{proof} We know from Corollary \ref{corolpbdd} that   $p_-(L)=1$ and from Corollary \ref{coroqplus>2} that $q_+(L)>2$. Hence,
${\cal I}(L)=(1,2+\ep)$ and  we have the isomorphism property from Theorem \ref{thinvertibility}.

Now, $r_-(L)=1$ and  $r_+(L)\ge p_+(L)=(p_-(L^*))'=\infty $. Hence 
$ \| L^{1/2}f\|_{p}\lesssim
\| \nabla f\|_{p}
$  for $1<p<\infty$ by Theorem \ref{riesztoreverse}.
\end{proof}

\begin{proposition}\footnote{\, 
The first inequality is proved in \cite{BK3}, \cite{HM} for the
range $\frac{2n}{n+2}<p<2$ and  the extension to this larger range is in
\cite{A}, 
the second inequality is in \cite{HM} for the range $2<p<\frac{2n}{n -2}$, and in
\cite{A} for the rest of the range.}\label{propRTnge3} If $n\ge 3$, then
$\|
\nabla L^{-1/2} f\|_{p}
\lesssim
 \| f\|_{p}
$ for  $ \frac{2n}{n +2} -\ep<p<2+\ep'$ and $
 \| L^{1/2}f\|_{p}\lesssim
\| \nabla f\|_{p}
$  for $\sup(1,\frac{2n}{n+4}-\ep_1)<p<\frac{2n}{n -2} +\ep_1'$
with $\ep,\ep',\ep_1, \ep_1'>0$ depending on dimension and  the ellipticity constants of $L$. Furthermore,  
 $L^{1/2}$ extends to an isomorphism from $ \dot W^{1,p}$ onto $L^p$ when
$ \frac{2n}{n +2} -\ep<p<2+\ep'$.  
\end{proposition}

\begin{proof} We know from Corollary \ref{corolpbdd} that   $p_-(L)<\frac{2n}{n +2}$ and from Corollary \ref{coroqplus>2} that $q_+(L)>2$.
 Hence,
${\cal I}(L)=(\frac{2n}{n +2} -\ep,2+\ep')$ and  the isomorphism property follows from Theorem \ref{thinvertibility}.
By Theorem \ref{riesztoreverse}, $r_-(L) \le \sup (1 , (p_-(L))_*)$ and $r_+(L)\ge
p_+(L)=(p_-(L^*))'$. Thus
  $
 \| L^{1/2}f\|_{p}\lesssim
\| \nabla f\|_{p}
$  for $\sup(1,\frac{2n}{n+4}-\ep_1)<p<\frac{2n}{n -2} +\ep_1'$ .
\end{proof}

\begin{remark}  If $n\ge 2$, the bound $p<2+\ep'$ for the Riesz transform $L^p$ boundedness  is sharp  as
$q_+(L)>2$ is optimal (see Section
\ref{sec:sharpness}).  
\end{remark}

A new fact  is the following negative result. 

\begin{corollary} If $n\ge 5$, there exists $p\in (1,2)$ and an operator $L$ for which
$\nabla L^{-1/2} 
$ is not bounded on $L^p$.
\end{corollary}

\begin{proof} If $n\ge 5$, there exists an operator $L$ 
for which $p_-(L)>1$ (see Section
\ref{sec:sharpness}). It remains to invoke   Theorem \ref{thrieszp<2}.
\end{proof}

Let us come to results where further hypotheses may be taken on $L$.

\begin{corollary}\footnote{\, This is first proved in \cite{AT}, Chapter 4, using
\cite{AHLMcT}, Theorem 1.4 and Aronson's estimates in \cite{Ar}.} If
$n\ge 3$ and $L$ has \textbf{real} coefficients, then   $\| \nabla L^{-1/2}
f\|_{p}
\lesssim
 \| f\|_{p}
$ for $1<p<2+\ep$ and $
 \| L^{1/2}f\|_{p}\lesssim
\| \nabla f\|_{p}
$  for $1<p<\infty$. Furthermore,  
 $L^{1/2}$ extends to an isomorphism from $ \dot W^{1,p}$ onto $L^p$ when
$1<p<2+\ep$.
\end{corollary}

\begin{proof}  The semigroup $(e^{-tL})_{t>0}$ is contracting on $L^1$. By Proposition  \ref{propequiv}, part 1, this implies 
$L^1-L^2$ boundedness. By Proposition \ref{propRTallp<2}, we obtain  $\| \nabla L^{-1/2} f\|_{p}
\lesssim
 \| f\|_{p}
$ for $1<p<2$. The rest of the proof is as the one of Proposition \ref{propRTn=2}.
\end{proof}

\begin{remark}  If  $n\ge 2$ and $L$ has continuous and periodic coefficients with common period then a careful spectral analysis
of the semigroup yields that  the Riesz  transform is bounded 
on $L^p$ for all $p\in (1,\infty)$.~\footnote{\, This is in \cite{ERS}, Theorem
1.1.} As a consequence, one obtains uniform gradient bounds on the heat operator from $L^p$ to  
$L^p$, for any $1<p<\infty$.
 
\end{remark}

\begin{remark}  If  $n\ge 2$ and $L$ has real, H\"older continuous and quasiperiodic coefficients
satisfying Koslov condition then by homogeneization techniques,  the Riesz  transform is proved to be 
bounded  on
$L^p$ for all
$p\in (1,\infty)$.~\footnote{\, This is in
\cite{Al}, Theorem 1.2.}  As a consequence, one obtains uniform gradient bounds on the heat operator from
$L^p$ to  
$L^p$, for any $1< p<\infty$.
 
\end{remark}

\begin{remark}  If  $n\ge 2$ and $L$ has almost-periodic coefficients  then 
 the semigroup is $L^1-L^2$ bounded (in fact, much more is true), which is enough to
conclude that  the Riesz  transform is bounded on $L^p$ for all
$p\in (1,2)$. The situation for $p>2$ is unclear. \footnote{\, This is in
\cite{DER}, Theorem 1.1.} 
 
\end{remark}

\subsection{Riesz transforms and Hodge decomposition}\label{sec:hodge}

An $L^p$ Hodge decomposition adapted to the operator $L$ consists in writing a field $f \in L^p$ into the sum
of  fields
$g+\nabla h$ where $g,\nabla h \in L^p$ with $\|g\|_p + \|\nabla h\|_p \le c\|f\|_p$ and $\dv (Ag)=0$.  This amounts to the boundedness of the Hodge
projection
$\nabla L^{-1} \dv $ on $L^p$, that is an inequality of the type 
\begin{equation}\label{eqhodgeprojection}
\| \nabla L^{-1} \dv f\|_{p} \lesssim \| f\|_{p}.
\end{equation}
Indeed, formally, one has
$$
\nabla L^{-1} \dv (A f) = \nabla h.
$$
Alternately, it can be seen as the boundedness for the second order Riesz transform, $\nabla L^{-1/2}
(\nabla (L^*)^{-1/2})^*$.

 For $p=2$, the inequality \eqref{eqhodgeprojection} in automatic by construction of $L$. For $p\ne 2$, we
see here the connection to  Riesz transforms estimates. Before, we put aside the one dimensional case.

\begin{proposition} If $n=1$, $\frac d{dx} L^{-1} \frac d {dx}$ extends to a bounded
operator on $L^p$, $1<p<\infty$, equal to the negative of the operator of
pointwise multiplication with $\frac 1 {a(x)}$.
\end{proposition}

\begin{proof}  Let $D$ be the space of compactly supported and $C^1$ functions $f$ 
with  $\int_{-\infty}^{\infty} \frac
{f(t)}{a(t)} \, dt = 0$. It is easy to see that $D$ is dense in $L^p$ when
$1<p<\infty$. For $f$ in $D$, the unique solution in $\dot W^{1,2}$ of $Lu=\frac
{df}{dx}$ is given by  
$$
u(x)= \int^{\infty}_x \frac
{f(t)}{a(t)} \, dt, \quad x\in \RR.
$$
We have, $u'(x)= -\frac
{f(x)}{a(x)}$ and $u' \in L^p$  with $\|u'\|_p \le \| \frac 1 a \|_\infty
\|f\|_p$ since
$\frac 1 a$ is a bounded function.  Thus the boundedness property holds from $D$
into $L^p$ and the conclusion of the proposition follows readily.
\end{proof}  

 Recall that 
   Riesz transforms inequalities hold in the full
range $1<p<\infty$ when $n=1$.

 We now restrict our attention to dimensions $n\ge 2$. 
In that case, 
\eqref{eqhodgeprojection} is always valid for 
$| \frac 1 2 - \frac 1 p | < \ep$, the value of $\ep$ depending on the ellipticity
constant of $L$ and dimension.

\begin{theorem} 
\begin{enumerate}
\item Let $2<p<\infty$ and $n\ge 2$. Then $\nabla L^{-1} \dv $ is bounded on $L^p$ if and only if  $
 \nabla L^{-1/2}$ is bounded on $L^p$.
 
\item Let $1<p<2$ and $n\ge 3$. If   $\nabla L^{-1} \dv $ is bounded on $L^p$ then  $
 \nabla L^{-1/2}$ is bounded on $L^q$ for $\sup (1,p_*) < q <2$.
\end{enumerate}
\end{theorem}

\begin{proof} Let us begin the argument with the second statement. Let $1<p<2$ and $n\ge 3$. Assume
$
\| \nabla L^{-1} \dv f\|_{p} \lesssim \| f\|_{p}$ for all $f \in L^p$. 
Using interpolation with the corresponding $L^2$-inequality, there is no loss of generality to assume 
$p_*>1$. By Sobolev embeddings we deduce that  $L^{-1}\colon L^{p_*} \to L^{p^*}$ for any such $p$. Take $p_0$ with $p_*<p_0<p$. Define $p_k$ by
$p_k=\big({(p_{k-1})^*}\big)^*$ if $p_0\le 2_*$ and stop when $2_*<p_k\le 2^*$.  

If $k=0$, then we already know that $(e^{-tL})_{t>0}$ is $L^{p_0}-L^2$ bounded,
hence $p_0 \in {\cal J}_-(L)$. By Theorem \ref{thrieszp<2}, we conclude that
$(p_0,2)\subset {\cal I}_-(L)$. Since $p_0$ was arbitrary, we have shown that $(p_*,2) \subset {\cal I}_-(L)$ in this case.

Assume now, that $k\ge 1$.  By construction, $L^{-k}$ is bounded
from
$L^{p_0}$ to
$L^{p_k}$. Since $2_*\le p_k\le 2^*$,
$(e^{-tL})_{t>0}$ is an analytic semigroup on $L^{p_k}$ by Proposition \ref{propequivz}. As in Step 2 of the proof of Theorem \ref{thrieszp<2}, we obtain that
$e^{-tL}$ is bounded from
$L^{p_0}$ to $L^{p_k}$. We also have that $(e^{-tL})_{t>0}$ satisfies $L^{p_k}$ off-diagonal estimates by interpolation between
$L^2$ off-diagonal estimates and $L^r$ boundedness for $r$ chosen so that $p_k$ is between $2$ and $r$. By interpolation again,  
we deduce that if $p_0<q<p_k$, $(e^{-tL})_{t>0}$ satisfies $L^q-L^{p_k}$ off-diagonal estimates. Using Lemma \ref{lemmaDa},
this implies that $(e^{-tL})_{t>0}$ is $L^q$ bounded. Hence, $q  \in {\cal J}_-(L)$. We conclude as above that $(p_*,2) \subset {\cal I}_-(L)$.

Let us now consider the case $p>2$. 
Assume
$
\| \nabla L^{-1} \dv f\|_{p} \lesssim \| f\|_{p}$.  We first claim that $\| {(L^*)}^{1/2} f\|_{p'} \lesssim \| \nabla f\|_{p'}$. Indeed, this always holds if $n=2$
by Proposition \ref{propRTn=2}.  If $n\ge 3$,  by duality we have $\| \nabla (L^*)^{-1} \dv f\|_{p'} \lesssim \| f\|_{p'}$ and the preceeding case tells us that
$\nabla(L^*)^{-1/2}$ is  $L^{p'}$ bounded. Applying Theorem \ref{riesztoreverse} proves the claim.
 
Secondly, let $h \in L^{p'} \cap \dot W^{1,2}$. Since $ \dv h \in L^2$, we have that $${(L^*)}^{-1/2} \dv h= {(L^*)}^{1/2}{(L^*)}^{-1} \dv h,$$ hence
$$
\|{(L^*)}^{-1/2} \dv h \|_{p'} = \|{(L^*)}^{1/2}{L^*}^{-1} \dv h \|_{p'} \lesssim \| \nabla {(L^*)}^{-1} \dv h  \|_{p'} \lesssim \|h\|_{p'}.
$$
Since $ L^{p'} \cap \dot W^{1,2}$ is dense in $L^{p'}$, we obtain the $L^{p'}$-boundedness of ${(L^*)}^{-1/2} \dv$
which, by duality, means the $L^p$-boundedness of the Riesz transform $\nabla L^{-1/2}$. 
 
For the converse, we deduce from the characterization of $\calI(L)$ and the relation between $p_\pm(L)$ and $q_\pm(L)$ that  
 the $L^p$-boundedness of the Riesz transform $\nabla L^{-1/2}$ and $p>2$ imply the $L^{p'}$-boundedness of the Riesz transform $\nabla {(L^*)}^{-1/2}$,
which by duality means that $L^{-1/2} \dv$ is bounded on $L^p$. Hence, $ \nabla L^{-1} \dv=\nabla L^{-1/2} L^{-1/2} \dv  $ is bounded on $L^p$.  
\end{proof}

\begin{remark} By Proposition \ref{propRTn=2}, the second statement is meaningless if $n=2$, hence the assumption $n\ge 3$. 
Its converse  is false for any $p<2$. For example,
for any
$\ep>0$, there exist real symmetric operators $L$  such that $
\| \nabla L^{-1} \dv f\|_{p} \lesssim \| f\|_{p}$ only when $|\frac 1 2 - \frac 1 p | <\ep$, yet its associated 
Riesz transform  is $L^p$ bounded when $1<p<2$.
\end{remark}

\begin{remark} Let us describe a geometric interpretation of the reverse
inequalities for square roots. Let ${\calH}_p(L)=\{ g \in L^p;
\dv(Ag)=0\}$ and ${\cal G}_p =\nabla(\dot W^{1,p})$.  These are closed subspaces of
$L^p$.  The Hodge decomposition in $L^p$ is equivalent to having ${\calH}_p(L) + {\calG}_p = L^p$ as a topological direct sum. 
By duality, it is also equivalent to $\calH_{p'}(L^*)+\calG_{p'}=L^{p'}$ is a topological direct sum.

Assume that $\| {L}^{1/2} f\|_{p} \lesssim \| \nabla f\|_{p}$ holds. Then, by duality $(L^*)^{1/2}$ is bounded from $L^{p'}$ into $\dot  W^{-1,p'}$. If $\nabla h
\in {\cal G}_{p'}$ then 
$(L^*)^{-1/2} \dv (A^*\nabla h) = (L^*)^{1/2} h$ makes  sense. Hence, the restriction of $(L^*)^{-1/2} \dv (A^*\cdot) $ to ${\cal G}_{p'}$ is bounded into
$L^{p'}$.  If $g \in {\calH}_{p'}(L^*)$ then $(L^*)^{-1/2} \dv (A^*g) = 0$, hence the restriction of $(L^*)^{-1/2} \dv (A^*\cdot) $ to ${\calH}_{p'}(L^*)$ is
bounded into
$L^{p'}$ (without any hypothesis). Conversely, these two facts  imply $\| {L}^{1/2} f\|_{p} \lesssim \| \nabla f\|_{p}$.  Thus, this inequality means that
$(L^*)^{-1/2}
\dv (A^*\cdot) $ is bounded on  closed subspaces of
$L^{p'}$ into
$L^{p'}$  even without knowing  whether they are in direct sum in
$L^{p'}$. 

 If the topological direct sum holds then $(L^*)^{-1/2} \dv (A^*\cdot)$
 is bounded on $L^{p'}$, hence $\nabla {L}^{-1/2}$ is bounded on $L^{p}$. This is what we proved above. 

This also illustrates why $\| {L}^{1/2} f\|_{p} \lesssim \| \nabla f\|_{p}$ is possible even when 
$\nabla {L}^{-1/2}$ is not bounded on $L^{p}$.
\end{remark}

Let us finish with another identification of $q_+(L)$.

\begin{corollary} $q_+(L)$ is the supremum of exponents $p$ for which one has the Hodge decomposition in $L^p$, 
or alternatively for which $L$ extends to
an isomorphism from $\dot W^{1,p}$ onto $\dot W^{-1,p}$.~\footnote{\, See
Lemma \ref{lemmaw1p}.} The interval of values of $p$ for which this holds is,
therefore, $((q_+(L^*))',q_+(L))$.
\end{corollary}

\begin{proof} If $n=1$ we have $q_+(L)=\infty$ and the Hodge projections are bounded on all $L^p$ spaces as recalled above. 
If $n \ge 2$, this follows right away from the previous theorem and the fact that $q_+(L)$ is the supremum of $\calI(L)$.
\end{proof}

\section{Riesz transforms and functional calculi}\label{sec:RTfc}

In this section, we present the theorem of Blunck \& Kunstmann concerning
$H^\infty$ functional calculus on $L^p$ spaces.  We also discuss
Hardy-Littlewood-Sobolev inequalities. Combining this with Riesz transform
estimates, we obtain a family of inequalities which we summed up in what we call
the  Hardy-Littlewood-Sobolev-Kato diagram.

\subsection{Blunck \& Kunstmann's theorem}\label{sec:BK}

Let $L$ be as usual and $\omega$ be the type of $L$ defined in  \eqref{type}. 
We know that $L$ admits a bounded holomorphic functional calculus on $L^2$. 

Let $p\in
(1,\infty)$. We say that
$L$ has a bounded holomorphic functional calculus on
$L^p$ if one has the following property: 
for any $\mu\in (\omega,\pi)$ and any $\varphi$ holomorphic and bounded in
$\Sigma_\mu$, and all $f \in L^2\cap L^p$
$$
\|\varphi(L)f\|_p \le c \|\varphi\|_\infty \|f\|_p
$$
the constant $c$ depending only on $p$, $\omega $ and $\mu$. The
operator $\varphi(L)$ is defined on $L^p$ by density.~\footnote{\, Strictly
speaking, one should incorporate a statement about convergence of operators to allow
limiting procedures. In fact, this follows from the  $L^2$-functional calculus  and
density.} 

We define $\calH(L)$ as the sets of those exponents $p$ with the above property.
By interpolation, these sets are intervals (if nonempty).

Recall that
$\calJ(L)$ is the set of exponents $p$ for which $(e^{-tL})_{t>0}$ is $L^p$ {bounded} 
and $\calJ(L)\cap [1,2)=\calJ_-(L)$.
 
\begin{theorem}\label{thBK} The sets ${\calJ}(L)$ and 
${\calH}(L)$
have  same interiors. 
\end{theorem}

By Theorem \ref{thrieszp<2}, this gives us  

\begin{corollary} The sets ${\calI}_-(L)$ and 
${\calH}_-(L)$
have  same interiors. 
\end{corollary}

We turn to the proof of the theorem.

\begin{proof} It is enough to   show that $\calJ_-(L)$ and ${\calH}_-(L)$ have
same interiors as one can use duality in this context.  

Let $p<2$. If $L$ has a bounded holomorphic functional calculus on
$L^p$ then the semigroup $(e^{-tL})_{t>0}$ is $L^p$ bounded. Hence ${\calH}_-(L) \subset \calJ_-(L)$.  

Conversely, it is enough to show that $p_0\in \mathrm{int}{\cal K}_-(L)$, that is $(e^{-zL})_{z\in \Sigma_\beta}$ satisfies
$L^{p_0}-L^2$ off-diagonal estimates for any $\beta \in (0,{\frac \pi 2}-\omega)$ (see Proposition \ref{propequivz}), implies that $\varphi(L)$
is weak-type 
$(p_0,p_0)$ whenever
$\varphi$ is holomorphic and bounded in
$\Sigma_\mu$ for 
$\mu>\omega$ and we may choose $\mu<{\frac \pi 2}$ by convenience. To this end, we apply Theorem \ref{lpp<2} to $T=\varphi(L)$.  It
is enough  to assume further on $\varphi$ the technical condition  
\eqref{eqtechnical}.~\footnote{\, This follows from  the Convergence Lemma in
\cite{Mc}.}

We set $A_r=I-(I-e^{-r^2L})^m$, $r>0$, for
some large integral number  $m$.

Assume   $p_0<2$  is such that $(e^{-zL})_{z\in \Sigma_\beta}$ satisfies
$L^{p_0}-L^2$ off-diagonal estimates for any $\beta \in (0,{\frac \pi 2}-\omega)$. It is enough to check \eqref{pointwisep<2} as \eqref{domap<2} is granted
from the assumption on the semigroup. 
Our goal is to establish the inequality
\begin{equation}\label{check} 
\left(\frac{1}{|2^{j+1}B|}\int_{C_j(B)}|\varphi (L)(I-e^{r^2L})^mf|^2\right)^{1/2}  \le  g(j) \left(\frac{1}{|B|}\int_B |f|^{p_0}\right)^{1/p_0}
\end{equation}
for   all ball $B$ with $r$ the radius of $B$ and all $f$ supported in $B$ and all $j\ge 2$, with $\sum 2^{nj}g(j)<\infty$. 
As before, $C_j(B)$ denotes the ring $2^{j+1}B\setminus
2^j B$.

To do this, let $\psi(z)=\varphi(z)(1-e^{-r^2z})^m$ so that $\psi(L)=\varphi (L)(I-e^{r^2L})^m$.  Represent 
$\psi(L)
$  using the representation formula \eqref{eqrepresentation}.  
Using the exact form of $\psi$ and the definition of $\eta_\pm$, it is easy to obtain
$$
|\eta_\pm(z)| \le C \|\varphi\|_\infty |z|^{-1}\inf \left(1, {r^{2m}} {|z|^{-m}} \right), \quad z \in \Gamma_\pm.
$$
Hence, by Minkowski integral inequality in \eqref{eqrepresentation}, the hypothesis on $p_0$ and on the support of $f$, 
$\|\varphi (L)(I-e^{r^2L})^mf \|_{L^2(C_j(B))}$ is bounded above by
 $$ C \int_{\Gamma_+} e^{-\frac{c4^j r^2}{|z|}}\frac 1 {  |z|^{\gamma_{p_0}/2}}  \inf\left(1,  {r^{2m}} {|z|^{-m}} \right) \frac{d|z|}{|z|} \|f\|_{p_0}
$$
plus the  term corresponding to integration on $\Gamma_-$. A calculation gives us a bound
$$
Cr^{-\gamma_{p_0}} 2^{-j(\gamma_{p_0}+2m)} \|f\|_{p_0}.
$$
Using the value of $\gamma_{p_0}$ we obtain \eqref{check} with $g(j)= C2^{-nj/2} 2^{-j(\gamma_{p_0}+2m)}$. Choosing $m$ with 
$\gamma_{p_0}+2m > n/2$ concludes the argument. 
\end{proof}

\subsection{Hardy-Littlewood-Sobolev estimates}

We prove the following result.

\begin{proposition}\label{propHLS} Let $p_-(L)<p<q<p_+(L)$. Then $L^{-\alpha}$ is bounded from $L^p$ into $L^q$ provided
$$
\alpha= \frac 1 2 \left(\frac n p - \frac n q\right).
$$
\end{proposition} 

\begin{proof} We first observe that with the following  choice of $p$ and $q$ then $(e^{-tL})_{t>0}$ is $L^p-L^q$ bounded.
Indeed, we know it if $p=2$ or if $q=2$ by Proposition \ref{propequiv}. If $p<2<q$, it suffices to use composition. If $p<q<2$,  interpolate between $L^p$
boundedness and
$L^p-L^2$ boundedness. And if $2<p<q$,  interpolate between $L^q$ boundedness and $L^2-L^q$ boundedness.

Next, by  $L^2$ functional calculus, we have
$$
 f = \frac{1}{\Gamma(\alpha)} \int_0^\infty t^{\alpha -1}L^{\alpha} e^{-tL}f\, dt
$$
for all $f\in  L^2$, where the integrals $\int_\ep^R \ldots$ converge in $L^2$ as $\ep \downarrow 0$ and $R\uparrow \infty$.  
Set $T_{\ep,R}= \Gamma(\alpha)^{-1}\int_\ep^R t^{\alpha -1} e^{-tL}\, dt.$
Let $f\in L^2\cap L^p$ with 
$\|f\|_p=1$. Fix $\alpha=
\frac 1 2 (\frac n p - \frac n q).
$
For $a>0$ and $p<q_0<q<q_1 <q_+(L)$, we easily obtain with uniform constant $C$, 
$$
\left\| 
\int_\ep^a t^{\alpha -1} e^{-tL}f \, dt\right \|_{q_1} \le Ca^{\frac 1 2 (\frac n q - \frac n {q_1})}
$$
and
$$
\left\| 
\int_a^R t^{\alpha -1} e^{-tL}f \, dt\right \|_{q_0} \le Ca^{\frac 1 2 (\frac n {q_0} - \frac n {q})}.
$$
Hence, by the argument of Marcinkiewic interpolation theorem, we have if $\lambda>0$ 
$$
|\{ x \in \RR^n; \ |T_{\ep,R}f(x)| > \lambda\}| \le   C\lambda^{-q_1}a^{\frac {q_1} 2 (\frac n q - \frac n {q_1})}  + C\lambda^{-q_0} a^{\frac {q_0} 2 (\frac n
{q_0} -
\frac n {q})}.
$$
Choosing $a^{  \frac n {2q}}= \lambda$ yields,
$$
|\{ x \in \RR^n; \ |T_{\ep,R}f(x)| > \lambda\}| \le 2 C\lambda^{-q}.
$$
Hence, $T_{\ep,R}f$ belongs to the Lorentz space $L^{q,\infty}$. Since this holds for all $q$ as above, by interpolation again, we conclude that
$T_{\ep,R}f\in L^q$ and
$$
\|T_{\ep,R}f\|_q \le C \|f\|_p
$$ 
whenever $f\in  L^2\cap L^p$ with uniform constant with respect to $\ep, R$.

 It remains to pass to the limit. 
If $0\le \alpha\le 1/2$, we claim that $L^{\alpha}$ defines an isomorphism between the fractional Sobolev space $\dot H^{2\alpha}$ (defined as the closure of
$C_0^\infty$ for the seminorm $\|f\|_{\dot H^{2\alpha}}\equiv \|(-\Delta)^\alpha f\|_2$)  and
$L^2$ with 
$$\|L^\alpha f\|_2 \sim  \|f\|_{\dot H^{2\alpha}}.
$$
Indeed, if $\Re \alpha=0$, $z\mapsto z^\alpha$ is bounded holomorphic in $\Sigma_\pi$, hence
$$\|L^\alpha f\|_2 \sim  \|f\|_2
$$
with  implicit constants  independent of $\Im\alpha$. For $\Re\alpha=1/2$,
combining the square root problem \eqref{eq5} with the $L^2$ functional calculus,
we have 
$$\|L^\alpha f\|_2 \sim  \|\nabla f\|_2
$$
 with  implicit constants  independent of $\Im\alpha$. The claim follows by complex interpolation.  
By construction, if $f\in L^2$ then
$$\lim_{\ep \downarrow 0, R\uparrow \infty} \| T_{\ep,R}f - L^{-\alpha}f\|_{\dot H^{2\alpha}},
$$
hence, for $\alpha$ sufficiently small so that Sobolev embeddings $\dot H^{2\alpha} \subset L^r$ applies with  $r<\infty$, we also  have the convergence in $L^r$.  

Now, if $f\in L^2\cap L^p$, we combine the uniform bound of $T_{\ep,R}f$ in $L^q$ and its  convergence in $L^r$ to conclude for the
$L^p-L^q$ bounded extension of $L^{-\alpha}$ for all small positive $\alpha$. 

We obtain all possible values of $\alpha$ by writing $L^{-\alpha}= (L^{-\alpha/k})^k$ for $k$ large enough.
\end{proof} 

\begin{remark}\label{remarkKL} The isomorphism property used in the proof holds $0\le \alpha<1/2$ without
knowing the solution of the Kato square root problem using that  the domain of $L^\alpha$ 
is an interpolation space between the form domain (\textit{i.e.} $ W^{1,2}$)
and $L^2$.~\footnote{\, This  complex interpolation result  follows  from Kato
\cite{K1} and   Lions
\cite{Li}.}  
\end{remark}

\subsection{The Hardy-Littlewood-Sobolev-Kato diagram}\label{sec:HLSK}

In the plane $\{(\frac 1 p, s)\}$, where $p$ is a Lebesgue exponent and $s$  a regularity
index, we introduce a convex set on which we have a rule for  boundedness  from $\dot W^{s,p}$ to $\dot W^{\sigma,q}$
of functions of $L$. We call it the Hardy-Littlewood-Sobolev-Kato diagram of $L$ because it includes all previously
seen estimates from functional calculus and Kato type inequalities (see below). 

    Pick exponents $p_0,q_0,p_1,q_1$ as follows:
 $$q_-(L)=p_-(L) < p_0 < 2 < q_0 < q_+(L)$$ and $$q_-(L^*)=p_-(L^*) < p_1 < 2 < q_1 <
q_+(L^*).$$ Hence, $p_0, q_0$ (resp. $p_1,q_1$) are in the range of
$L^p$-boundedness of the Riesz transform $\nabla L^{-1/2}$  (resp.  $\nabla
(L^*)^{-1/2}$). Consider 
the closed convex polygon ${\calP}=ABDFEC$ in the  plane $\{(\frac 1 p, s)\}$ with $A=(\frac 1 {q_0}, 1)$,
$B= (\frac 1 {p_0}, 1) $,  $D=(\frac 1 {p_0}, 0) $
$F=(\frac 1 {(q_1)^\prime}, -1)$, $E=(\frac 1 {(p_1)'}, -1)$ and $C= (\frac 1 {(p_1)'}, 0)$.  See Figure 1. 
Let  $M=(\frac 1 p, s), N=(\frac 1 q, \sigma) $ in $ {\calP}$. We call 
$\overrightarrow{MN}$  an authorized arrow if $p\le q$.  
Set in this case,
$$\alpha (M,N) =  \frac {\sigma -s} 2 + \frac 1 2 \left(\frac n p - \frac n q\right).$$
In other words, we are allowed to move in $\calP$ horizontally from the right to the
left, vertically up and down and all possible combinations. In fact, an accurate
correspondance would be  the convex three dimensional set of authorized
arrows $\overrightarrow{MN}$ between two parrallel copies of $\calP$ in $\RR^3$  whenever 
$M$ in the first copy and $N$ in the second copy.

\begin{figure}[h]
\begin{center}
\resizebox*{0.3\textwidth}{!}{\includegraphics{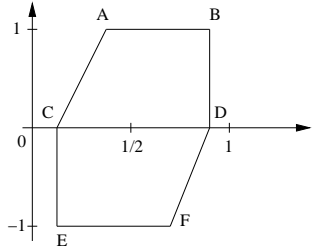}}
\caption{{\small \bf The Hardy-Littlewood-Sobolev-Kato diagram $\calP$.}
\label{Figure1}}
\end{center}
\end{figure}

\begin{theorem}\label{thHLSK} Let $\mu \in (\omega, \pi)$. For any 
$M=(\frac 1 p, s), N=(\frac 1 q, \sigma) \in  {\calP}$ with $p\le q$ and any $\varphi
\in \calF(\Sigma_\mu)$,
$$
\|\varphi(L)f \|_{\dot W^{\sigma, q}} \lesssim \|z^{\alpha(M,N)} \varphi\|_{L^\infty(\Sigma_\mu)} \|f\|_{\dot W^{s, p}},
$$
provided the quantity $\|z^{\alpha(M,N)} \varphi\|_{L^\infty(\Sigma_\mu)}$ is finite. The implicit constant depends on $\mu,\omega, p,q,s,\sigma$.
\end{theorem}

A few explanations are necessary. First, for $0<s$, $\dot W^{s,p}$ is defined as the closure of
$C_0^\infty$ for $\|(-\Delta)^{s/2}f\|_p$ and $\dot W^{-s,p'}$ is its dual space.
 If
$0<
\alpha< 1$,
$z\mapsto z^\alpha$  is the analytic continuation in $\Sigma_\pi$ of $t\mapsto
t^\alpha$ defined on
$(0,\infty)$. If $k<\alpha<k+1$ with $k$ integral number, then $z \mapsto
z^\alpha=z^kz^{\alpha -k}$ is analytic in $\Sigma_\pi$. Hence, if $\varphi\in
\calF(\Sigma_\mu)$ then $z^\alpha \varphi\in \calF(\Sigma_\mu)$ and $L^\alpha
\varphi(L)$ is well-defined (See Section \ref{sec:fcL2}). Also, we shall write
\textit{a priori} inequalities for suitable functions $f$ and we shall  leave to the
reader the care of providing the density arguments. 

 When $s=\sigma=0$, this includes $L^p-L^q$ estimates for negative fractional powers of $L$ (Hardy-Littlewood-Sobolev estimates).
For $p=q$, we have a comparison between fractional powers
 of $L$ with the ones of the Laplacian $-\Delta$ (Kato type estimates): the operators
$(-\Delta)^{\sigma/2} L^{(s-\sigma )/2} (-\Delta)^{-s/2}$ are bounded on $L^p$
provided
$(\frac 1 p, s)$ and $(\frac 1 p, \sigma)$ belong to
 ${\calP}$. More precisely

\begin{corollary} If  $(\frac 1 p, s)$ and $(\frac 1 p, \sigma)$ belong to
 ${\calP}$, then $L^{(s-\sigma )/2}$ extends to an isomorphism from $\dot W^{s,p}$
onto
$\dot W^{\sigma,p}$
\end{corollary}
 
\paragraph{Proof:} The previous remark means that 
$$
\|L^{(s-\sigma )/2}f \|_{\dot W^{\sigma, p}} \sim  \|f\|_{\dot W^{s, p}}
$$
which implies that the extension from $C_0^\infty$ to $\dot W^{s, p}$ is bounded and
one-to-one. Since the same thing is true for $L^*$ in the dual range, we conclude
that this extension is onto by classical arguments.
\qed
 
In particular we recover  the simultaneous $L^p$ boundedness of the  Riesz
transform  and the Hodge projection when $p>2$. We pursue the discussion on Kato
type estimates in the next section.

\paragraph{Proof of Theorem \ref{thHLSK}:}
By convexity and complex interpolation, its suffices to prove the result for the following extremal authorized arrows:
\begin{enumerate}
\item the minimal ones: $\overrightarrow{MM}$ where $M$ is one of the six vertices of $\calP$. 
\item  the  maximal horizontal  ones:  $\overrightarrow{BA}$, $\overrightarrow{DC}$, $\overrightarrow{FE}$.
\item  the maximal vertical ones: 
${\overrightarrow{DB}}$, ${\overrightarrow{FH}}$, ${\overrightarrow{GA}}$,
${\overrightarrow{EC}}$   and their respective opposite where $G=(\frac {1}{p_0},
-1)$ and $H=(\frac {1}{(q_1)'}, 1)$.
\end{enumerate}
 
\paragraph{Step 1. Minimal arrows}\hfill\break

$\bullet$ {$\overrightarrow{AA}$:} Write
$$
\|\nabla \varphi(L) f \|_{q_0} \lesssim \|L^{1/2} \varphi(L) f\|_{q_0} \sim \| \varphi(L) L^{1/2} f\|_{q_0} \lesssim \|\varphi\|_\infty \|L^{1/2}f\|_{q_0} \lesssim
\|\nabla f\|_{q_0}.
$$
The first inequality holds because the Riesz transform is $L^{q_0}$ bounded, the second by the commutative property of the functional calculus, the third by Theorem
\ref{thBK} and the last by the reverse inequalities at $q_0$.

$\bullet${$\overrightarrow{BB}$:} Same as $\overrightarrow{AA}$ by changing $q_0$ to $p_0$.

$\bullet${$\overrightarrow{CC}$:} We have $(p_1)' \in \mathrm{int}\calJ(L)$,
 hence Theorem \ref{thBK} applies and yields
$\| \varphi(L) f \|_{(p_1)'} \lesssim
\|\nabla f\|_{(p_1)'}.
$

$\bullet${$\overrightarrow{DD}$:} Same as $\overrightarrow{CC}$ by changing $(p_1)'$ to $p_0$.

$\bullet${$\overrightarrow{EE}$:} By duality, we have to show 
$
\|\nabla \varphi(L^*) f \|_{p_1} 
\lesssim
\|\nabla f\|_{p_1}.
$
This is the same as $\overrightarrow{BB}$ by changing $L$ to $L^*$ and $p_0$ to $p_1$.

$\bullet${$\overrightarrow{FF}$:} By duality, this is the same as $\overrightarrow{AA}$ by changing $L$ to $L^*$ and $q_0$ to $q_1$. 

\paragraph{Step 2. Maximal horizontal arrows}\hfill\break

$\bullet${$\overrightarrow{BA}$:}
As  for  $\overrightarrow{AA}$, we begin with 
$$
\|\nabla \varphi(L) f \|_{q_0} \lesssim \|L^{1/2} \varphi(L) f\|_{q_0} \sim \| \varphi(L) L^{1/2} f\|_{q_0}.
$$ 
Next, we continue with 
$$
\| \varphi(L) L^{1/2} f\|_{q_0} \lesssim \|z^{\alpha(B,A)}\varphi\|_\infty
\|L^{-\alpha(B,A) }L^{1/2}f\|_{q_0}
$$
by Theorem \ref{thBK}, then 
$$\|L^{-\alpha(B,A) }L^{1/2}f\|_{q_0} \lesssim
\|L^{1/2}f\|_{p_0} \lesssim
 \|\nabla f\|_{p_0},
$$
by the Hardy-Littlewood-Sobolev inequality for $L$ and the definition of $\alpha (B,A)$, and the reverse inequality at $p_0$.

$\bullet${$\overrightarrow{DC}$:} Since $p_-(L)<p_0<(p_1)'<p_+(L)$, by Proposition \ref{propHLS} and Theorem \ref{thBK}, 
$$\| \varphi(L) f \|_{(p_1)'} \lesssim  \|z^{\alpha(D,C)}\varphi\|_\infty
\|L^{-\alpha(D,C) }f\|_{(p_1)'} \lesssim \|f\|_{p_0}.
$$

$\bullet${$\overrightarrow{FE}$:} By duality, this is the same as $\overrightarrow{BA}$ by changing 
$L$ to $L^*$, $p_0$ to $p_1$ and $q_0$ to $q_1$.

\paragraph{Step 3. Maximal vertical arrows}\hfill\break

Recall that $G$ and $H$ are points    in $\calP$ given by $G=(\frac {1}{p_0}, -1)$ and $H=(\frac {1}{(q_1)'}, 1)$.

$\bullet${$\overrightarrow{DB}$:} Using that the Riesz transform is $L^{p_0}$ bounded and Theorem \ref{thBK}, one has
$$
\|\nabla \varphi(L) f \|_{p_0} \lesssim \|L^{1/2} \varphi(L) f\|_{p_0}  \lesssim \|z^{1/2}\varphi\|_\infty \|f\|_{p_0}.
$$

$\bullet${$\overrightarrow{BD}$:} We have
$$
\| \varphi(L) f \|_{p_0} \lesssim \|L^{-1/2} \varphi(L) L^{1/2} f\|_{p_0}  \lesssim
\|z^{-1/2}\varphi\|_\infty \| L^{1/2}f\|_{p_0}\lesssim \|\nabla f\|_{p_0}.
$$

$\bullet${$\overrightarrow{FH}$:} 
Since $(q_1)'>p_-(L)$, one has the $L^{(q_1)'}$ boundedness of the Riesz transform
associated to $L$:
$$
\|\nabla \varphi(L) \dv f \|_{(q_1)'} \lesssim \|L^{1/2} \varphi(L) \dv f\|_{(q_1)'}. 
$$
Then writing $L^{1/2}\varphi(L)= L\varphi(L) L^{-1/2}$ and using Theorem \ref{thBK} yield
$$
\|L^{1/2} \varphi(L) \dv f\|_{(q_1)'}\lesssim
\|z\varphi\|_\infty
\|L^{-1/2}\dv f\|_{(q_1)'}.
$$
Next, using that the Riesz transform associated to $L^*$ is bounded on $L^{q_1}$, one concludes by
$$
\|L^{-1/2}\dv f\|_{(q_1)'}\lesssim
\|f\|_{(q_1)'}.
$$

$\bullet${$\overrightarrow{HF}$:} Since $(q_1)'>p_-(L)$, $L^{1/2}$ is bounded from 
$L^{(q_1)'}$  to $\dot W^{-1,(q_1)'}$, hence
$$
\| \varphi(L)  f \|_{\dot W^{-1,(q_1)'}}= \| L^{1/2} \varphi(L) L^{-1/2}  f
\|_{\dot W^{-1,(q_1)'}}
\lesssim \| \varphi(L) L^{-1/2} f\|_{(q_1)'}.
$$
Next, by Theorem \ref{thBK},
$$
\| \varphi(L) L^{-1/2} f\|_{(q_1)'} \lesssim \| z^{-1} \varphi\|_\infty 
\|L^{1/2}f\|_{(q_1)'} 
$$
We finish with $L^{1/2}$  bounded from 
 $\dot W^{(q_1)'}$ to $L^{(q_1)'}$ since $(q_1)'>p_-(L)$.

$\bullet${$\overrightarrow{GA}$:} Same as $\overrightarrow{FH}$ by changing $(q_1)'$ to $q_0$.

$\bullet${$\overrightarrow{AG}$:} Same as $\overrightarrow{HF}$ by changing $(q_1)'$
to $q_0$.

$\bullet${$\overrightarrow{EC}$:} By duality, this is the same as $\overrightarrow{DB}$ by changing $L$ to $L^*$ and $p_0$ to $p_1$.

$\bullet${$\overrightarrow{CE}$:} By duality, this is the same as
$\overrightarrow{BD}$ by changing $L$ to $L^*$ and $p_0$ to $p_1$.

\bigskip 

We conclude this section with the following corollary.

\begin{corollary} If $q_-(L)<p<q_+(L)$, then $L$ has a bounded holomorphic
functional calculus on $\dot W^{1,p}$.~\footnote{\, We do not know if the converse
holds.}
\end{corollary}

\begin{proof} $(q_-(L),q_+(L)) =\calI(L)$, so this functional calculus corresponds
to minimal arrows $\overrightarrow{MM}$ for $M=(\frac 1 p, 1)$.
\end{proof}

\subsection{More on the Kato diagram}\label{sec:more}

As we have seen, we can move vertically up and down in $\calP$. However, one can
authorize more downward arrows. Doing this and modifying slightly the definition of
the numbers  $r_\pm(L)$ in Theorem \ref{riesztoreverse}, we shall show that
$r_+(L)=p_+(L)$   and obtain a lower bound on $r_-(L)$.

First we claim that the proof given for the reverse inequalities
$\|L^{1/2}f\|_p \lesssim \|\nabla f\|_p$ extends as follows.

\begin{proposition} For all $\mu \in (\omega, \pi)$ and all
$\varphi \in H^\infty(\Sigma_\mu)$,  
$$\|L^{1/2}\varphi(L)f\|_p \lesssim \|\varphi\|_\infty\|\nabla f\|_p$$
whenever 
$\tilde p_-(L)=r_-(L) <p< p_+(L)$.
\end{proposition}

\begin{proof} The case $p>2$ is a simple consequence of this inequality for $L^{1/2}$
together with the bounded holomorphic functional calculus on $L^p$ for $p<p_+(L)$.

We now turn to $p<2$. We begin with the representation formula
\eqref{eqrepresentation} for $\psi(z)=z^{1/2} \varphi(z)$. One has  
$$
\psi(L)= \int_{\Gamma_+} e^{-zL} \eta_+(z)\, dz + \int_{\Gamma_-} e^{-zL}
\eta_-(z)\, dz 
$$
and
$$
\eta_\pm(z)= 
\frac 1 {2\pi i} \int_{\gamma_\pm} e^{\zeta z} \psi(\zeta)\, d\zeta,
\quad z
\in \Gamma_\pm,
$$
valid provided the technical assumption $|\varphi(\zeta)| \le C (1+|\zeta|)^{-1/2-s}$
for some
$s>0$. If one defines  the primitive $N_\pm(z)$ of respectively $\eta_\pm(z)$ which
vanishes at infinity then, under the technical assumption, one may integrate by parts
and, since the terms at 0 cancel each other, one finds
$$
\psi(L)= \int_{\Gamma_+} Le^{-zL} N_+(z)\, dz + \int_{\Gamma_-} Le^{-zL}
N_-(z)\, dz.
$$ 
Furthermore, one has
$$
|N_\pm(z)| \le C |z|^{-1/2} \|\varphi\|_\infty, \quad z\in \Gamma_\pm.
$$
Hence, by repeating the argument for $L^{1/2}$ where integration on the positive
axis is replaced by integration on half-rays $\Gamma_\pm$, we obtain 
$$\|\psi(L)f\|_p \lesssim \|\varphi\|_\infty\|\nabla f\|_p$$
whenever 
$\tilde p_-(L) <p< 2$ and $\varphi$ sastisfies the technical condition, which is
removed by a limiting argument. 
\end{proof}

 Let us discuss the sharpness of the bounds obtained above.

\begin{proposition} Assume that for all $\mu \in (\omega, \pi)$ and all
$\varphi \in H^\infty(\Sigma_\mu)$,  
$$\|L^{1/2}\varphi(L)f\|_p \lesssim \|\varphi\|_\infty\|\nabla f\|_p.$$
Then  
$\sup(\ol{ p_-}(L),1) \le p \le p_+(L)$ where $\ol{ p_-}(L)$ is defined as follows:
the point $\ol B=(\frac 1 {\ol{ p_-}(L)},1)$ is symmetric to $F=(\frac 1 
{(q_+(L^*))'},-1)$ with respect to $D=(\frac 1 
{(p_+(L^*))'},0)=(\frac 1 
{p_-(L)},0)$.
\end{proposition}

 One can see that  $\ol{p_-}(L)\le \tilde p_-(L)$.   This lower bound ${\ol{ p_-}(L)}$ is the
best one can obtain by a convexity method. If $n\le 4$, $\tilde p_-(L)=1$ so there is nothing more to
say. In dimensions
$n\ge 5$, the lower bound is optimal for any operator
$L$ for which  this inequality cannot be improved, which is the same as saying that
$p_+(L^*) \ge (q_+(L^*))^{*}$ is optimal (See Section \ref{sec:sharpness} for this).

\begin{proof}  Fix $\varphi \in H^\infty(\Sigma_\mu)$ with 
$\|\varphi\|_\infty=1$. 
We have for all $t\in
\RR$ 
$$\|L^{1/2+it}\varphi(L)f\|_p \lesssim \|\nabla f\|_p=\|f\|_{\dot W^{1,p}}.$$
We also have from the HLSK diagram that for all $(q_+(L^*))'<q<p_+(L)$ and $t\in
\RR$,
$$\|L^{-1/2+it}\varphi(L)f\|_q \lesssim \| f\|_{\dot W^{-1,q}}.$$
By complex interpolation, we obtain 
$$
\|\varphi(L)f\|_r \lesssim \|f\|_r
$$
for $1/r$ the middle of $1/p$ and $1/q$. By Theorem \ref{thBK}, it is necessary that 
$p_-(L)\le r \le p_+(L)$. 

Choosing $q$ arbitrarily close to $p_+(L)$ 
forces $p\le p_+(L)$. Choosing $q$ arbitrarily close to $(q_+(L^*))'$ forces
$p\ge {\ol{ p_-}(L)}$ given the definition of this number.
\end{proof}

Geometrically, this provides us with a family of authorized downarrows that may not
be contained in 
$\cal P$.  We  assume that $\tilde p_0 =(p_0)_*>1$ otherwise we exclude from this discussion points
$\{(\frac 1 p, s)\}$ with $\frac 1 p \ge 1$.   Consider  the closed convex  polygon
${\calP}_{in}=A\tilde BDFEC$ in the  plane
$\{(\frac 1 p, s)\}$ with  
$A,D,F,E,C$ as before and 
$\tilde B= (\frac 1 {\tilde p_0}, 1) $ and also the closed convex polygon
${\calP}_{ext}=\tilde A\tilde B\tilde DFE$, where $
\tilde A=(\frac 1 {(p_1)'}, 1)$ and $\tilde D=(\frac 1 {\tilde p_0}, 0) $.  See Figure 2.

\begin{figure}[h]
\begin{center}
\resizebox*{0.3\textwidth}{!}{\includegraphics{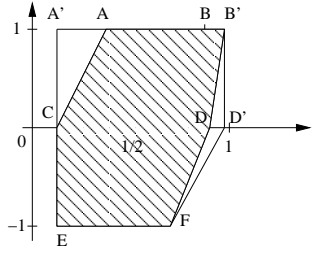}}
\caption{{\small  \bf The convex polygons  ${\calP}_{in}$ (shaded) and ${\calP}_{ext}$.   The points $A',B',D'$ correspond to $\tilde A,\tilde B,\tilde D$ in the text. Downarrows ${MN}$
with   $M\in {\calP}_{in}$ and   $N\in {\calP}_{ext}$ or  $M\in {\calP}_{ext}$ and   $N\in {\calP}_{in}$ are authorized.}
\label{Figure 2}}
\end{center}
\end{figure}

\begin{proposition} For any downarrow {$\overrightarrow{MN}$}, $M=(\frac 1 { p}, s)$ 
$N=(\frac 1 { p}, \sigma)$ with one extremity ${\calP}_{in}$ and the other in
${\calP}_{ext}$ we have for all $\mu \in (\omega, \pi)$ and all
$\varphi \in H^\infty(\Sigma_\mu)$,
$$
\|\varphi(L)f \|_{\dot W^{\sigma, p}} \lesssim \|z^{\alpha(M,N)} \varphi\|_{L^\infty(\Sigma_\mu)}
\|f\|_{\dot W^{s, p}},
$$
provided the quantity $\|z^{\alpha(M,N)} \varphi\|_{L^\infty(\Sigma_\mu)}$ is finite. The implicit constant depends on $\mu,\omega,
p,q,s,\sigma$.
\end{proposition} 

\begin{proof} The previous proposition shows that  for the maximal downarrows   {$\overrightarrow{\tilde A
C}$} and  {$\overrightarrow{\tilde B \tilde D}$} the corresponding inequality of the statement are
valid. Combining this with all possible vertical authorized arrows in $\calP$, the convexity property of
complex interpolation gives us the  desired result. 
\end{proof}

Of course, the downarrows can be reversed exactly when $M$ and $N$ belong to the HLSK diagram $\calP$.

The expected maximal convex set $\calP_{in}$ would be with $\tilde B$ replaced by
$\ol B$ where $\ol B$ is the symmetric point of $F$ with respect to $D$. But our arguments do not
suffice  to prove the inequality corresponding to the vertical downarrow {$\overrightarrow{\ol {B} \ol
{D}}$} with $\ol D$ the vertical projection of $\ol B$ on the $\frac 1 p $ axis. 

\section{Square function estimates}\label{sec:RTsf}

In this section, we study vertical square functions of two different types which are representative of larger classes of square functions. Then, we prove some weak type
and strong type inequalities for non-tangential square functions.

\subsection{Necessary and sufficient conditions for boundedness of vertical square functions}

Define the quadratic functionals for $f \in L^2$
$$g_{L}(f)(x)= \left( \int_0^\infty |(L^{1/2}e^{-tL}f)(x)|^2 \, {dt}\right)^{1/2}
$$
and 
$$
G_{L}(f)(x)= \left( \int_0^\infty |(\nabla e^{-tL}f)(x)|^2 \, {dt}\right)^{1/2}.
$$
The $L^2$ theory of  quadratic estimates for  operators  having a  bounded holomorphic functional
calculus on
$L^2$ implies 
\begin{equation}\label{gLL2bdd}
\|g_{L}(f)\|_2 \sim \|f\|_2.
\end{equation}
In fact, such an inequality is equivalent to the bounded holomorphic functional
calculus on
$L^2$. Moreover, the family  $((tL)^{1/2}e^{-tL})_{t>0}$ in $g_L$ can be replaced by more general functions
of $L$. \footnote{\, All this is due to McIntosh and Yagi, \cite{McY} and
\cite{Y}.} 

As for $G_L$,  
\begin{equation}\label{GLL2bdd}
\|G_{L}(f)\|_2 \sim \|f\|_2
\end{equation}
is a consequence of ellipticity and 
$$
\int_{\RR^n}|f(x)|^2\, dx = 2 \Re \iint_{\RR^n \times (0,\infty)} A(x) (\nabla e^{-tL} f) (x)\cdot \ol{(\nabla e^{-tL} f)(x)}\, dxdt.
$$
This equality is obtained as follows. We have
\begin{align*}
\|f\|_2^2 &= -\int_0^\infty \frac d{dt} \|e^{-tL}f\|^2_2\, dt
\\
& = 2\Re \iint_{\RR^n \times (0,\infty)}  (Le^{-tL}f)(x) \ol{( e^{-tL} f)(x)}\, dxdt
\end{align*}
and it remains to integrate by parts in the $x$ variable using the definition of $L$.

We are interested in the $L^p$ counterparts of these results. Define
$$\calS(L)=  \{1<p< \infty;\ \forall f \in L^2\cap 
L^p \quad \| g_L(f)\|_p \sim \| f\|_p.\}
$$
and $\calS_-(L)=\calS(L)\cap (1,2)$ and
$$\calG(L)=  \{1<p <\infty; \  \forall f \in L^2\cap 
L^p\quad \| G_L(f)\|_p\sim \| f\|_p.\} , 
$$
 $\calG_-(L)= \calG(L) \cap (1,2)$ and $\calG_+(L)= \calG(L) \cap (2,\infty)$.  Recall 
$\mathrm{int}\calJ(L)=(p_-(L),p_+(L))$ and $\mathrm{int}\calN(L)=(q_-(L),q_+(L))$.

\begin{theorem} \begin{enumerate}
\item  $\mathrm{int}\calS(L)=(p_-(L),p_+(L))$.~\footnote{\, We develop an argument
based on the extension of Calder\'on-Zygmund theory in Section \ref{sec:CZ}.
However, there is an other argument as follows: we have shown that the
interior of $\calJ(L)$ is the maximal open interval of exponents $p$ for which
$L$ has a bounded holomorphic functional calculus of $L^p$. Take such a $p$.
Then, apply Le Merdy's Theorem 3 in \cite{LeM} to conclude that $g_L$ defines a
new norm on $L^p$. More is proved in \cite{LeM}: the operators $(tL)^{1/2}e^{-tL}$
in the definition of $g_L$ can be replaced by more general functions of $L$. We
choose  not to develop this here.}  
\item $\mathrm{int}\calG(L)=(q_-(L),q_+(L))$.~\footnote{\, The statement is valid
for square functions where the semigroup in the definition of $G_L$ is replaced
by more general functions of $L$. Again, we do not go into such considerations.}
\end{enumerate}
\end{theorem}

Roughly, this theorem says that, up to endpoints, the range of $p \in (1,\infty)$ for which 
 $g_L$ defines a new norm on $L^p$ is the same as the one of boundedness of the semigroup
and that the range of $p\in (1,\infty)$ for which $G_L$ defines a new norm on $L^p$ is the same as the
one of boundedness of $(\sqrt t\, \nabla e^{-tL})_{t>0}$.

Due to Theorems \ref{thrieszp<2} and \ref{thrieszp>2}, the connections to  intervals of exponents $p$ for
which one has  $L^p$ boundedness for the Riesz  transform $\nabla L^{-1/2}$ is as follows.

\begin{corollary} \begin{enumerate}
\item The intervals $\calI_-(L)$,  $\calS_-(L)$ and $\calG_-(L)$ have same interiors. 
\item The intervals $\calI_+(L)$ and $\calG_+(L)  $  have same interiors.
\end{enumerate}
\end{corollary}

Hence, there is again a dichotomy $p<2$ vs $p>2$ in the description. We turn to the proof of the theorem.

\begin{proof} The argument has several steps. By Proposition \ref{propequiv}, we
may freely replace
$\calJ_\pm(L)$ by one of the intervals  $\calK_\pm(L)$, $\calM_\pm(L)$.

\paragraph{Step 1. $p_0>2$ and $(\sqrt t\, \nabla e^{-tL})_{t>0}$ satisfies $L^2-L^{p_0}$ off-diagonal
estimates imply $\|G_L(f)\|_p \lesssim \|f\|_p$ for
$2<p<p_0$.}

\paragraph{Step 2. $p_0>2$ and $(e^{-tL})_{t>0}$ satisfies $L^2-L^{p_0}$ off-diagonal estimates imply
$\|g_L(f)\|_p \lesssim \|f\|_p$ for $2<p<p_0$.}

\paragraph{Step 3. $p_0<2$ and $(e^{-tL})_{t>0}$ satisfies $L^{p_0}-L^2$ off-diagonal estimates imply
$\|g_L(f)\|_p \lesssim \|f\|_p$ for $p_0<p<2$.  }

\paragraph{Step 4. $p_0<2$ and $(\sqrt t\, \nabla e^{-tL})_{t>0}$ satisfies $L^{p_0}-L^2$ off-diagonal
estimates imply
$\|G_L(f)\|_p \lesssim \|f\|_p$ for $p_0<p<2$.  }

\paragraph{Step 5. Reverse  $L^p$ inequality for $g_L$ when $p_-(L)<p<p_+(L)$. }

\paragraph{Step 6. 
$\|g_L(f)\|_p \sim \|f\|_p$ implies $(e^{-tL})_{t>0}$  $L^p$ bounded.
}

\paragraph{Step 7. $\|G_L(f)\|_p \lesssim \|f\|_p$   implies $(\sqrt t\, \nabla e^{-tL})_{t>0}$
$L^q$-bounded for $q$ in the interval between 2 and $p$. } 

\paragraph{Step 8. Reverse $L^p$ inequality for $G_L$ when $1<p<\infty$.  } \hfill\par\bigskip

The combination of steps 2, 3 and 5 shows that $ \mathrm{int}\calJ(L)$ is contained
in 
 $ \calS(L)$. Step 6 implies
that $\calS(L)$ is contained in $\calJ(L)$. 

The combination of steps 1 and 4 shows that $\mathrm{int}\calN(L) $ is contained
in   $\mathrm{int}\calG(L)$. Step 7 and Step 8 show the converse.

\paragraph{Step 1. $p_0>2$ and $(\sqrt t\, \nabla e^{-tL})_{t>0}$ satisfies $L^2-L^{p_0}$ off-diagonal
estimates implies $\|G_L(f)\|_p \lesssim \|f\|_p$ for
$2<p<p_0$.}

We apply Theorem \ref{lpp>2} to the sublinear operator $T=G_L$. We choose as usual $A_r=I-(I-e^{-r^2L})^m
$ for $m$ a large enough integral number. 

We have to check \eqref{doma}. Let $B$ be a ball and $r=r(B)$ be its radius and $k=1, \ldots, m$. 
By Minkowski integral inequality and $p_0>2$, we have
$$
\left(\frac 1{|B|}\int_B |G_L(e^{-kr^2L}f)|^{p_0}\right)^{2/p_0} \le \int_0^\infty \left(\frac
1{|B|}\int_B | \nabla e^{-tL}(e^{-kr^2L}f)|^{p_0}\right)^{2/p_0}\,  {dt}$$
Using the hypothesis of $p_0$,  the commutativity property of the semigroup and applying the 
scalar inequality \eqref{maxp} to $e^{-tL}f$ for each $t>0$, we have
$$\left(\frac
1{|B|}\int_B | \nabla e^{-tL}(e^{-kr^2L}f)|^{p_0}\right)^{1/p_0}
\le \sum_{j\ge 1} g(j) \left( \frac 1{|2^{j+1}B|} \int_{2^{j+1}B} |\nabla e^{-tL}f|^2 \right)^{1/2}
$$
Squaring this and using that $\sum g(j)<\infty$, we obtain
$$
\int_0^\infty \left(\frac
1{|B|}\int_B | \nabla e^{-tL}(e^{-kr^2L}f)|^{p_0}\right)^{2/p_0}\,  {dt} \le 
C\int_0^\infty \sum_{j\ge 1}  \frac {g(j) }{|2^{j+1}B|} \int_{2^{j+1}B} |\nabla e^{-tL}f|^2
\,  {dt}
$$
Exchanging the sum and the integral, the latter is equal to 
$$
C \sum_{j\ge 1} g(j) \frac 1{|2^{j+1}B|} \int_{2^{j+1}B} |G_L(f)|^2 
$$
which is controlled by $C M( G_L(f)^2)(y)$ for any $y\in B$.

Next, we establish \eqref{pointwise}. Let $B$ be a ball and $r=r(B)$ be its radius. Let $f \in
L^2$. Decompose $f=f_1+f_2+f_3 \ldots$ where $f_j=f\chi_{C_j}$ and $C_j$ are defined as usual.
We start
from
$$
\left(\frac 1{|B|}\int_B |G_L((I-e^{-r^2L})^mf)|^2\right)^{1/2} \le \sum_{j\ge 1} \left(\frac
1{|B|}\int_B |G_L((I-e^{-r^2L})^mf_j)|^2\right)^{1/2}.
$$
For $j=1$, we merely use the $L^2$ boundedness of $G_L$ in \eqref{GLL2bdd} and that of $(I-e^{-r^2L})^m$ 
to obtain
$$
\left(\frac
1{|B|}\int_B |G_L((I-e^{-r^2L})^mf_1)|^2\right)^{1/2}  \le  C\left(\frac
1{|4B|}\int_{4B} |f|^2\right)^{1/2}.
$$
For $j\ge 2$, we write
$$
\frac 1{|B|}\int_B |G_L((I-e^{-r^2L})^mf_j)|^2 = \int_0^\infty \frac 1{|B|}\int_B |\sqrt t \, \nabla
(I-e^{-kr^2L})^m e^{-tL}f_j|^2\, dx \frac {dt}t.
$$
and we use the estimates in the second argument of the proof of Theorem \ref{thrieszp>2} to obtain a bound
$$
\frac {C2^{nj}}{|2^{j+1}B|} \int_{2^{j+1}B} |f|^2 \left(\int_0^\infty \frac{1}{4^{jm}}  \inf \left(
\left(\frac {t} {4^jr^2}\right)^{1/2}, 
 \left(\frac  {4^jr^2}{t}\right)^{m-1/2}\,\, \right) \frac {dt}t\right)^{2},
$$
hence as $m\ge1$,
$$C2^{nj}2^{-4mj} \frac 1{|2^{j+1}B|} \int_{2^{j+1}B} |f|^2.
$$
Choosing further $m>n/4$ allows to sum in $j\ge 2$ and to conclude for \eqref{pointwise}.

\paragraph{Step 2. $p_0>2$ and $(e^{-tL})_{t>0}$ satisfies $L^2-L^{p_0}$ off-diagonal estimates imply
$\|g_L(f)\|_p \lesssim \|f\|_p$ for $2<p<p_0$.}

We apply Theorem \ref{lpp>2} to the sublinear operator $T=g_L$. We choose as usual $A_r=I-(I-e^{-r^2L})^m
$ for $m$ a large enough integral number. 

We have to check \eqref{doma}. Let $B$ be a ball and $r=r(B)$ be its radius and $k=1, \ldots, m$. 
By Minkowski integral inequality and $p_0>2$, we have
$$
\left(\frac 1{|B|}\int_B |g_L(e^{-kr^2L}f)|^{p_0}\right)^{2/p_0} \le \int_0^\infty \left(\frac
1{|B|}\int_B | (tL)^{1/2} e^{-tL}(e^{-kr^2L}f)|^{p_0}\right)^{2/p_0}\, \frac {dt}t.$$
 Using the hypothesis of $p_0$,  and following already used arguments, we
have that 
$$
\left(\frac
1{|B|}\int_B | (e^{-kr^2L}g)|^{p_0}\right)^{2/p_0}
\le \sum_{j\ge 1}   \frac {c_j}{|2^{j+1}B|} \int_{2^{j+1}B} |g|^2 
$$
with $c_j=Ce^{-c4^j}$ for some positive constants $c,C$.
Applying
the commutativity property of the
semigroup and this inequality to $g=(tL)^{1/2}e^{-tL}f$ for each $t>0$, we have
$$\left(\frac 1{|B|}\int_B |g_L(e^{-kr^2L}f)|^{p_0}\right)^{2/p_0}
\le \int_0^\infty \left\{\sum_{j\ge 1}   \frac {c_j}{|2^{j+1}B|} \int_{2^{j+1}B} |(tL)^{1/2}
e^{-tL}f|^2\right\}  \frac {dt} t.
$$ 
As the latter expression equals
$$
 \sum_{j\ge 1}  \frac {c_j}{|2^{j+1}B|} \int_{2^{j+1}B} |g_L(f)|^2 
$$
we obtain a bound in $C M( g_L(f)^2)(y)$ for any $y\in B$ as desired.

Next, we establish \eqref{pointwiseriesz}. Let $B$ be a ball and $r=r(B)$ be its radius. Let $f \in
L^2$. Decompose $f=f_1+f_2+f_3 \ldots$ where $f_j=f\chi_{C_j}$ and $C_j=C_j(B)$ are defined as usual.
We start
from
$$
\left(\frac 1{|B|}\int_B |g_L((I-e^{-kr^2L})^mf)|^2\right)^{1/2} \le \sum_{j\ge 1} \left(\frac
1{|B|}\int_B |g_L((I-e^{-kr^2L})^mf_j)|^2\right)^{1/2}.
$$
For $j=1$, we merely use the $L^2$ boundedness of $g_L$ in \eqref{gLL2bdd} and that of $(I-e^{r^2L})^m$ 
to obtain
$$
\left(\frac
1{|B|}\int_B |g_L((I-e^{-kr^2L})^mf_1)|^2\right)^{1/2}  \le  C\left(\frac
1{|4B|}\int_{4B} |f|^2\right)^{1/2}.
$$
For $j\ge 2$, we write
$$
\frac 1{|B|}\int_B |g_L((I-e^{-kr^2L})^mf_j)|^2 = \int_0^\infty \frac 1{|B|}\int_B |
(tL)^{1/2} e^{-tL}(I-e^{-kr^2L})^m f_j|^2\, dx \frac {dt}t.
$$
As in Section \ref{sec:RTp<2},   one may use the representation \eqref{eqrepresentation} with the  
function $\varphi(z)=(tz)^{1/2}e^{-tz}(1-e^{r^2z})^m$. The corresponding functions 
$\eta_\pm$ satisfy the estimates 
\begin{equation}\label{eqetasf}
|\eta_\pm(z)| \le \frac {Ct^{1/2}}{(|z|+t)^{3/2}} \inf \left(1 , \frac {r^{2m}}{(|z|+t)^m}\right), \quad
z\in
\Gamma_\pm.
\end{equation} Since 
for any $0<\beta< {\frac \pi 2}-\omega$,  $(e^{-zL})_{z \in \Sigma_{\beta}}$ satisfies $L^{2}$
off-diagonal estimates, using 
\eqref{eqrepresentation} and the above estimate for $\eta_\pm$,  $\| (tL)^{1/2} e^{-tL}(I-e^{-kr^2L})^m
f_j\|_{L^2(B)}$ is bounded by 
$$
  C \int_{\Gamma_+} e^{-\tfrac{c4^jr^2} {|z|}}\frac { t^{1/2}}
{(|z|+t)^{3/2}} \frac {r^{2m}}{(|z|+t)^{m}}\, |dz| \ \| f\|_{L^{2}(C_j)}
$$
plus the similar term corresponding to integration on $\Gamma_-$.   Using \eqref{eqintegrale}, this gives us the bound
$$
  \frac{C}{4^{jm}}  \inf \left( \left(\frac {t} {4^jr^2}\right)^{1/2}, 
 \left(\frac  {4^jr^2}{t}\right)^{m}\,\, \right)  \ \| f\|_{L^{2}(C_j)}.
$$
Squaring and integrating with respect to  $t$, we find 
$$\frac 1{|B|}\int_B |g_L((I-e^{-kr^2L})^mf_j)|^2
\le  C2^{jn} 4^{-mj} \frac 1{|2^{j+1}B|} \int_{2^{j+1}B} |f|^2
$$
and this readily implies \eqref{pointwiseriesz}.

\paragraph{Step 3. $p_0<2$ and $(e^{-tL})_{t>0}$ satisfies $L^{p_0}-L^2$ off-diagonal estimates imply
$\|g_L(f)\|_p \lesssim \|f\|_p$ for $p_0<p<2$.  }

We apply Theorem \ref{lpp<2} to $T=g_L$. We choose as usual $A_r=I-(I-e^{-r^2L})^m
$ for $m$ a large enough integral number. 

 Equation \eqref{domap<2} is a direct consequence of the assumption on the semigroup. 
We then turn to the verification of \eqref{pointwisep<2}. Let $B$ be a ball and
$r=r(B)$ its radius and  $j\ge 2$. Let $f$ be a function supported in $B$. We have
$$
\frac 1{|2^{j+1}B|} \int_{C_j} |g_L( (I-e^{-kr^2L})^mf)|^2 = \int_0^\infty
\frac 1{|2^{j+1}B|} \int_{C_j} |
(tL)^{1/2} e^{-tL}(I-e^{-kr^2L})^m f|^2\, dx \frac {dt}t.
$$
Using similar arguments as in the previous step and using the hypothesis on the semigroup,  $\| (tL)^{1/2}
e^{-tL}(I-e^{-r^2L})^m f\|_{L^2(C_j)}$ is bounded by 
$$
  C \int_{\Gamma_+} e^{-\tfrac{c4^jr^2} {|z|}} \frac 1 {|z|^{\gamma/2}}   \frac
{t^{1/2}} {(|z|+t)^{3/2}} \frac {r^{2m}}{(|z|+t)^{m}}\, |dz| \ \| f\|_{L^{p_0}(B)}
$$
plus the similar term corresponding to integration on $\Gamma_-$. Here $\gamma=\gamma_{p_0}=|\frac n 2 -\frac n {p_0}|$.  Using
\eqref{eqintegrale}, this gives us the bound
$$
  \frac{C}{4^{jm}(2^jr)^{\gamma}}  \inf \left( \left(\frac {t} {4^jr^2}\right)^{1/2}, 
 \left(\frac  {4^jr^2}{t}\right)^{m}\,\, \right)  \ \| f\|_{L^{p_0}(B)}.
$$
Squaring and integrating with respect to  $t$, we find 
$$\left(\frac 1{|2^{j+1}B|} \int_{C_j} |g_L( (I-e^{-r^2L})^mf)|^2\right)^{1/2}
\le  C2^{-jn/2} 4^{-{mj}}2^{-j\gamma} \left( \frac 1{|B|} \int_{B} |f|^{p_0}\right)^{1/p_0}
$$
and this readily implies \eqref{pointwisep<2} and $m$ is chosen with $2m+\gamma>n/2$.

\paragraph{Step 4. $p_0<2$ and $(\sqrt t\, \nabla e^{-tL})_{t>0}$ satisfies $L^{p_0}-L^2$ off-diagonal
estimates imply
$\|G_L(f)\|_p \lesssim \|f\|_p$ for $p_0<p<2$.  }

We apply Theorem \ref{lpp<2} to $T=G_L$. We choose as usual $A_r=I-(I-e^{-r^2L})^m
$ for $m$ a large enough integral number. 

 Equation \eqref{domap<2} is a direct consequence of the assumption on the semigroup. 
We then turn to the verification of \eqref{pointwisep<2}. Let $B$ be a ball and
$r=r(B)$ its radius and  $j\ge 2$. Let $f$ be a function supported in $B$. 
We have
$$
\frac 1{|2^{j+1}B|} \int_{C_j} |G_L( (I-e^{-r^2L})^mf)|^2 = \int_0^\infty
\frac 1{|2^{j+1}B|} \int_{C_j} |
\sqrt t\, \nabla e^{-tL}(I-e^{-kr^2L})^m f|^2\, dx \frac {dt}t.
$$
   Proposition \ref{propequivgrafdz} shows that we may replace the hypothesis by
$(\sqrt z\, \nabla e^{-zL})_{z\in \Sigma_\beta}$ satisfies $L^{p_0}-L^2$ off-diagonal estimates for any 
$\beta \in (0,{\frac \pi 2}-\omega)$ (up to changing $p_0$ to an arbitrary larger value). Using then similar arguments as in the step 1,  $\| \sqrt t \, \nabla
e^{-tL}(I-e^{-kr^2L})^m f\|_{L^2(C_j)}$ is bounded by 
$$
  C \int_{\Gamma_+} e^{-\tfrac{c4^jr^2} {|z|}} \frac 1 {|z|^{\gamma/2+1/2}}   \frac
{t^{1/2}} {(|z|+t)} \frac {r^{2m}}{(|z|+t)^{m}}\, |dz| \ \| f\|_{L^{p_0}(B)}
$$
plus the similar term corresponding to  integration on $\Gamma_-$. Here $\gamma=\gamma_{p_0}=|\frac n 2 -\frac n {p_0}|$.  Using
\eqref{eqintegrale}, this gives us the bound
$$
  \frac{C}{4^{jm}(2^jr)^{\gamma }}  \inf \left( \left(\frac {t} {4^jr^2}\right)^{1/2}, 
 \left(\frac  {4^jr^2}{t}\right)^{m-1/2}\,\, \right)  \ \| f\|_{L^{p_0}(B)}.
$$
Squaring and integrating with respect to  $t$, we find 
$$\left(\frac 1{|2^{j+1}B|} \int_{C_j} |G_L( (I-e^{-kr^2L})^mf)|^2\right)^{1/2}
\le  C2^{-jn/2} 4^{-{mj}}2^{-j\gamma} \left( \frac 1{|B|} \int_{B} |f|^{p_0}\right)^{1/p_0}
$$
and this readily implies \eqref{pointwisep<2} if $m$ is chosen with $2m+\gamma>n/2$.

\paragraph{Step 5. Reverse $L^p$ inequality for $g_L$ when $p_-(L)<p<p_+(L)$.  }
By functional calculus for $L$, we have the representation formula for the identity on $L^2$
\begin{equation}\label{eqrepidentity}
f= 2\int_0^\infty tL e^{-2tL} f\frac {dt}t.
\end{equation}
Hence, dualizing against $g$ in $L^2$, writing
$$\langle tL e^{-2tL} f, g\rangle = \langle (tL)^{1/2} e^{-tL} f, (tL^*)^{1/2} e^{-tL^*} g\rangle$$
and using Cauchy-Schwarz inequality yield
$$ \left|\int_{\RR^n} f \ol g\right| \le \int_{\RR^n} g_L(f) g_{L^*}(g).
$$
Recall that $\mathrm{int}\calJ(L)=(p_-(L),p_+(L))$ is the largest open interval of exponents $p$ in $(1,\infty)$ for
which  the semigroup
$(e^{-tL})_{t>0}$ is
$L^p$ bounded. Then, by steps 2 and 3,  $g_L$ is $L^p$ bounded for $p \in \mathrm{int}\calJ(L)$.
Let $p \in \mathrm{int}\calJ(L)$. Applying this to $L^*$  since $p' \in \mathrm{int}\calJ(L^*)$, we have  
$$
\|g_{L^*}(g)\|_{p'} \le C \|g\|_{p'},
$$
hence,
$$ \left|\int_{\RR^n} f \ol g\right| \le C\|g_L(f)\|_p \|g\|_{p'}.$$
In conclusion, 
$$
\|f\|_p \le  C\|g_L(f)\|_p$$
for $p \in \mathrm{int}\calJ(L)$.

\paragraph{Step 6. 
$\|g_L(f)\|_p \sim \|f\|_p$ implies $(e^{-tL})_{t>0}$  $L^p$ bounded. }

It follows easily form the definition of $g_L$ and the commutativity properties of the semigroup
that for all $s>0$,
$$
g_L(e^{-sL}f) \le g_L(f)
$$
in the pointwise sense, hence in  $L^p$ norm. In particular, this and the hypothesis yield
$\|e^{-sL}f\|_p \le C \|f\|_p$ uniformly in $s>0$.

\paragraph{Step 7. $\|G_L(f)\|_p \lesssim \|f\|_p$  implies $(\sqrt t\, \nabla e^{-tL})_{t>0}$
$L^q$-bounded for $q$ between $2$ and $p$.}
\footnote{\, In practice, quadratic functionals are built
from operators whose individual boundedness properties are known. The converse, thus, is never
considered, so that this result appears new. 
Note that the argument could be written for arbitrary
quadratic functionals made after an analytic family of operators. We leave to the reader the care of
stating this general result.}

The argument of step 5 does not apply for the semigroup and the gradient operator do not commute.
We rely instead on some trick using complex interpolation.
Define 
$$
\widetilde G_{L}(f)(x)= \left( \int_0^\infty |(\nabla (tL)e^{-tL}f)(x)|^2 \, {dt}\right)^{1/2}.
$$
We assume $p>2$. The argument for $p<2$ is entirely similar. We  show that $\|G_L(f)\|_p \lesssim
\|f\|_p$ implies 
$\|\widetilde G_L(f)\|_q \lesssim \|f\|_q$ for $2<q<p$.
Assume this is done. 
We prove the $L^q$ boundedness of $\sqrt s\, \nabla e^{-sL}$. Without loss of generality, we may asume
$s=1$.
Write
\begin{align*}
|\nabla e^{-L} f| &\le  \int_1^2 |
\nabla e^{-sL} f| \, ds +\int_1^2 |\nabla e^{-L} f - \nabla e^{-sL} f| \, ds
\\
&\le  \int_1^2 |
\nabla e^{-sL} f| \, ds + \int_1^2 \left| \int_{1}^s \nabla Le^{-tL}f \, dt\right| \,ds
\\
&
\le \int_1^2 |
\nabla e^{-sL} f| \, ds + \int_1^2  |\nabla (tL)e^{-tL}f| \, dt.
\end{align*}
Hence,
$$
|\nabla e^{-L} f| \le G_L(f) + \widetilde G_L(f).
$$
and the $L^q$ boundedness of $\nabla e^{-L}$ follows from that of $G_L$ 	and $\widetilde G_L$. 

It remains to prove the $L^q$ boundedness of the latter. To this end, we follow the proof of Stein's
complex interpolation theorem after dualizing. Fix
$q$ with
$2<q<p$.

Let $f$ be a $\CC$-valued simple function  with $\|f\|_q=1$ and $g=(g_t)_{t>0}$ be a $H$-valued  simple
function where $H=L^2((0,\infty), \CC^n;  {dt} )$ with $\|g\|_{q'}=1$. Write
$f(x)=\sum_k a_k \chi_{E_k}(x)$ and $g(t,x)=\sum_k b_k(t) \chi_{E_k}(x)$ where $E_k$ are pairwise
disjoint measurable sets,
$a_k$ are complex numbers and $b_k(t)$ are $\CC^n$-valued, and  set $B_k= \left(\int_0^\infty |b_k(t)|^2
 {dt} \right)^{1/2}$.

Let 
$$f_z(x)= \sum_k |a_k|^{\alpha(z)} \frac {a_k}{|a_k|} \chi_{E_k}(x)
$$
and 
$$
g_z(t,x) = \sum_k {B_k}^{\beta(z)} \frac {b_k(t)}{B_k} \chi_{E_k}(x)$$
 with 
$$\alpha(z)= \frac q p (1-z) + \frac q 2 z \quad \mathrm{ and} \quad \beta(z)=\frac {q'}{p'} (1-z) + \frac {q'} 2 z.$$

 Pick $\beta \in (0,
{\frac \pi 2}-\omega)$ and  consider the function 
$$
F(z)= \int_{\RR^n}\int_0^\infty \nabla  e^{-te^{i\beta z}L}f_z(x)\cdot g_z(t,x) \, dx {dt}.
$$
defined for $z$ in the strip  $0\le \Re z \le 1$. This function is clearly  continuous, and it is 
analytic in the interior  of this  strip. Moreover, for $z=iy$ with $y\in \RR$, one finds easily from 
 H\"older inequalities and the   change of variable $te^{-\beta y} \mapsto t$, 
$$
|F(iy)| \le \|G_{e^{-\beta y }L}(f_{iy})\|_p \|g_{iy}\|_{p'} \le C_0 e^{\beta y /2} \|f_{iy}\|_p\|g_{iy}\|_{p'} = C_0 e^{\beta y /2} 
$$
and for $z=1+iy$ with $y\in \RR$, by the square function estimate \eqref{GLL2bdd} for $e^{i\beta } L$ and the same change of variable,
$$
|F(1+iy)| \le \|G_{e^{i\beta (1+iy)}L}f_{1+iy}\|_2 \|g_{1+iy}\|_2 \le C_1 e^{\beta y /2}  \|f_{1+iy}\|_2\|g_{1+iy}\|_2=C_1 e^{\beta y /2} .
$$
Here, $C_{1}$ depends on $\beta$ in a non decreasing manner. 
Hence,  the extension of three lines theorem in \cite[Chapter V, Lemma 4.2]{SW} yields for $0<x<1$ and $y\in \RR$, 
$$
|F(x+iy)| \le C_0^{1-x}C_1^x e^{\beta y /2} .
$$
This implies that  $L^q$ boundedness of $G_{e^{i\beta c_q} L}$ when $\Re c_q= \frac {\frac 1 q - \frac 1
p } {\frac 1 2 -\frac 1 p}$ with a bound  controlled by $C e^{\beta \Im c_{q} /2} $. Now, this is true for all
$\beta$ in $(0,{\frac \pi 2} -\omega)$ and also  in $(-{\frac \pi 2}+\omega, 0)$ by changing   $\beta$ to $-\beta$.

Keep $q$ as before and let $f$, $g$ be simple functions as above. If $\beta$ is small enough,
the function 
$$G(z)= \int_{\RR^n}\int_0^\infty \nabla  e^{-te^{i\beta z}L}f(x)\cdot g(t,x) \, dx{dt}
$$
is  continuous in the strip $-1\le \Re z \le 1$ and, by the previous argument, analytic in the open strip
minus the real axis. Thus it is analytic in the open strip by Morera's theorem. Applying the three lines
theorem again (this times, $f$ and $g$ are fixed and the exponent remains $q$) or Cauchy's theorem about 0 (since we are not after optimal bounds), we deduce a  bound for
$G'(0)$, which is equivalent to the
$L^q$ boundedness of
$\widetilde G_L$.

\paragraph{Step 8. Reverse $L^p$ inequality for $G_L$ when $1<p<+\infty$. }

What we have proved so far applies to any operator $L$ in our class, and in particular, 
to $L=-\Delta$. The explicit formula for the heat kernel implies that $p_-(-\Delta)=1$ and $q_+(-\Delta)=\infty$.
Hence, we recover from our method the well-known estimate~\footnote{\, due to
Stein. See \cite{St1}.}
$$
\|G_{-\Delta} (f)\|_p \lesssim \|f\|_p
$$
for all $1<p<\infty$.

Now, let $f,g \in L^2$ and observe that~\footnote{\, we learned this trick from J.
Pipher in an unpublished manuscript. See \cite{AMcN1} where it is used.} 
\begin{align*}
\int_{\RR^n} f \, \ol g &= \lim_{\ep\downarrow 0} \int_{\RR^n} e^{-\ep L}f\, \ol{e^{\ep\Delta}g} -  \lim_{R\uparrow \infty} \int_{\RR^n} e^{-R L}f\,
\ol{e^{R\Delta}g}
\\
&
= - \int_0^\infty \frac d {dt} \int_{\RR^n} (e^{-t L}f)(x)\, \ol{(e^{t\Delta}g)(x)}\, dx\, dt
\\
&
= \iint_{\RR^n \times (0,\infty)} (A(x)+I) (\nabla e^{-tL}f)(x) \cdot 
 \ol{(\nabla e^{t\Delta}g)(x)}\, dxdt.
\end{align*}
The last equality is 
 obtained by integration by parts in the $x$ variable after computing the time
derivative. Hence, we obtain
$$
\left| \int_{\RR^n} f \, \ol g\, \right| \le (\|A\|_\infty +1) \int_{\RR^n} G_L(f) G_{-\Delta}(g).
$$
Thus, if $1<p<\infty$,  the $L^{p'}$ boundedness of $G_{-\Delta}$ yields
$$
\left| \int_{\RR^n} f \, \ol g\, \right| \lesssim   \|G_L(f)\|_p \|g\|_{p'}
$$
and it follows
$$
\|f\|_p \lesssim \|G_L(f)\|_p.$$
Of course this inequality is meaningfull whenever the right hand side is finite. 
\end{proof}

Let us draw some consequences of our results.

\begin{corollary} If $n=1$, we have for $1<p<\infty$,  $$
\|G_L(f)\|_p \sim \|f\|_p \sim \|g_L(f)\|_p.
$$
\end{corollary}

\begin{corollary}  If $n=2$, we have $\|g_L(f)\|_p \sim \|f\|_p$ for $1<p<\infty$ and 
$\|G_L(f)\|_p \sim \|f\|_p$ for $1<p<2+\ep'$. 
\end{corollary}

\begin{corollary} If $n\ge 3$, we have $\|g_L(f)\|_p \sim \|f\|_p$ for $\frac{2n}{n+2} -\ep<p<\frac{2n}{n-2} +\ep'$ and  
$\|G_L(f)\|_p \sim \|f\|_p$ for $\frac{2n}{n+2} -\ep<p<2+\ep'_1$. 
\end{corollary}

\begin{corollary}  If $n\ge 3$ and $L$ has \textbf{real} coefficients, then we have $\|g_L(f)\|_p \sim \|f\|_p$ for $1<p<\infty$
and $\|G_L(f)\|_p \sim \|f\|_p$ for
$1<p<2+\ep$.
\end{corollary}

\begin{corollary} If $n\ge 5$, there exists $L$ such that $g_L$ and $G_L$ are not bounded for some $p$
close to 1 (hence $g_{L^*}$ is not bounded on $L^p$ for some $p$ close to $\infty$).
\end{corollary}

We finish this section by proving Lemma \ref{lemmafc}.

\paragraph{Proof of Lemma \ref{lemmafc}:} The assumption of the lemma is that 
$(e^{-tL})_{t>0}$ satisfies $L^\rho-L^2$ off-diagonal estimates. 
By duality, $(e^{-tL^*})_{t>0}$ satisfies $L^2-L^{\rho^\prime}$ off-diagonal estimates. 
Hence, it follows from the method of step 2 above (applied to discrete times $4^k$, $k\in \ZZ$, and changing
$(tL)^{1/2}e^{-tL}$ to $\psi(tL)$ with $\psi$ given by \eqref{eqpsi}) that for $2<q'<\rho'$,
$$
\left\| \left( \sum_{k \in \ZZ} |\psi(4^kL^*)f|^2 \right)^{1/2}\right\|_{q'} \le C \|f\|_{q'}.
$$
By duality, we obtain \eqref{eq22}.

\subsection{On inequalities of Stein and Fefferman for non-tangential square functions} 
 
Consider the functional
$$g_\lambda^*(f)(x)= \left(\iint_{\RR^{n+1}_+} \left(\frac
t{|x-y|+t}\right)^{n\lambda} |t\nabla u(y,t)|^2 \frac
{dydt}{t^{n+1}}\right)^{1/2}
$$
where $u(x,t)$ is the harmonic extension of $f$ and $1<\lambda$. 
It is  bounded on $L^p$ ($1<p<\infty$) if and only if
$\lambda> \frac 2 p$.~\footnote{\, This   is due to Stein \cite{St1}.} At
the critical case
$\lambda=\frac 2 p$,  it is weak-type $(p,p)$.~\footnote{\, This is in C. Fefferman's thesis
\cite{Fef}.} Of course, the main tool for this is  Calder\'on-Zygmund decomposition for $L^p$ functions.

We show how to use our technology in this situation by
replacing functions of the Laplacian  by functions of $L$. Again, the main point is that $L^p$ boundedness
of the semigroup suffices. We also separate the cases $p<2$ and $p>2$.

Define
$$g_\lambda^*(f)(x)^2
= \iint_{\RR^{n+1}_+} \left(\frac
{\sqrt t}{|x-y|+\sqrt t}\right)^{n\lambda} |\calT (e^{-tL}f)(y)|^2 \frac
{dydt}{{\sqrt t\, }^{n}},
$$ 
where for fixed $t>0$, 
$$
|\calT g(y,t)|^2 = |\nabla_y\, g(y,t)|^2 + |L^{1/2}g(y,t)|^2.
$$
In particular, this square function controls non tangential square functions
where integration is performed on parabolic cones $|x-y| \le c\sqrt t$:
$$
g_{NT}(f)(x)^2=\iint_{|x-y| < c\sqrt t}  |\calT e^{-tL}f(y)|^2 \frac
{dydt}{{\sqrt t\, }^{n}}.
$$

\begin{proposition}\label{propsteinp<2} If $q_-(L)<p<2$ and $\lambda = \frac 2 p$ then $g_\lambda^*$  has
weak type (p,p).
\end{proposition}

\begin{proposition}\label{propsteinp>2} If $2<p<q_+(L)$ and $\lambda >1$ then $g_\lambda^*$  has
strong  type (p,p).
\end{proposition}

\begin{corollary} If $q_-(L)<p<q_+(L)$ and $\lambda > \frac 2 p$, then  $g_{NT}, g_\lambda^*$ are bounded on $L^p$ and one
has~\footnote{\,  While Le Merdy's theorem mentioned above works well for the vertical square function of
abstract operators on
$L^p$ spaces, it is not clear it applies to non-tangential square functions which are more geometrical
objects in an abstract setting.}
$$\|g_\lambda^*f\|_p \sim \|g_{NT}(f)\|_p \sim \|f\|_p.$$
\end{corollary} 

Note that  the result  for  $g_{NT}$ holds for $q_-(L)<p<q_+(L)$ since it is independent of $\lambda$.  The proof of the corollary is simple: we have a pointwise control
of $g_{NT}$ by $g_\lambda^*$ for any $\lambda>\frac 2 p $ and we obtain boundedness. The reverse inequalities are
 obtained as for $g_L$ and $G_L$ following the heuristic  idea below.
Suppose we can write $\int f(y) g(y)\, dy = \iint f_t(y) g_t(y) dydt$ and that 
$\int h(x,y,t)\, dx= 1$ for all $y,t$. Then
$$
\int f(y) g(y)\, dy = \int \left( \iint f_t(y) g_t(y) h(x,y,t)\, dydt\right) dx.
$$
It suffices to apply Cauchy-Schwarz in the variables $y,t$ to obtain non-tangential
square function by choosing $h$ and then H\"older in the $x$ variable to deduce 
reverse inequalities from direct ones. We skip further details.
 
\begin{remark} The limitations on $p$ are only due to the presence of the spatial gradient $\nabla_y$ in
the definition of $g_\lambda^*$. If one drops this gradient to keep only the  $L^{1/2}$ part  then the
range of $p$ becomes $p_-(L)<p<p_+(L)$.
\end{remark}

Let us turn to the proof of  Proposition \ref{propsteinp<2}.
  Due to the fact that $x$ and $y$ may be far apart, Theorem \ref{lpp<2}
does not apply directly and one has to do some transformations. The key
identity of this proof is that   for any closed set $F$,
$$
\int_F g_\lambda^*(f)(x)^2\, dx= \iint_{\RR^{n+1}_+}  |\calT
(e^{-tL}f)(y)|^2 J_{\lambda,F}(y,t) \, {dydt}
$$
with  
$$
J_{\lambda,F}(y,t)=\frac 1 {{\sqrt t \,}^n}\ \int_F \left(\frac
{\sqrt t}{|x-y|+\sqrt t}\right)^{n\lambda}\, dx.
$$
First
 $J_{\lambda,F} \le C$  so that  
$$
\| g_\lambda^*(f)\|_2^2 \le C \|g_L(f)\|_2^2 + C\|G_L(f)\|_2^2
$$
where $g_L$ and $G_L$ are the square functions defined earlier.
The $L^2$ boundedness of  $g_L$ and $G_L$ implies
the $L^2$ boundedness of $g_\lambda^*$.

Second, we have also  if $y$ lies in some cube $Q$ and $F=\RR^n\setminus (2Q)$
that
$$ J_{\lambda,F}(y,t) \le C \sqrt t\,^{n(\lambda-1)} |Q|^{-(\lambda -1)/2}.$$

We begin as in Theorem \ref{lpp<2}, by looking at $\{ g_\lambda^*(f)>\alpha\}$ and
decomposing $f=g+\sum b_i$ according to the threshold $\alpha^p$ for $|f|^p$. For
$g$  use the $L^2$ boundedness of $g_\lambda^*$. 

Next write again $$b_i=A_{r_i}b_i +
(1-A_{r_i})b_i$$ where $A_r$ is the operator that works for $g_L$ and $G_L$ in the previous section.
The term $\sum A_{r_i}b_i$ is again in $L^2$ with the right bound so that the
$L^2$ boundedness of $g_\lambda^*$ suffices again. It remains to estimate the size of  the
set 
$$\{x\in \RR^n; g_\lambda^*(\sum(1-A_{r_i})b_i)(x)>\alpha/3\}.$$
 Again we
take away the union of the dilated Whitney cubes $4Q_i$, whose mass is under control. 
It remains to estimate what is left on  its complement $F$. By Tchebytchev's inequality,
it is enough to estimate $\int_F g_\lambda^*(\sum(1-A_{r_i})b_i)(x)^2\, dx$
which we rewrite as
$$
\iint_{\RR^{n+1}_+}  |\sum_i \calT
(e^{-tL}(1-A_{r_i})b_i)(y)|^2
J_{\lambda,F}(y,t) \, {dydt}. 
$$
The non local part of the ith summand is when $y \notin 2Q_i$. We  bound
$J_{\lambda,F}$ by a constant and we are back to the estimates performed to
obtain the weak type $(p,p)$ for $g_L$ and $G_L$. We refer the reader to steps 
3 and 4 in the previous section.

It remains to localise each ith summand on $2Q_i$.
By the bounded overlaps  of $2Q_i$'s~\footnote{\, We have to make sure in the construction of the Whitney
cubes that  this actually holds and may be twice the cubes is not appropriate but certainly $cQ_i$
with some $c>1$ is.}  and the
second upper bound on
$J_{\lambda, F}$, we have an upper bound  
$$
N \sum_i \iint_{2Q_i\times \RR^+} \frac{\sqrt t\,\,^{n(\lambda-1)}}{
|Q_i|^{(\lambda -1)/2}} | \calT (e^{-tL}(1-A_{r_i})b_i)(y)|^2  \,
{dydt}.
$$
Now for each $i$, we integrate on the full upper half space: the integral in
$y$ (with $t$ fixed) and the solution of the Kato problem~\footnote{\, One can proceed
also using Remark \ref{remarkKL}.} allow us to  bound the term with $\nabla_y$
by the one with
$L^{1/2}$. Next if $a=\frac n 2(\lambda-1)= \frac  n p - \frac n 2$, the $L^p$ bounds of the
semigroup implies   by Proposition \ref{propHLS}   the Hardy-Littlewood-Sobolev inequality
$$
\|L^{-a}f\|_2  \le C \|f\|_p.
$$
Using the square function estimate of McIntosh-Yagi based on $(t \, L)^a 
e^{-tL}$ we obtain that the ith term is bounded by
$$
C |Q_i|^{1-2/p} \|L^{-a}(1-A_{r_i})b_i\|_2^2 \le C |Q_i|^{1-2/p} \|(1-A_{r_i})b_i\|_p^2 \le C |Q_i|^{1-2/p}
\|b_i\|_p^2. 
$$ 
It remains to sum other $i$ and we are done.  \qed 
\bigskip

Let us turn to the proof of  the Proposition \ref{propsteinp>2}.  Here also, we cannot apply
directly Theorem \ref{lpp>2} but rather its spirit and its proof.
We let $f\in  L^2(\RR^n)$ and $B$ be a ball with radius $r$. We also let $2<p_0<q_+(L)$ and $m>n/4$ an
integer. Assume that we have proved that for $k=1,2, \ldots, m$ 
\begin{equation}\label{eqsteinp>2}
\left(\frac{1}{|B|}\int_B |g_\lambda^*(e^{-kr^2L}f)|^{p_0}\right)^{1/p_0} \le C \inf_{x\in B}
M(g_\lambda^*(f)^2)^{1/2}(x)
\end{equation}
then we can argue as follows. For $A_r=I-(I-e^{r^2L})^m$, we have
$$ 
g_\lambda^*(f)^2(x) \le 2 g_\lambda^*(A_rf)^2(x) +  2 g_\lambda^*((I-A_r)f)^2(x).
$$
 Write
\begin{align*}
g_\lambda^*((I-A_r)f)^2(x)&= \iint_{y\in 2B} h(x-y,t) |\calT (e^{-tL}(I-A_r)f)(y)|^2 dydt
\\
&+ 
\iint_{y\notin 2B} h(x-y,t) |\calT (e^{-tL}(I-A_r)f)(y)|^2 dydt
\end{align*}
with $h(x,t)=\sqrt t\, ^n \left( \frac {\sqrt t} {|x|+\sqrt{t}}\right)^{n\lambda}.$
As $h(x-y,t) \sim h(z-y,t)$ for $x,z \in B$ and $y\notin 2B$, the second integral is
bounded by 
$$2 g_\lambda^*(A_rf)^2(x) + 2\inf_{z\in B} g_\lambda^*(f)^2(z).$$ 
Hence, we can apply Proposition \ref{goodlambda} as in the proof of Theorem \ref{lpp>2} with
$$
G_B(x) = 2\iint_{y\in 2B} h(x-y,t) |\calT (e^{-tL}(I-A_r)f)(y)|^2 dydt 
$$
and 
$$H_B(x)= 4 g_\lambda^*(A_rf)^2(x) + 2\inf_{z\in B} g_\lambda^*(f)^2(z)
$$ 
provided we show that 
$$
\frac{1 }{|B|} \int_B G_B \le C\inf_{z\in B} M(|f|^2)(z).
$$
But,  $\int_B h(x-y,t)\, dx\le 1$ since $\lambda>1$, hence
$$
\int_B G_B(x)\, dx \le 2 \iint_{2B\times \RR^+}  |\calT (e^{-tL}(I-A_r)f)(y)|^2 dydt
$$
and we are back to the calculations made in steps 1 and 2 of the previous section that give us
$$
\iint_{2B\times \RR^+}  |\calT (e^{-tL}(I-A_r)f)(y)|^2 dydt \le C\inf_{z\in B} M(|f|^2)(z).
$$

Hence, it remains to establish \eqref{eqsteinp>2}. To do that, assume that $B$ is the unit
ball and $r=1$ (One can treat the general case by rescaling and translation, which changes $L$ to another
operator with the same properties and the same critical numbers). Assume also for simplicity that
$k=1$. 
We write for $x\in B$,
$$g_\lambda^*(e^{-r^2L}f)^2(x) = I + II
$$
where
$$I= \iint_{(y,t) \notin E} h(x-y,t) |\calT (e^{-tL}e^{-r^2L}f)(y)|^2 dydt
$$
and 
$$
II= 
\iint_{(y,t) \in E} h(x-y,t) |\calT (e^{-tL}e^{-r^2L}f)(y)|^2 dydt
$$
and $(y,t)\notin E$ means $y\in 2B$ and $t\le r(2B)^2$. Let us treat the first term.

Using Minkowksi integral inequality (with respect to $t$) and  $\int_{2B} h(x-y,t) \, dx\le 1$ whenever
$y\in 2B$, we obtain
$$
\left(\int_B I^{p_0/2} \right)^{2/p_0}  \le \int_0^4 \left( \int_{2B} |\calT
(e^{-r^2L}e^{-tL}f)(y)|^{p_0} dy\right)^{2/p_0} dt.
$$
Following again the calculations in steps 1 and 2 of the previous section
$$
\left( \int_{2B} |\calT
(e^{-tL}(I-A_r)f)(y)|^{p_0} dy\right)^{2/p_0} \le \sum_{j\ge 2} \frac {c_j} {|2^{j+1}B|}\int_{
2^{j+1}B}  |\calT
(e^{-tL}f)(y)|^2 dy
$$ with $c_j=Ce^{-c4^j}$ and with the limitation $p_0<q_+(L)$ from $\nabla_y$ and $p_0<p_+(L)$ from
$L^{1/2}$ in the definition of $\calT$.
Now since $t\le 4$ and $j\ge 2$, 
$$\int_{
2^{j+1}B}  |\calT
(e^{-tL}f)(y)|^2 dydt \le \frac A {\sqrt t\, ^n } \int_{
2^{j+1}B}  \int_{|x-y| \le \sqrt t} |\calT
(e^{-tL}f)(y)|^2 dy\ dx
$$
for some $A>0$. Hence, we obtain readily
$$
\left(\int_B I^{p_0/2} \right)^{2/p_0} \le AC \inf_{x\in B} M((g_{NT}f)^2)^{1/2}(x).
$$

For the second term $II$, we first observe that 
for $x\in B$ and $(y,t)\in E$, then $h(x-y,t)\sim h(y,t)$ (recall that $B$ is the unit ball so that $0$ is
its center), thus $\left(\int_B II^{p_0/2} \right)^{2/p_0} \le \sup_B II$. Next, decompose $E$  as the
union of
$E_k$, $k\ge -1$, as follows: $E_{-1}= E \cap\{ |y] < \sqrt t\}$ and $E_k=E\cap \{2^k\sqrt t < |y| \le
2^{k+1}\sqrt t\}$. Then for $(y,t) \in E_k$, $h(y,t) \le 2^{-kn\lambda} \sqrt t\, ^{-n}.$  
 A crucial geometrical observation is that if $(y,t) \in E_k$ then  $2^k\sqrt t \ge 1$.
Using the method of the second argument  in the proof of Theorem \ref{thrieszp>2} (recall that $r=1$), 
\begin{align*}
\int_{|y| \le 2^{k}\sqrt t} |\calT e^{-r^2L}
(e^{-tL}f)(y)|^2 dy
&\le C\int_{|z| \le 2^{k+1}\sqrt t} |\calT
(e^{-tL}f)(z)|^2 dz
\\
&
+
\sum_{j\ge k+2} \int_{|z| \le 2^{j}\sqrt t} Ce^{-c4^jt}|\calT
(e^{-tL}f)(z)|^2 dz.
\end{align*}
Hence
$$II \le 
\iint m(z,t)  |\calT
(e^{-tL}f)(z)|^2 dzdt
$$
where
$$m(z,t)= C\sum_{k\ge -1} 2^{-kn\lambda} \sqrt t \, ^{-n} 1_{2^k\sqrt t \ge  1}( 1_{|z|\le 2^{k+1} \sqrt t
} + 
\sum_{j\ge k+2} 
 e^{-c4^jt} 1_{|z|\le 2^{j}\sqrt t}).
$$
Tedious but elementary verifications show that $m(z,t) \le C h(x-z,t)$ for all
$x\in B$ and $(z,t) \in \RR^{n+1}_+$ using only $\lambda >0$. Hence,
$II \le C \inf_{x\in B} g_\lambda^*(f)^2(x)$. \qed

\section{Miscellani}

\subsection{Local theory}

Let $L$ be as in the Introduction.  We have developed a global
(or homogeneous) $L^p$ theory by making global in time assumptions on the semigroup. 
Reasons for this theory not to apply to a particular $L$  at $p$ are that the
semigroup is not $L^p$ bounded for some (or all) $t>0$ (in which case this is
the end of the story) or that the semigroup operators are bounded on $L^p$ but
not uniformly, ususally with an exponential blow up. In the second case, adding a
large
$s$ to
$L$ gives us back the uniformity. The local 
$L^p$ theory consists in working with $L+s$ instead of $L$.   Hence, the above
results may be adapted with minor modifications in the proofs by changing 
systematically
$L$ to $L+s$ for
$s>0$.  One may define the four critical exponents
$p_\pm(L+s)$ and $q_\pm(L+s)$ which may depend on $s$ or not. Indeed, the $L^2$ theory developed in Section \ref{sec:basicL2} works 
with $s=0$ and the numbers $p_-(L+s)$, $q_-(L+s)$ are non increasing, and
$p_+(L+s)$,  $q_+(L+s)$ non decreasing as $s$ grows.

We have the following assertions for $s>0$.
\begin{enumerate}
\item $p_-(L+s)=q_-(L+s)$ and $(q_+(L+s))^* \le p_+(L+s)$.
\item $q_+(L+s)$ is the supremum of exponents $p$ for which one has the invertibility of 
$L+s$ from $W^{1,p}$ onto $W^{-1,p}$. 
\item For the Riesz transform $\nabla (L+s)^{-1/2}$ the range for $L^p$ boundedness is the open interval
$(p_-(L+s), q_+(L+s))$. For $p$ in this range $L+s$ is an isomorphism from 
$W^{1,p}$ onto $L^p$. 
\item There is bounded holomorphic functional calculus for $L+s$ on $L^p$ essentially when $ p_-(L+s)<p<p_+(L+s)$.
\item There is a Hardy-Littlewood-Sobolev-Kato diagram.
\item There is an  equivalent $L^p$ norm defined by 
$g_{L+s}$ essentially for $ p_-(L+s)<p<p_+(L+s)$.
\item There is an  equivalent $L^p$ norm defined by 
$G_{L+s}$ essentially for $ q_-(L+s)'<p<q_+(L+s)$.
\end{enumerate}

This applies to operators whose coefficients have some smoothness. If the coefficients are, in addition,  BUC (bounded uniformly continuous) 
or in the closure of $BUC$ for the $bmo$ norm, then it is known that 
$p_-(L+s) =q_-(L+s)=1$ and $p_+(L+s) =q_+(L+s)=\infty$  
for $s$ large enough.~\footnote{\, This is a consequence of \cite{AMT} (see also
\cite{IS})}  This gives  $L^p$ estimates for Riesz transforms, functional
calculi, square functions in the range 
$1<p<\infty$.

\bigskip

One can also add to $L$ perturbation by lower order terms and develop the similar theory. 

\bigskip

One can probably develop this theory for operators on domains with  Lipschitz boundaries at least with 
Dirichlet or Neumann boundary conditions. This is left to the interested reader.  
\bigskip

Another interesting direction is to test this theory for other classes of
elliptic operators such as Schr\"odinger operators for which criteria for the
determination of $p_-(L)$ and $p_+(L)$ have been given.~\footnote{\, 
 See \cite{LSV}} This theory already applies in the range $p<2$.
~\footnote{\,\cite{BK2} and \cite{BK3} for results on the functional calculus and
Riesz transforms} It remains to study the range $p>2$.

\subsection{Higher order operators and systems}\label{sec:higherorder}

Consider   an homogeneous elliptic operator $L$ of order $m$, $m\in \NN$, $m\ge 2$, defined by
\begin{eqnarray}\label{eqLhomho}
Lf= (-1)^{m}\sum_{|\alpha|=|\beta|=m} 
\partial^\alpha (a_{\alpha\beta} \partial^\beta f ),
\end{eqnarray}
where the  coefficients  $a_{\alpha\beta}$ are  
complex-valued  $L^\infty$  functions on $\RR^n$,
and we assume
 \begin{equation}\label{eqLhombounded}
\left| \sum_{|\alpha|=|\beta|=m} \int_{\RR^n} a_{\alpha\beta}(x) \partial^\beta
f(x) 
\partial^\alpha \bar g(x) \, dx \right| \le \Lambda \| \nabla^m f\|_2 \|
\nabla^m g\|_2
\end{equation}
and the strong G\aa rding inequality   
\begin{equation}\label{eqstronggarding}
\R \sum_{|\alpha|=|\beta|=m} \int_{\RR^n} a_{\alpha\beta}(x) \partial^\beta
f(x) 
\partial^\alpha \bar f(x) \, dx \ge \lambda \| \nabla^m f\|_2^2
\end{equation}
for some  $\lambda>0$ and $\Lambda <+\infty$ independent of $f,g
\in  W^{m,2}$. Here, $\nabla^k$ is
the array of all $k$th order derivatives.

\bigskip

One can also generalize second order or higher order operators to elliptic systems of
any  even order verifying the strong G\aa rding inequality.  For simplicity of
exposition we stick to the scalar case but all works similarly for systems. 

\bigskip

The $L^2$ theory for the semigroup is analogous. There are  
\begin{enumerate}
\item bounded holomorphic functional
calculus on $L^2$,
\item  
$L^2$ off-diagonal estimates for the families $( t^{k/2m}\,  \nabla^k e^{-tL})_{t>0}$ and their analytic extensions for $0\le k \le m$ where the gaussian decay
$e^{-cu^2}$ is changed to
$e^{-c u^{\frac {2m} {2m-1}}}$ and  the homogeneity changes from $\sqrt t$ to $t^{1/2m}$.  
\end{enumerate}
These estimates yield the generalized conservation property
$$e^{-tL}P=P$$ for all $t>0$ in the $L^2_{loc}$ sense for $P$ polynomial of degree less 
than $m$.~\footnote{\, See  \cite{AHMcT} for a proof under additional hypotheses.
The argument has the same structure in the general case.} Since
$L$ is constructed as before as a maximal-accretive operator, it has a square root and    
one has in all dimensions~\footnote{\, This is   \cite{AHMcT}, Theorem 1.1.}
\begin{equation}
\label{eqkato}
\| {L^{1/2}} f\|_2
\sim \| \nabla^m f\|_2.
\end{equation}
Moreover, the square functions $g_L$ and $G_L$ define equivalent norms on $L^2$ (in $G_L$ replace $\nabla
$ by $\nabla ^m$). 

Then one can develop the $L^p$ theory of the semigroup, introducing the 
limits $p_\pm(L)$ for the $L^p$ boundedness of  $(e^{-tL})_{t>0}$
and the limits $q_\pm(L)$ for the $L^p$ boundedness of  
 $(\sqrt t\,  \nabla^m e^{-tL})_{t>0}$.
  The results are similar
with more technical burden in the arguments as one 
often has to control intermediate families $( t^{k/2m}\,  \nabla^k
e^{-tL})_{t>0}$ for $1\le k \le m-1$.~\footnote{\, see \cite{AT}, Chapter I,
\cite{AQ} and  \cite{Da}} One has that 
\begin{align*} p_-(L)&=q_-(L)
\\
p_+(L)&\ge (q_+(L))^{*m}
\end{align*}
 where $p^{*m}$ means $m$ times the operation $p\mapsto p^*$ . 

 By Sobolev embeddings plus perturbation results (such
as Lemma
\ref{lemmaw1p}) we have
$$ p_-(L) \begin{cases} =1, & {\rm if}\ n\le 2m \\
< \frac{2n}{n+ 2m},& {\rm if}\ n>2m.
\end{cases}
$$
Furthermore, this upper bound is
sharp for the class of {\bf all}  higher order
operators with $n>2m\ge 4$ : for any $n$ and $m  \ge 2$ with $2m<n$ and $p<\frac{2n}{n+ 2m}$, there 
exists  an operator (scalar or system)
$L$ of order
$2m$ such that
$p_-(L)>p$.~\footnote{\, This is in Davies \cite{Da1} based on examples of
Maz'ya and de Giorgi.} 
One has also 
$$
q_+(L) \begin{cases} =+\infty, & {\rm if}\ n=1 \\
>2,& {\rm if}\ n\ge 2.
\end{cases}
$$
 If $L^*$ is an example with $p_-(L^*)\sim \frac{2n}{n+ 2m}$ ($n>2m$), then 
$p_+(L)\sim \frac{2n}{n- 2m}$ and one has 
$p_+(L)\sim (q_+(L))^{*m}$ (Here, $\sim$ means within some arbitrary small $\ep$).
Thus the inequality $p_+(L)\ge (q_+(L))^{*m}$ is best possible.

 Next, the Riesz transform becomes 
$\nabla^m L^{-1/2}$ and 
$$\|\nabla^m L^{-1/2}f\|_p \lesssim \|f\|_p \quad\hbox{if and only if}\quad q_-(L)<p<q_+(L)$$ and one can
show the reverse inequalities (with $t\in \RR$)
$$
\| L^{1/2+it}f\|_p \lesssim \|\nabla^m f\|_p \quad\hbox{whenever}\quad 
\sup(1, (p_-(L))_{*m}) <p <p_+(L).$$
The bounds for
$p>2$ are merely obtained by duality from the Riesz transform bounds and $p_+(L)$ is best possible. A tool
to obtain the estimates for
$p<2$ is the extension of the Calder\'on-Zygmund decomposition to  Sobolev $W^{m,p}$
functions.~\footnote{\, See \cite{A} for the proof.} The lower limit $(p_-(L))_{*m}$ (if not $\le 1$,
which implies large dimensions), is best possible if $L$ is an operator for which
$p_+(L^*)\ge q_+(L^*)^{*m}$ is best possible. We have seen there exist such operators.
The connexion with Hodge theory is analogous to the second order case and the dichotomy
$p>2$ vs
$p<2$  appears again. 

The 
 critical numbers $\frac{2n}{n\pm 2}$ and $\frac{2n}{n+ 4}$ which appear in the
 $L^p$ theory of square roots for second order operators (namely Propositions \ref{propRTn=2} and
\ref{propRTnge3})  become
$\frac{2n}{n\pm 2m}$ and $\frac{2n}{n+ 4m}$.~\footnote{\, See \cite{A} where all
this is explicited. The case $\frac{2n}{n+ 2m}<p<2$ when $n> 2m$ is due to Blunck
\& Kunstmann \cite{BK2}. The other cases follow from the methods in \cite{AT}
although this is not explicited.} Hence, the range of
exponents $p$ for   the Riesz transform 
$L^p$ estimate is
\begin{enumerate}
\item[] $1<p<\infty$, when $n=1$, $m\ge 1$, 

\item[] $1<p< 2+\ep$, if  $1< n\le 2m$,

\item[] $\frac{2n}{n+ 2m}-\ep<p<2+\ep'$, when  $n> 2m$.
 \end{enumerate}
The discussion above show that these open ranges are best possible.
The range of exponents
$p$ for the reverse  inequality~\footnote{\, This is proved in \cite{A} for $p<2$ and
$n>4m$. The other cases were done earlier as a consequence of the methods in
\cite{AT} for $n\le 2m$ and in \cite{AHLMcT} for $2m\le n < 4m$.} is 
\begin{enumerate}
\item[]$1<p<\infty$ if
$n\le 2m$,

\item[]$1<p< \frac{2n}{n- 2m}+\ep$ if $2m< n \le 4m$,

 \item[] $\frac{2n}{n+ 4m}-\ep_1 <p<
\frac{2n}{n- 2m}+\ep$ if
$n>4m
$.
\end{enumerate}
Again,  these  open ranges  are best possible.~\footnote{\, In particular, this means that the methods used
here cannot improve the ranges of $p$ for second order operators (as improved methods would apply for
higher order as well) unless they use specific features of second order operators.} 

The bounded holomorphic functional calculus extends on $L^p$ 
for  $p_-(L)<p<p_+(L)$.~\footnote{\, This  is due to  Blunck \& Kunstmann
\cite{BK3}.}

The theory of square functions also generalize similarly. \footnote{\, 
For $g_L$, as for second order operators, one can also combine  results of Blunck \&
Kunstmann and  of Le Merdy; for
$G_L$ this is new.} For $g_L$ the range of $p$'s is $p_-(L)<p<p_+(L)$ and for
$G_L$, $q_-(L)<p<q_+(L)$. One could treat also variants of the non-tangential square functions.

\bigskip

If the strong G\aa rding inequality  is weakened by a term $+\kappa \|f\|_2^2$ in the right hand side (in both the operator and system cases), then one has to
replace 
$L+ \lambda$ for $\lambda \ge \kappa$ and the local theory applies. 
\bigskip

One can also add pertubation by lower order terms with bounded measurable coefficients without any harm
to the theory.

\section{Calder\'on-Zygmund decomposition for Sobolev functions}

Here we prove the Calder\'on-Zygmund decomposition for Sobolev functions in Lemma \ref{lemmaCZD}. The notation are those
of the statement. 

\begin{proof} If $p=\infty$, set $g=f$. Assume next that $p<\infty$. 
Let $\Omega= \{x \in \RR^n; M(|\nabla f|^p)(x) >\alpha^p\}$ where $M$ is the uncentered maximal operator over cubes of $\RR^n$.
If $\Omega$ is empty, then set $g=f$. Otherwise, the maximal theorem gives us
\begin{equation*} 
|\Omega| \le C\alpha^{-p} \int_{\RR^n} |\nabla f|^p.  \end{equation*} 
Let $F$ be the complement of $\Omega$. By the Lebesgue differentiation
theorem, $|\nabla f| \le \alpha$  almost everywhere on $F$. We also have,

\begin{lemma}\label{lemmalipschitzonF} One can redefine $f$ on a null set of $F$ so
that  for all $x \in F$, for all cube $Q$ centered at $x$,
\begin{equation}\label{eq19}
|f(x) - m_Qf| \le C\alpha \ell(Q) \end{equation}
where $\ell(Q)$ is the sidelength of $Q$ and for all $x,y \in F$, 
\begin{equation}\label{eq20}
|f(x) - f(y)| \le C\alpha |x-y|.
\end{equation}
The constant $C$ depends only on dimension and $p$.
\end{lemma}

Here $m_Ef$ denotes the mean of $f$ over $E$. It is well-defined if $E$ is a cube as $f$ is locally integrable.  Let us postpone the proof of 
this lemma and continue the argument. 

Let $(Q_i)$ be a Whitney decomposition of $\Omega$ by dyadic cubes. Hence, $\Omega$ is the disjoint union of the $Q_i$'s, 
the cubes $2Q_i$ are contained in  $\Omega$ and have the bounded overlap property, but the cubes $4Q_i$ intersect $F$.
As usual, $\lambda Q$ is the cube co-centered with $Q$ with sidelength $\lambda$ times that of $Q$. Hence \eqref{eqcsds4} and 
\eqref{eqcsds5} are satisfied by the cubes $2Q_i$. Let us now define the functions $b_i$. Let $(\calX_i)$ be a partition of unity on $\Omega$
associated to the covering $(Q_i)$ so that for each $i$, $\calX_i$ is a $C^1$ function supported in $2Q_i$ with $\|\calX_i\|_\infty +
\ell_i \|\nabla \calX_i\|_\infty \le c(n)$, $\ell_i$ being the sidelength of $Q_i$. Pick a point $x_i \in 4Q_i \cap F$. 
Set
$$
b_i = (f-f(x_i))\calX_i.
$$
It is clear that $b_i$ is supported in $2Q_i$. Let us estimate $\int_{2Q_i} |\nabla b_i|^p$. Introduce $\widetilde Q_i$ the cube centered at $x_i$ with 
sidelength $8\ell_i$. Then $2Q_i \subset \widetilde Q_i$. Set $c_i=m_{2Q_i}f$ and $\tilde c_i= m_{\widetilde Q_i}f$ and write
$$
b_i= (f-c_i)\calX_i + (c_i-\tilde c_i)\calX_i + (\tilde c_i -f(x_i))\calX_i.
$$
By \eqref{eq19} and \eqref{eqcsds5} for the cubes $2Q_i$,  $|\tilde c_i -f(x_i)| \le C\alpha \ell_i$, hence
$\int_{2Q_i} |\tilde c_i -f(x_i)|^p |\nabla \calX_i|^p \le C\alpha^p |2Q_i|$.
Next, using the $L^p$-Poincar\'e inequality and the fact that $\widetilde Q_i \cap F$ is not empty,
$$
|c_i-\tilde c_i| \le \frac{1}{|2Q_i|} \int_{\widetilde Q_i} |f - \tilde c_i| \le C \ell_i \left(\frac{1}{|\widetilde Q_i|} \int_{\widetilde Q_i} |\nabla
f|^p\right)^{1/p} \le C\alpha \ell_i.$$
Hence,  $\int_{2Q_i} |c_i-\tilde c_i|^p |\nabla \calX_i|^p \le C\alpha^p |2Q_i|$.
Lastly, since $\nabla \big((f-c_i)\calX_i\big) = \calX_i \nabla f + (f-c_i)\nabla\calX_i$, we have again by the $L^p$-Poincar\'e inequality 
and the fact that the average of $|\nabla f|^p$ on $2Q_i$ is controlled by $C\alpha^p$ that 
$$\int_{2Q_i} |\nabla \big((f-c_i)\calX_i\big)|^p  \le C\alpha^p |2Q_i|.$$ Thus \eqref{eqcsds3} is proved. 

Set $h(x)= \sum_i f(x_i) \nabla\calX_i(x)$. Note that this sum is locally finite and $h(x)=0$ for $x \in F$. Note also that 
$\sum_i \calX_i(x)$ is 1 on $\Omega$ and 0 on $F$. Since it is also locally finite we have $\sum_i  \nabla\calX_i(x)=0$ for $x \in \Omega$.  We claim that $h(x) \le
C\alpha$. Indeed, fix $x \in \Omega$. Let
$Q_j$ be the Whitney cube containing
$x$ and let $I_x$ be the set of indices $i$ such that $x \in 2Q_i$. We know that $\sharp I_x \le N$.  Also for $i \in I_x$ we have that 
$C^{-1}\ell_i \le  \ell_j \le C\ell_i$ and $|x_i-x_j| \le C\ell_j$ where the constant $C$ depends only on dimension 
(see \cite{St1}).  We have
$$
|h(x)| = \left|\sum_{i \in I_x} (f(x_i) -f(x_j)) \nabla\calX_i(x)\right| \le C \sum_{i \in I_x} |f(x_i) -f(x_j)| \ell_i^{-1} \le C N \alpha,$$
by the previous observations. 

It remains to obtain \eqref{eqcsds1} and \eqref{eqcsds2}. We easily have  using $\sum_i  \nabla\calX_i(x)=0$ for $x \in \Omega$, that 
$$
\nabla f = (\nabla f) {\bf 1}_F + h + \sum_i \nabla b_i, \quad \text{a.e.}.
$$
Now $ \sum_i b_i$ is a well-defined distribution on $\RR^n$. Indeed, for a test function $u$, using the properties of the Whitney cubes,
$$\sum_i \int |b_i u|\le C \int \left(\sum_i |b_i(x)|\ell_i^{-1}\right) |u(x)| d(x,F)\, dx$$ and the last sum converges in
$L^p$
 as a consequence of \eqref{eqcsds4} and 

\begin{lemma}\label{lemmasumbi} Set $p^*= \frac{np}{n-p}$ if $p<n$ and $p^*=\infty$ otherwise, then for all real numbers $r$ with $p\le r\le p^*$,
\begin{equation}\label{eq20bis}
\|\sum_i |b_i|\ell_i^{-1}\|_r^r \le C \alpha^r \sum_i |Q_i|.
\end{equation}
\end{lemma}

Admit this lemma and set $g= f- \sum_i b_i$. Then 
$\nabla g = (\nabla f ){\bf 1}_F + h$ in the sense of distributions and, hence, $\nabla g$ is a bounded function with $\|\nabla g\|_\infty \le
C\alpha$.
\end{proof}

\paragraph{Proof of Lemma \ref{lemmasumbi}:} By \eqref{eqcsds5} and the Poincar\'e-Sobolev inequality:
$$
\|\sum_i |b_i|\ell_i^{-1}\|_r^r \le N \sum_i \| |b_i|\ell_i^{-1}\|_r^r \le N C \sum_i \ell_i^{r\theta}\|\nabla b_i\|^r_p$$
where $\theta= \frac {n} r - \frac n p$. By  \eqref{eqcsds3}, $\ell_i^{r\theta}\|\nabla b_i\|^r_p \le \alpha^r \ell_i^{nr/p}$, hence
$$
\|\sum_i |b_i|\ell_i^{-1}\|_r^r  \le CN \alpha^r \sum_i \ell_i^n.$$

\paragraph{Proof of  Lemma \ref{lemmalipschitzonF}:}

Let $x$ be a point in $F$. Fix such  cube $Q$ with center $x$ and let $Q_k$ be co-centered cubes with $\ell(Q_k) =2^k\ell(Q)$ for $k$ a negative  integer. Then, by
Poincar\'e's inequality
\begin{align*}|m_{Q_{k+1}}f - m_{Q_{k}}f| &
\le 2^n |m_{Q_{k+1}}(f - m_{Q_{k+1}}f)| 
\\
&\le C 2^n \ell(Q_k) (m_{Q_{k+1}}|\nabla f|^p)^{1/p} 
\\
&\le
C2^k\ell(Q)\alpha
\end{align*}
 since $Q_{k+1}$ contains $x \in F$. It easily follows that $m_Qf$ has a limit as $|Q|$ tends to 0. If, moreover, $x$ is
 in the Lebesgue set of $f$, then this limit is equal to $f(x)$. Redefine $f$ on the complement of the Lebesgue set in $F$ so 
that
$m_Qf$ tends to
$f(x)$ with $Q$ centered at
$x$ with
$|Q|
\to 0$.   Moreover,  summing over $k$ the previous inequality gives us \eqref{eq19}. To see \eqref{eq20}, let $Q_x$ be the cube centered at $x$ with 
sidelength $2|x-y|$ and $Q_y$ be the cube centered at $y$ with sidelength $8|x-y|$. It is easy to see that $Q_x \subset Q_y$. As before, one can see
that $|m_{Q_{x}}f - m_{Q_{y}}f| \le C\alpha |x-y|$. Hence by the triangle inequality and \eqref{eq19}, one obtains \eqref{eq20} readily.

\end{document}